\documentclass[10pt,a4paper]{article}
\usepackage{amsmath,amssymb,amsthm,graphics,graphicx,epic,eepic,multicol,ascmac}
\usepackage[mathscr]{euscript}
\usepackage[all]{xy}
\usepackage{fancybox}
\usepackage{url}
\usepackage{ytableau}

\numberwithin{equation}{section}
\theoremstyle{theorem}
\newtheorem{thm}{Theorem}[section]
\newtheorem{prop}[thm]{Proposition}
\newtheorem{lem}[thm]{Lemma}
\newtheorem{rem}[thm]{Remark}
\newtheorem{cor}[thm]{Corollary}

\theoremstyle{definition}
\newtheorem{defn}[thm]{Definition}
\newtheorem{ex}[thm]{Example}

\def\al{\alpha}

\def\wht(#1){\widehat{\ #1\ }}

\newcommand{\frg}{\mathfrak g}

\newcommand{\frt}{\mathfrak t}

\newcommand{\bbQ}{\mathbb Q}

\newcommand{\bbZ}{\mathbb Z}

\newcommand{\ch}{\mathrm{ch}}

\newcommand{\lbr}{\begin{bmatrix}}
\newcommand{\rbr}{\end{bmatrix}}

\newcommand{\cd}{commutative diagram }

\def\ge{\frg}

\def\al{\alpha}

\def\beneme{\begin{enumerate}}
\def\beq{\begin{equation}}
\def\beqn{\begin{eqnarray}}
\def\beqnn{\begin{eqnarray*}}

\def\bfii0{{\bf i_0}}

\def\bbra#1,#2,#3{\left\{\begin{array}{c}\hspace{-5pt}
#1;#2\\ \hspace{-5pt}#3\end{array}\hspace{-5pt}\right\}}
\def\cd{\cdots}
\def\ci(#1,#2){c_{#1}^{(#2)}}
\def\Ci(#1,#2){C_{#1}^{(#2)}}
\def\mpp(#1,#2,#3){#1^{(#2)}_{#3}}
\def\bCi(#1,#2){\ovl C_{#1}^{(#2)}}
\def\ch(#1,#2){c_{#2,#1}^{-h_{#1}}}
\def\cc(#1,#2){c_{#2,#1}}

\def\di(#1,#2){D_{#1}^{(#2)}}
\def\dbi(#1,#2){\ovl D_{#1}^{(#2)}}

\def\eit{\tilde{e}_i}
\def\eneme{\end{enumerate}}

\def\eeq{\end{equation}}
\def\eeqn{\end{eqnarray}}
\def\eeqnn{\end{eqnarray*}}
\def\fit{\tilde{f}_i}

\def\gau#1,#2{\left[\begin{array}{c}\hspace{-5pt}#1\\
\hspace{-5pt}#2\end{array}\hspace{-5pt}\right]}

\def\ify{\infty}
\def\io{\iota}
\def\ji(#1,#2){j_{#1}^{(#2)}}
\def\kp{k^{(+)}}
\def\km{k^{(-)}}

\def\lan{\langle}
\def\lar{\longrightarrow}
\def\lm{\lambda}
\def\Lm{\Lambda}

\def\nd{\noindent}

\def\ovl{\overline}

\def\qq{\qquad}

\def\qed{\hfill\framebox[2mm]{}}
\def\QQ{\mathbb Q}

\def\ran{\rangle}

\def\tt{\frt}
\def\TY(#1,#2,#3){#1^{(#2)}_{#3}}

\def\uq{U_q(\ge)}

\def\uqm{U^-_q(\ge)}

\def\vp{\varphi}

\def\xxi(#1,#2,#3){\displaystyle {}^{#1}\Xi^{(#2)}_{#3}}
\def\xsi(#1,#2,#3){\displaystyle {}^{#1}\Sigma^{(#2)}_{#3}}
\def\xE(#1,#2,#3){\displaystyle {}^{#1}E_{#2}[#3]}
\def\xF(#1,#2){\displaystyle {}^{#1}F_{#2}}
\def\xx(#1,#2){\displaystyle {}^{#1}\Xi_{#2}}
\def\W1{W(\varpi_1)}

\def\ZZ{\mathbb Z}

\def\m@th{\mathsurround=0pt}
\def\fsquare(#1,#2){
\hbox{\vrule$\hskip-0.4pt\vcenter to #1{\normalbaselines\m@th
\hrule\vfil\hbox to #1{\hfill$\scriptstyle #2$\hfill}\vfil\hrule}$\hskip-0.4pt
\vrule}}

\newcommand{\ba}{\begin{array}}
\newcommand{\ea}{\end{array}}

\newcommand{\eq}{\begin{eqnarray}}
\newcommand{\eneq}{\end{eqnarray}}

\title{\textbf{\large{Adapted Sequence for Polyhedral Realization of Crystal Bases}}}
\author{\normalsize{YUKI KANAKUBO\thanks{Division of Mathematics, 
Sophia University, Kioicho 7-1, Chiyoda-ku, Tokyo 102-8554,
Japan: {j\_chi\_sen\_you\_ky@eagle.sophia.ac.jp. }}
\ and\ 
TOSHIKI NAKASHIMA\thanks{Division of Mathematics, 
Sophia University, Kioicho 7-1, Chiyoda-ku, Tokyo 102-8554,
Japan: {toshiki@sophia.ac.jp}.}
}}
\date{}

\begin{document}

\maketitle
\vspace{-10pt}

\begin{abstract}
The polyhedral realization of crystal base 
has been introduced by A.Zelevinsky and the second author(\cite{NZ}), 
which describe 
the crystal base $B(\infty)$ as a polyhedral
convex cone in the infinite $\bbZ$-lattice $\mathbb{Z}^{\infty}$.
To construct the polyhedral realization, we need to fix
an infinite sequence $\iota$ from the indices of the simple roots.
According to this $\io$, one has certain set of linear functions
defining a polyhedral convex cone and 
under the `positivity condition' on $\iota$, 
it has been shown that the polyhedral convex cone 
is isomorphic to the crystal base $B(\infty)$. 
To confirm the positivity condition for a given $\iota$, 
we need to obtain the whole feature of the set of linear functions, 
which requires, in general, a bunch of explicit calculations. 
In this article, 
we introduce the notion of the {\it adapted sequence} 
and show that if $\iota$ is an adapted sequence then
the positivity condition holds for classical Lie algebras.
Furthermore, we reveal the explicit
forms of the polyhedral realizations associated with arbitrary 
adapted sequences $\iota$ in terms of column tableaux.
\end{abstract}

\section{Introduction}

The crystal $B(\infty)$ is, roughly speaking, a basis of the subalgebra $U^-_q(\mathfrak{g})$ of
the quantum group $U_q(\mathfrak{g})$ at $q=0$, where $\ge$ is
a symmetrizable Kac-Moody Lie algebra with an index set $I=\{1,2,\cdots,n\}$.
It had been invented by 
G.Lusztig for A, D, E cases \cite{L} and 
by M.Kashiwara for Kac-Moody cases \cite{K1}.
Since then the theory of crystal base has influenced many areas of mathematics and physics, e.g.,
algebraic combinatorics, modular representation theory, automorphic forms, statistical mechanics, cellular automaton, etc. 
In order to apply the theory of crystal bases to these areas, it is required to realize the crystal bases in suitable forms,
like as tableaux realizations, path realizations, geometric realizations, etc.
In \cite{NZ}, the polyhedral realization has been introduced as an image of  `Kashiwara embedding'
$\Psi_{\iota}:B(\infty)\hookrightarrow \mathbb{Z}^{\infty}$, where $\iota$ is an infinite sequence of entries in $I$
and $\mathbb{Z}^{\infty}$ is an infinite $\mathbb{Z}$-lattice with certain crystal structure. More precisely, we will associate some set of
linear functions $\Xi_{\iota}\subset (\mathbb{Q}^{\infty})^*$ with the sequence $\iota$ and define the subset
$\Sigma_{\iota}\subset \mathbb{Z}^{\infty}$ by
\[ \Sigma_{\io}
=\{\textbf{x}\in\mathbb{Z}^{\infty} | \varphi(\textbf{x})\geq0, \ \ \forall\varphi\in\Xi_{\iota} \}.\]
Then under some condition on $\iota$ called `positivity condition', we find that 
${\rm Im}(\Psi_{\iota})=\Sigma_{\io}\cong B(\infty)$, which implies that the crystal $B(\infty)$ is realized as a polyhedral
convex cone in $\mathbb{Z}^{\infty}$. To confirm the positivity condition for a given $\iota$, it is necessary
to obtain the whole linear functions of $\Xi_{\iota}$, which requires a lot of calculations.
So far, in \cite{H1,H2,NZ}  
it has been shown that the specific sequence $\iota=(\cdots,2,1,n,\cdots,2,1,n\cdots,2,1)$ satisfies the positivity condition
for all simple Lie algebras and almost all affine Lie algebras.

The aim of this article is to give a sufficient condition of the positivity condition for a
sequence $\iota$ and
construct the explicit forms of $\Xi_{\iota}$ in terms of column tableaux for sequences which satisfy the sufficient condition
in the case $\mathfrak{g}$ is a classical Lie algebra. Here, note that the column tableaux do not correspond to
crystal base elements themselves but represent linear functions in $\Xi_{\iota}$.
In Definition \ref{adapt}, we introduced the notion of sequences {\it adapted} to the Cartan matrix of $\mathfrak{g}$.
A similar notion `adapted to a quiver' has been defined for reduced words of Weyl group elements of type A,D,E in \cite{HL, L,R}.
If $\mathfrak{g}$ is of type A, B, C or D, then it is the sufficient condition for the positivity condition (Theorem \ref{thm1}).
To give an explicit form of $\Xi_{\iota}$ for a sequence $\iota$ adapted to the Cartan matrix of type X (X$=$A, B, C or D),
we define a set ${\rm Tab}_{{\rm X},\iota}$ of column tableaux, which stand for linear functions in Definition \ref{box-def}, \ref{tab-def}.
In Theorem \ref{thm2}, we show that $\Xi_{\iota}={\rm Tab}_{{\rm X},\iota}$.
Using this description of $\Xi_{\iota}$, it is proved that the polyhedral realization $\Sigma_{\iota}$ is expressed 
by finitely many inequalities of finitely many variables (Corollary \ref{cor1}).
For example, in the case $\mathfrak{g}$ is of type ${\rm A}_2$ and $\iota=(\cdots,2,1,2,1,2,1)$, we obtain
\begin{equation}\label{intro-eq}
\Xi_{\iota}={\rm Tab}_{{\rm A},\iota}
=\{\begin{ytableau}
 i 
\end{ytableau}^{\rm A}_{s} | 1\leq i\leq 3,\ s\in\mathbb{Z}_{\geq1}\}\cup
\{\begin{ytableau}
 i \\
j
\end{ytableau}^{\rm A}_{s} | 1\leq i<j\leq 3,\ s\in\mathbb{Z}_{\geq1}\}.
\end{equation}
We define $x_k\in (\QQ^{\ify})^*$ as the linear function as $x_k((\cd,a_k,\cd,a_2,a_1))=a_k$
for $k\in\mathbb{Z}_{\geq1}$ and rewrite them as double indexed variables
$x_{2l-1}=x_{l,1}$, $x_{2l}=x_{l,2}$ for $l\in\mathbb{Z}_{\geq1}$. Then 
the above boxes and tableaux are as follows:
\[
\begin{ytableau}
 i 
\end{ytableau}^{\rm A}_{s}=x_{s,i}-x_{s+1,i-1},
\]
\[
\begin{ytableau}
 i \\
j
\end{ytableau}^{\rm A}_{s}=x_{s+1,i}-x_{s+2,i-1} +x_{s,j}-x_{s+1,j-1},
\]
where $x_{m,0}=x_{m,3}=0$ for $m\in\mathbb{Z}$.
Thus, we get
\begin{eqnarray*}
 \hspace{-7mm} {\rm Im}(\Psi_{\iota})
&=&\left\{\textbf{x}\in\mathbb{Z}^{\infty} | x_{s,i}-x_{s+1,i-1}\geq 0\ (1\leq i\leq 3,\ s\in\mathbb{Z}_{\geq1}),\right.\\ 
& &\left. x_{s+1,i}-x_{s+2,i-1} +x_{s,j}-x_{s+1,j-1}\geq0\ \ (1\leq i<j\leq 3,\ s\in\mathbb{Z}_{\geq1})  \right\}.\\
\end{eqnarray*}
Simplifying the inequalities, we get
\[
 {\rm Im}(\Psi_{\iota})=
\{\textbf{x}\in\mathbb{Z}^{\infty} |  
x_{1,2}\geq x_{2,1}\geq0,\ x_{1,1}\geq0,\ x_{m+1,2}=x_{m+2,1}=0,\ \forall m\in\mathbb{Z}_{\geq1}\}.\]

The organization of this article is as follows. In Sect.2, we review on 
the theory of crystals and polyhedral realizations. In Sect.3, we introduce the
column tableaux and give our main results. In Sect.4, we prove
the sets ${\rm Tab}_{{\rm X},\iota}$ (X$=$A, B, C or D) 
of the tableaux are closed under the actions of operators $S_k$ ($k\in\mathbb{Z}_{\geq1}$). 
In Sect.5, we will prove our main theorem.

\vspace{2mm}

\nd
{\bf Acknowledgements} 
Y.K. is supported by JSPS KAKENHI Grant Number JP17H07103, 
and T.N. is supported in part by JSPS KAKENHI Grant Number JP15K04794.

\section{Crystal and its polyhedral realization}

Let us recall the definition of {\it crystals} \cite{K3}. 

\subsection{Notations}

We list the notations used in this paper.
Let $\ge$ be
a symmetrizable Kac-Moody Lie algebra over $\bbQ$
with a Cartan subalgebra
$\tt$, a weight lattice $P \subset \tt^*$, the set of simple roots
$\{\al_i: i\in I\} \subset \tt^*$,
and the set of coroots $\{h_i: i\in I\} \subset \tt$,
where $I=\{1,2,\cdots,n\}$ is a finite index set.
Let $\lan h,\lm\ran=\lm(h)$ be the pairing between $\tt$ and $\tt^*$,
and $(\al, \beta)$ be an inner product on
$\tt^*$ such that $(\al_i,\al_i)\in 2\bbZ_{\geq 0}$ and
$\lan h_i,\lm\ran={{2(\al_i,\lm)}\over{(\al_i,\al_i)}}$
for $\lm\in\tt^*$ and $A:=(\lan h_i,\al_j\ran)_{i,j}$ be the associated generalized symmetrizable Cartan matrix.
Let $P^*=\{h\in \tt: \lan h,P\ran\subset\ZZ\}$ and
$P_+:=\{\lm\in P:\lan h_i,\lm\ran\in\ZZ_{\geq 0}\}$.
We call an element in $P_+$ a {\it dominant integral weight}.
The quantum algebra $\uq$
is an associative
$\QQ(q)$-algebra generated by the $e_i$, $f_i \,\, (i\in I)$,
and $q^h \,\, (h\in P^*)$
satisfying the usual relations.
The algebra $\uqm$ is the subalgebra of $\uq$ generated 
by the $f_i$ $(i\in I)$.

For the irreducible highest weight module of $\uq$
with the highest weight $\lm\in P_+$, we denote it by $V(\lm)$
and its {\it crystal base} we denote $(L(\lm),B(\lm))$.
Similarly, for the crystal base of the algebra $\uqm$ we denote 
$(L(\ify),B(\ify))$ (see \cite{K0,K1}).
For positive integers $l$ and $m$ with $l\leq m$, we set $[l,m]:=\{l,l+1,\cdots,m-1,m\}$.

\subsection{Crystals}

By the terminology {\it crystal } we mean some combinatorial object 
obtained by abstracting the properties of crystal bases:

\begin{defn}
A {\it crystal} is a set $\mathcal{B}$ together with the maps
${\rm wt}:\mathcal{B}\rightarrow P$,
$\varepsilon_i,\varphi_i:\mathcal{B}\rightarrow \mathbb{Z}\cup \{-\infty\}$
and $\tilde{e}_i$,$\tilde{f}_i:\mathcal{B}\rightarrow \mathcal{B}\cup\{0\}$
($i\in I$) satisfying the following: For $b,b'\in\mathcal{B}$, $i,j\in I$,
\begin{enumerate}
\item[$(1)$] $\varphi_i(b)=\varepsilon_i(b)+\langle{\rm wt}(b),h_i\rangle$,
\item[$(2)$] ${\rm wt}(\tilde{e}_ib)={\rm wt}(b)+\alpha_i$ if $\tilde{e}_i(b)\in\mathcal{B}$,
\quad ${\rm wt}(\tilde{f}_ib)={\rm wt}(b)-\alpha_i$ if $\tilde{f}_i(b)\in\mathcal{B}$,
\item[$(3)$] $\varepsilon_i(\tilde{e}_i(b))=\varepsilon_i(b)-1,\ \ 
\varphi_i(\tilde{e}_i(b))=\varphi_i(b)+1$\ if $\tilde{e}_i(b)\in\mathcal{B}$, 
\item[$(4)$] $\varepsilon_i(\tilde{f}_i(b))=\varepsilon_i(b)+1,\ \ 
\varphi_i(\tilde{f}_i(b))=\varphi_i(b)-1$\ if $\tilde{f}_i(b)\in\mathcal{B}$, 
\item[$(5)$] $\tilde{f}_i(b)=b'$ if and only if $b=\tilde{e}_i(b')$,
\item[$(6)$] if $\varphi_i(b)=-\infty$ then $\tilde{e}_i(b)=\tilde{f}_i(b)=0$.
\end{enumerate}
We call $\tilde{e}_i$,$\tilde{f}_i$ {\it Kashiwara operators}.
\end{defn}

\begin{defn}
A {\it strict morphism} $\psi:\mathcal{B}_1\rightarrow\mathcal{B}_2$
of crystals $\mathcal{B}_1$, $\mathcal{B}_2$ is a map 
$\mathcal{B}_1\bigsqcup\{0\}\rightarrow\mathcal{B}_2\bigsqcup\{0\}$
satisfying the following conditions:
 $\psi(0)=0$, 
${\rm wt}(\psi(b))={\rm wt}(b)$,
 $\varepsilon_i(\psi(b))=\varepsilon_i(b)$,
 $\varphi_i(\psi(b))=\varphi_i(b)$,
 if $b\in \mathcal{B}_1$, $\psi(b)\in \mathcal{B}_2,$ $i\in I$,
 and
$\psi: \mathcal{B}_1\sqcup\{0\} \lar \mathcal{B}_2\sqcup\{0\}$
commutes with all $\eit$ and $\fit$, where $\eit(0)=\fit(0)=0$.
An injective strict morphism is said to be {\it embedding} of crystals.
\end{defn}

\subsection{Polyhedral Realization of $B(\ify)$}
\label{poly-uqm}
Let us recall the results in \cite{NZ}.

Consider the infinite $\bbZ$-lattice
\begin{equation}
\ZZ^{\ify}
:=\{(\cd,x_k,\cd,x_2,x_1): x_k\in\ZZ
\,\,{\rm and}\,\,x_k=0\,\,{\rm for}\,\,k\gg 0\};
\label{uni-cone}
\end{equation}
we will denote by $\ZZ^{\ify}_{\geq 0} \subset \ZZ^{\ify}$
the subsemigroup of nonnegative sequences.
For the rest of this section, we fix an infinite sequence of indices
$\io=\cd,i_k,\cd,i_2,i_1$ from $I$ such that
\begin{equation}
{\hbox{
$i_k\ne i_{k+1}$ and $\sharp\{k: i_k=i\}=\ify$ for any $i\in I$.}}
\label{seq-con}
\end{equation}

Given $\io$, we can define a crystal structure
on $\ZZ^{\ify}$ and denote it by $\ZZ^{\ify}_{\io}$ 
(\cite[2.4]{NZ}).

\begin{prop}[\cite{K3}, See also \cite{NZ}]
\label{emb}
There is a unique strict embedding of crystals
$($called Kashiwara embedding$)$
\begin{equation}
\Psi_{\io}:B(\ify)\hookrightarrow \ZZ^{\ify}_{\geq 0}
\subset \ZZ^{\ify}_{\io},
\label{psi}
\end{equation}
such that
$\Psi_{\io} (u_{\ify}) = (\cd,0,\cd,0,0)$, where 
$u_{\ify}\in B(\ify)$ is the vector corresponding to $1\in \uqm$.
\end{prop}

\begin{defn}
The image ${\rm Im} \Psi_{\iota} (\cong B(\infty))$ is called a \it{polyhedral realization} of $B(\infty)$.
\end{defn}

Consider the infinite dimensional vector space
$$
\QQ^{\ify}:=\{\textbf{a}=
(\cd,a_k,\cd,a_2,a_1): a_k \in \QQ\,\,{\rm and }\,\,
a_k = 0\,\,{\rm for}\,\, k \gg 0\},
$$
and its dual space $(\QQ^{\ify})^*:={\rm Hom}(\QQ^{\ify},\QQ)$.
Let $x_k\in (\QQ^{\ify})^*$ be the linear function defined as $x_k((\cd,a_k,\cd,a_2,a_1)):=a_k$
for $k\in\mathbb{Z}_{\geq1}$.
We will also write a linear form $\vp \in (\QQ^{\ify})^*$ as
$\vp=\sum_{k \geq 1} \vp_k x_k$ ($\vp_j\in \QQ$).

For the fixed infinite sequence
$\io=(i_k)$ and $k\geq1$ we set $\kp:={\rm min}\{l:l>k\,\,{\rm and }\,\,i_k=i_l\}$ and
$\km:={\rm max}\{l:l<k\,\,{\rm and }\,\,i_k=i_l\}$ if it exists,
or $\km=0$  otherwise.
We set $\beta_0=0$ and
\begin{equation}
\beta_k:=x_k+\sum_{k<j<\kp}\lan h_{i_k},\al_{i_j}\ran x_j+x_{\kp}\in (\QQ^{\ify})^*
\qq(k\geq1).
\label{betak}
\end{equation}
We define the piecewise-linear operator 
$S_k=S_{k,\io}$ on $(\QQ^{\ify})^*$ by
\begin{equation}
S_k(\vp):=
\begin{cases}
\vp-\vp_k\beta_k & {\mbox{ if }}\vp_k>0,\\
 \vp-\vp_k\beta_{\km} & {\mbox{ if }}\vp_k\leq 0.
\end{cases}
\label{Sk}
\end{equation}
Here we set
\begin{eqnarray}
\Xi_{\io} &:=  &\{S_{j_l}\cd S_{j_2}S_{j_1}x_{j_0}\,|\,
l\geq0,j_0,j_1,\cd,j_l\geq1\},
\label{Xi_io}\\
\Sigma_{\io} & := &
\{\textbf{x}\in \ZZ^{\ify}\subset \QQ^{\ify}\,|\,\vp(\textbf{x})\geq0\,\,{\rm for}\,\,
{\rm any}\,\,\vp\in \Xi_{\io}\}.
\end{eqnarray}
We impose on $\io$ the following {\it positivity condition}:
\begin{equation}
{\hbox{if $\km=0$ then $\vp_k\geq0$ for any 
$\vp=\sum_k\vp_kx_k\in \Xi_{\io}$}}.
\label{posi}
\end{equation}
\begin{thm}[\cite{NZ}]\label{polyhthm}
Let $\io$ be a sequence of indices satisfying $(\ref{seq-con})$ 
and $(\ref{posi})$. Then we have 
${\rm Im}(\Psi_{\io})(\cong B(\ify))=\Sigma_{\io}$.
\end{thm}

\begin{ex}

Let $\mathfrak{g}$ be of type ${\rm A}_2$ and $\iota=(\cdots,2,1,2,1,2,1)$.
It follows
\[
1^-=2^-=0,\ k^->0\ (k>2).\]
We rewrite a vector $(\cdots,x_6,x_5,x_4,x_3,x_2,x_1)$ as 
\[
(\cdots,x_{3,2},x_{3,1},x_{2,2},x_{2,1},x_{1,2},x_{1,1}), \]
that is, $x_{2l-1}=x_{l,1}$, $x_{2l}=x_{l,2}$ for $l\in\mathbb{Z}_{\geq1}$.
Recall that positivity condition means that the coefficients of $x_1=x_{1,1}$ and $x_2=x_{1,2}$ 
in each $\varphi\in\Xi_{\iota}$ are non-negative.
Similarly, we rewrite $S_{2l-1}=S_{l,1}$, $S_{2l}=S_{l,2}$.
For $k\in\mathbb{Z}_{\geq1}$, the action of the operators are the following:
\[
x_{k,1}\overset{S_{k,1}}{\rightleftarrows}x_{k,2}-x_{k+1,1}\overset{S_{k,2}}{\rightleftarrows} -x_{k+1,2},
\]
\vspace{-8mm}
\[\hspace{-6mm}
 {\scriptstyle S_{k+1,1}}\quad \quad \quad \quad \quad {\scriptstyle S_{k+1,2}}
\]

\[
x_{k,2}\overset{S_{k,2}}{\rightleftarrows}x_{k+1,1}-x_{k+1,2}\overset{S_{k+1,1}}{\rightleftarrows} -x_{k+2,1},
\]
\vspace{-8mm}
\[\hspace{-7mm}
 {\scriptstyle S_{k+1,2}}\quad \quad \quad \quad \quad \quad {\scriptstyle S_{k+2,1}}
\]
and other actions are trivial. Thus we obtain
\[
\Xi_{\iota}=\{x_{k,1},\ x_{k,2}-x_{k+1,1},\ -x_{k+1,2},\ x_{k,2},\ x_{k+1,1}-x_{k+1,2}, -x_{k+2,1} | k\geq1\}.
\]
Note that
the coefficients of $x_{1,1}$ and $x_{1,2}$ in each $\varphi\in\Xi_{\iota}$ are non-negative.
Therefore $\iota$ satisfies the positivity condition and
\begin{eqnarray*}
& &{\rm Im}(\Psi_{\iota})=\Sigma_{\io}
=\{\textbf{x}\in\mathbb{Z}^{\infty}_{\iota} | \varphi(\textbf{x})\geq0, \ \ \forall\varphi\in\Xi_{\iota} \}.
\end{eqnarray*}
For $\textbf{x}=(\cdots,x_{3,2},x_{3,1},x_{2,2},x_{2,1},x_{1,2},x_{1,1})\in {\rm Im}(\Psi_{\iota})$,
combining inequalities $x_{k,1}\geq0$, $-x_{k+2,1}\geq0$ $(k\geq1)$, we obtain
$x_{k,1}=0$ $(k\geq3)$. Similarly, by $x_{k,2}\geq0$, $-x_{k+1,2}\geq0$ $(k\geq1)$,
we get $x_{k,2}=0$ $(k\geq2)$. Hence, we obtain
\[
 {\rm Im}(\Psi_{\iota})=\Sigma_{\io}=\{\textbf{x}\in \mathbb{Z}^{\infty}_{\iota} | x_{k+1,1}=x_{k,2}=0\ {\rm for}\ k\in\mathbb{Z}_{\geq2},\ 
 \ x_{1,2}\geq x_{2,1}\geq 0,\ x_{1,1}\geq0 \}.
\]

\end{ex}

\begin{ex}\cite{N99}\label{counterofpos}
Let $\mathfrak{g}$ be of type ${\rm A}_3$ and $\iota=(\cdots,2,1,2,3,2,1)$.
We obtain
\[
x_{1} \overset{S_1}{\rightarrow} -x_5+x_4+x_2 \overset{S_2}{\rightarrow} -x_5 +x_3  
\overset{S_5}{\rightarrow} -x_4+x_3-x_2+x_1.
\]
Thus, $-x_4+x_3-x_2+x_1\in\Xi_{\iota}$ and $2^-=0$.
Therefore $\iota$ does not satisfy the positivity condition.
\end{ex}


\subsection{Infinite sequences adapted to $A$}

\begin{defn}\label{adapt}
Let $A=(a_{i,j})$ be the generalized symmetrizable Cartan matrix of $\mathfrak{g}$ and $\io$
a sequence of indices satisfying $(\ref{seq-con})$.
If $\iota$ satisfies the following condition, we say $\iota$ is {\it adapted} to $A$ : 
For $i,j\in I$ with $i\neq j$ and $a_{i,j}\neq0$, the subsequence of $\iota$ consisting of all $i$, $j$ is
\[
(\cdots,i,j,i,j,i,j,i,j)\quad {\rm or}\quad (\cdots,j,i,j,i,j,i,j,i).
\]
If the Cartan matrix is fixed then the sequence $\iota$ is shortly said to be {\it adapted}.
\end{defn}

\begin{ex}
Let us consider the case
 $\mathfrak{g}$ is of type ${\rm A}_3$, $\iota=(\cdots,2,1,3,2,1,3,2,1,3)$.
\begin{itemize}
\item
The subsequence consisting of $1$, $2$ : $(\cdots,2,1,2,1,2,1)$.
\item
The subsequence consisting of $2$, $3$ : $(\cdots,2,3,2,3,2,3)$.
\item
Since $a_{1,3}=0$ we do not need consider the pair $1$, $3$.
\end{itemize}

Thus $\iota$ is an adapted sequence.
\end{ex}

\begin{ex}
Let us consider the case
 $\mathfrak{g}$ is  of type ${\rm A}_3$, $\iota=(\cdots,2,1,2,3,2,1)$,
which is the same setting as in Example \ref{counterofpos}.
The subsequence consisting of $1$, $2$ is $(\cdots,2,1,2,2,1)$.
Thus $\iota$ is not an adapted sequence.
\end{ex}

\section{Tableaux descriptions of Polyhedral realizations}

In this section, we take $\mathfrak{g}$ as a finite dimensional simple Lie algebra of type ${\rm A}_n$,
${\rm B}_n$, ${\rm C}_n$ or ${\rm D}_n$. In the rest of article, we follow Kac's notation \cite{Kac}.

In what follows, we suppose $\iota=(\cdots,i_3,i_2,i_1)$ satisfies
$(\ref{seq-con})$ and is adapted to the Cartan matrix $A$ of $\mathfrak{g}$.
Let $(p_{i,j})_{i\neq j,\ a_{i,j}\neq0}$ be the set of integers such that
\begin{equation}\label{pij}
p_{i,j}=\begin{cases}
1 & {\rm if}\ {\rm the\ subsequence\ of\ }\iota{\rm\ consisting\ of}\ i,j\ {\rm is}\ (\cdots,j,i,j,i,j,i),\\
0 & {\rm if}\ {\rm the\ subsequence\ of\ }\iota{\rm\ consisting\ of}\ i,j\ {\rm is}\ (\cdots,i,j,i,j,i,j).
\end{cases}
\end{equation}
For $k$ $(2\leq k\leq n)$, we set
\[
P(k):=\begin{cases}
p_{2,1}+p_{3,2}+\cdots+p_{n-2,n-3}+p_{n,n-2} & {\rm if}\ k=n\ {\rm and}\ \mathfrak{g}\ {\rm is\ of}\ {\rm type\ D}_n, \\
p_{2,1}+p_{3,2}+p_{4,3}+\cdots+p_{k,k-1} & {\rm if}\ {\rm otherwise}
\end{cases}
\]
and $P(0)=P(1)=P(n+1)=0$. Since each $p_{i,j}$ is in $\{0,1\}$,
it holds for $k$, $l\in I$ such that $k\geq l$.
\begin{equation}\label{pineq}
P(k)\geq P(l),
\end{equation}
\begin{equation}\label{pineq2}
(k-l)+P(l)\geq P(k),
\end{equation}
except for the case $\mathfrak{g}$ is of type ${\rm D}_n$, $k=n$ and $l=n-1$. 

For $k\in\mathbb{Z}_{\geq1}$, we rewrite $x_k$, $\beta_k$ and $S_k$ in \ref{poly-uqm} as
\begin{equation}\label{rewrite}
x_k= x_{s,j},\quad S_k=S_{s,j},\quad \beta_k= \beta_{s,j} 
\end{equation}
if $i_k=j$ and $j$ is appearing $s$ times in $i_k$, $i_{k-1}$, $\cdots,i_1$.
For example, if $\iota=(\cdots,2,1,3,2,1,3,2,1,3)$ then
we rewrite $(\cdots,x_6,x_5,x_4,x_3,x_2,x_1)=(\cdots,x_{2,2},x_{2,1},x_{2,3},x_{1,2},x_{1,1},x_{1,3})$.

\begin{rem}\label{pos-rem}
Note that the positivity condition (\ref{posi}) implies that
for $T\in \Xi_{\iota}$ the coefficients of $x_{1,j}$ $(j\in I)$ in $T$
are non-negative.
\end{rem}

\nd
We will use the both notation $x_k$ and $x_{s,j}$.

\begin{defn}
Let us define the following (partial) ordered sets $J_{\rm A}$, $J_{\rm B}$, $J_{\rm C}$ and $J_{\rm D}$:

\begin{itemize}
\item $J_{\rm A}:=\{1,2,\cdots,n,n+1\}$ with the order $1<2<\cdots<n<n+1$.
\item
$J_{\rm B}=J_{\rm C}:=\{1,2,\cdots,n,\overline{n},\cdots,\overline{2},\overline{1}\}$ with the order
\[
1<2<\cdots<n<\overline{n}<\cdots<\overline{2}<\overline{1}.\]
\item
$J_{\rm D}:=\{1,2,\cdots,n,\overline{n},\cdots,\overline{2},\overline{1}\}$ with the partial order 
\[
 1< 2<\cdots< n-1<\ ^{n}_{\overline{n}}\ < \overline{n-1}< \cdots< \overline
 2< \overline 1.\]
\end{itemize}
For $j\in\{1,2,\cdots,n\}$, we set $|j|=|\overline{j}|=j$.

\end{defn}

\begin{defn}\label{box-def}
\begin{enumerate}
\item For $1\leq j\leq n+1$ and $s\in\mathbb{Z}$, we set
\[
\fbox{$j$}^{\rm A}_{s}:=x_{s+P(j),j}-x_{s+P(j-1)+1,j-1}\in (\QQ^{\ify})^*,\]
where $x_{m,0}=x_{m,n+1}=0$ for $m\in\mathbb{Z}$, and $x_{m,i}=0$ for $m\in\mathbb{Z}_{\leq0}$ and $i\in I$.
\item
For $1\leq j\leq n$ and $s\in\mathbb{Z}$, we set
\[
\fbox{$j$}^{\rm B}_{s}
:=x_{s+P(j),j}-x_{s+P(j-1)+1,j-1}\in (\QQ^{\ify})^*,\]
\[
\fbox{$\overline{j}$}^{\rm B}_{s}
:=x_{s+P(j-1)+n-j+1,j-1}-x_{s+P(j)+n-j+1,j}\in (\QQ^{\ify})^*,
\]
where $x_{m,0}=0$ for $m\in\mathbb{Z}$, and $x_{m,i}=0$ for $m\in\mathbb{Z}_{\leq0}$ and $i\in I$.
\item
For $1\leq j\leq n-1$ and $s\in\mathbb{Z}$, we set
\[
\fbox{$j$}^{\rm C}_{s}
:=x_{s+P(j),j}-x_{s+P(j-1)+1,j-1},
\quad
\fbox{$n$}^{\rm C}_{s}:=2x_{s+P(n),n}-x_{s+P(n-1)+1,n-1}\in (\QQ^{\ify})^*,\]
\[
\fbox{$\overline{n}$}^{\rm C}_{s}
:=x_{s+P(n-1)+1,n-1}-2x_{s+P(n)+1,n},
\quad
\fbox{$\overline{j}$}^{\rm C}_{s}
:=x_{s+P(j-1)+n-j+1,j-1}-x_{s+P(j)+n-j+1,j}\in (\QQ^{\ify})^*,
\]
\[
\fbox{$\overline{n+1}$}^{\rm C}_{s}
:=x_{s+P(n),n}
\in (\QQ^{\ify})^*,
\]
where $x_{m,0}=0$ for $m\in\mathbb{Z}$, and $x_{m,i}=0$ for $m\in\mathbb{Z}_{\leq0}$ and $i\in I$.
\item 
For $s\in\mathbb{Z}$, we set
\[
\fbox{$j$}^{\rm D}_{s}
:=x_{s+P(j),j}-x_{s+P(j-1)+1,j-1}\in (\QQ^{\ify})^*,\ \ (1\leq j\leq n-2,\ j=n),
\]
\[
\fbox{$n-1$}^{\rm D}_{s}
:=x_{s+P(n-1),n-1}+x_{s+P(n),n}-x_{s+P(n-2)+1,n-2}\in (\QQ^{\ify})^*,
\]
\[
 \fbox{$\overline{n}$}^{\rm D}_{s}
:=x_{s+P(n-1),n-1}-x_{s+P(n)+1,n}\in (\QQ^{\ify})^*,
\]
\[
\fbox{$\overline{n-1}$}^{\rm D}_{s}
:=x_{s+P(n-2)+1,n-2}-x_{s+P(n-1)+1,n-1}-x_{s+P(n)+1,n}\in (\QQ^{\ify})^*,
\]
\[
\fbox{$\overline{j}$}^{\rm D}_{s}
:=x_{s+P(j-1)+n-j,j-1}-x_{s+P(j)+n-j,j}\in (\QQ^{\ify})^*,\ \ (1\leq j\leq n-2),\]
\[
\fbox{$\overline{n+1}$}^{\rm D}_{s}
:=x_{s+P(n),n}\in (\QQ^{\ify})^*,
\]
where $x_{m,0}=0$ for $m\in\mathbb{Z}$, and $x_{m,i}=0$ for $m\in\mathbb{Z}_{\leq0}$ and $i\in I$.
\end{enumerate}
\end{defn}

\begin{lem}\label{box-lem}
\begin{enumerate}
\item In the case $\mathfrak{g}$ is of type A, the boxes $\fbox{$j$}^{\rm A}_{s}$ satisfy the following:
\begin{equation}\label{A-box}
\fbox{$j+1$}^{\rm A}_{s}=\fbox{$j$}^{\rm A}_{s} -\beta_{s+P(j),j} \qquad (1\leq j\leq n,\ s\geq 1-P(j)).
\end{equation}
\item In the case $\mathfrak{g}$ is of type B, the boxes $\fbox{$j$}^{\rm B}_{s}$
satisfy the following:
\begin{eqnarray}
\fbox{$j+1$}^{\rm B}_{s}&=&\fbox{$j$}^{\rm B}_{s} -\beta_{s+P(j),j} \qquad (1\leq j\leq n-1,\ s\geq 1-P(j)), \label{B-box1}\\
\fbox{$\overline{n}$}^{\rm B}_{s}&=&\fbox{$n$}^{\rm B}_{s} -\beta_{s+P(n),n} \qquad (s\geq 1-P(n)), \label{B-box2}\\
\fbox{$\overline{j-1}$}^{\rm B}_{s}&=&\fbox{$\overline{j}$}^{\rm B}_{s} -\beta_{s+P(j-1)+n-j+1,j-1}
 \qquad (2\leq j\leq n,\ s\geq j-P(j-1)-n). \label{B-box3}
\end{eqnarray}
\item In the case $\mathfrak{g}$ is of type C, the boxes $\fbox{$j$}^{\rm C}_{s}$
satisfy the following:
\begin{eqnarray}
\fbox{$j+1$}^{\rm C}_{s}&=&\fbox{$j$}^{\rm C}_{s} -\beta_{s+P(j),j} \qquad (1\leq j\leq n-1,\ s\geq 1-P(j)), \label{C-box1}\\
\fbox{$\overline{n}$}^{\rm C}_{s}&=&\fbox{$n$}^{\rm C}_{s} -2\beta_{s+P(n),n} \qquad (s\geq 1-P(n)), \label{C-box2}\\
\fbox{$\overline{j-1}$}^{\rm C}_{s}&=&\fbox{$\overline{j}$}^{\rm C}_{s} -\beta_{s+P(j-1)+n-j+1,j-1} \ 
 (2\leq j\leq n,\ s\geq j-P(j-1)-n),\ \ \ \ \ \label{C-box3}\\
\fbox{$\ovl{n+1}$}^{\rm C}_{l+1}+ \fbox{$\ovl{n}$}^{\rm C}_{l}
&=&\fbox{$\ovl{n+1}$}^{\rm C}_{l}-\beta_{l+P(n),n}
 \qquad (l\geq 1-P(n)). \label{BC-pr3}
\end{eqnarray}

\item In the case $\mathfrak{g}$ is of type D, the boxes $\fbox{$j$}^{\rm D}_{s}$
satisfy the following:
\begin{eqnarray}
\fbox{$j+1$}^{\rm D}_{s}&=&\fbox{$j$}^{\rm D}_{s} -\beta_{s+P(j),j} \qquad (1\leq j\leq n-1,\ s\geq 1-P(j)), \label{D-box1}\\
\fbox{$\overline{n}$}^{\rm D}_{s}&=&\fbox{$n-1$}^{\rm D}_{s} -\beta_{s+P(n),n} \qquad (s\geq 1-P(n)), \label{D-box2}\\
\fbox{$\overline{n-1}$}^{\rm D}_{s}&=&\fbox{$n$}^{\rm D}_{s} -\beta_{s+P(n),n} \qquad (s\geq 1-P(n)), \label{D-box3}\\
\fbox{$\overline{j-1}$}^{\rm D}_{s}&=&\fbox{$\overline{j}$}^{\rm D}_{s} -\beta_{s+P(j-1)+n-j,j-1} \ (2\leq j\leq n,\ s\geq 1+j-P(j-1)-n),
\qquad \quad \label{D-box4} \\
\fbox{$\overline{n+1}$}^{\rm D}_{l+2}
+\fbox{$\overline{n}$}^{\rm D}_{l+1}
+\fbox{$\overline{n-1}$}^{\rm D}_{l}
&=&\fbox{$\overline{n+1}$}^{\rm D}_{l} -\beta_{l+P(n),n} \qquad (l\geq 1-P(n)). \label{D-box5}
\end{eqnarray}
\end{enumerate}
\end{lem}

\nd
{\it Proof.}

\nd
(i) The definitions (\ref{betak}), (\ref{rewrite}) of $\beta_{s+P(j),j}$
and the definition (\ref{pij}) of $p_{i,j}$ imply that
\begin{eqnarray*}
\beta_{s+P(j),j}&=&x_{s+P(j),j} + x_{s+P(j)+1,j}-x_{s+P(j)+p_{j-1,j},j-1} -x_{s+P(j)+p_{j+1,j},j+1}\\
&=&x_{s+P(j),j} + x_{s+P(j)+1,j}-x_{s+P(j-1)+1,j-1} -x_{s+P(j+1),j+1},
\end{eqnarray*}
where we use the definition of $P(j)$ and $p_{j,j-1}+p_{j-1,j}=1$
in the second equality. Hence 
\begin{eqnarray*}
& &\fbox{$j$}^{\rm A}_{s} -\beta_{s+P(j),j}\\
&=&x_{s+P(j),j}-x_{s+P(j-1)+1,j-1} - (x_{s+P(j),j} + x_{s+P(j)+1,j}-x_{s+P(j-1)+1,j-1} -x_{s+P(j+1),j+1})\\
&=& x_{s+P(j+1),j+1}- x_{s+P(j)+1,j}=\fbox{$j+1$}^{\rm A}_{s}.
\end{eqnarray*}

\nd
(ii) The relation (\ref{B-box1}) can be shown by the same way as in (i).
We now turn to (\ref{B-box2}). It follows from (\ref{betak}), (\ref{rewrite}) that
\begin{eqnarray*}
\beta_{s+P(n),n}&=&x_{s+P(n),n} + x_{s+P(n)+1,n}-2x_{s+P(n)+p_{n-1,n},n-1} \\
&=&x_{s+P(n),n} + x_{s+P(n)+1,n}-2x_{s+P(n-1)+1,n-1},
\end{eqnarray*}
where we use the definition of $P(n)$ and $p_{n,n-1}+p_{n-1,n}=1$
in the second equality. Thus,
\begin{eqnarray*}
& &\fbox{$n$}^{\rm B}_{s} -\beta_{s+P(n),n}\\
&=&x_{s+P(n),n}-x_{s+P(n-1)+1,n-1} - (x_{s+P(n),n} + x_{s+P(n)+1,n}-2x_{s+P(n-1)+1,n-1})\\
&=& x_{s+P(n-1)+1,n-1}- x_{s+P(n)+1,n}=\fbox{$\ovl{n}$}^{\rm B}_{s}.
\end{eqnarray*}
We also get
\begin{eqnarray*}
& &\fbox{$\ovl{j}$}^{\rm B}_{s} -\beta_{s+P(j-1)+n-j+1,j-1}\\
&=&x_{s+P(j-1)+n-j+1,j-1}-x_{s+P(j)+n-j+1,j} \\
& & - (x_{s+P(j-1)+n-j+1,j-1} + x_{s+P(j-1)+n-j+2,j-1}
-x_{s+P(j-2)+n-j+2,j-2} -x_{s+P(j)+n-j+1,j})\\
&=& 
x_{s+P(j-2)+n-j+2,j-2}-x_{s+P(j-1)+n-j+2,j-1}=\fbox{$\ovl{j-1}$}^{\rm B}_{s}.
\end{eqnarray*}

\nd
(iii) Similar calculations to the proof of (i), (ii) show (\ref{C-box1})-(\ref{C-box3}).
The relation (\ref{BC-pr3}) is an easy consequence of the following calculation:
\begin{eqnarray*}
\fbox{$\ovl{n+1}$}^{\rm C}_{l}-\beta_{l+P(n),n}&=&x_{l+P(n),n}-(x_{l+P(n),n}+x_{l+P(n)+1,n}-x_{l+P(n-1)+1,n-1})\\
&=&x_{l+P(n-1)+1,n-1}-x_{l+P(n)+1,n}=\fbox{$\ovl{n+1}$}^{\rm C}_{l+1}+ \fbox{$\ovl{n}$}^{\rm C}_{l}.
\end{eqnarray*}
We can prove (iv) by a way similar to the proof of (ii) and (iii). \qed

\begin{defn}\label{tab-def}
\begin{enumerate}
\item
For ${\rm X}={\rm A}$, ${\rm B}$, ${\rm C}$ or ${\rm D}$,

\[
\begin{ytableau}
j_1 \\
j_2 \\
\vdots \\
\scriptstyle j_{k-1}\\
j_k
\end{ytableau}^{\rm X}_{s}
:=
\fbox{$
j_k$}_s^{X}
+\fbox{$j_{k-1}$}_{s+1}^{\rm X}+\cdots
+
\fbox{$
j_2 
$}^{\rm X}_{s+k-2}
+
\fbox{$
j_1
$}^{\rm X}_{s+k-1} \in (\QQ^{\ify})^*.
\]
\item
For ${\rm X}={\rm A}$, ${\rm B}$,
\[
{\rm Tab}_{{\rm X},\iota} := \{ \begin{ytableau}
j_1 \\
j_2 \\
\vdots \\
j_k
\end{ytableau}^{\rm X}_{s} | k\in I,\ j_i\in J_{X},\ s\geq 1-P(k),\ (*)^{\rm X}_k \},
\]

$(*)^{\rm A}_k : 1\leq j_1<j_2<\cdots<j_k\leq n+1$, 

$(*)^{\rm B}_k : 
\begin{cases}
1\leq j_1<j_2<\cdots<j_k\leq \overline{1} & {\rm for}\ k<n,\\
1\leq j_1<j_2<\cdots<j_n\leq \overline{1},\ \ |j_l|\neq |j_m|\ (l\neq m) & {\rm for}\ k=n.
\end{cases}$

\[
{\rm Tab}_{{\rm C},\iota}:=
\{ 
\begin{ytableau}
j_1 \\
j_2 \\
\vdots \\
j_k
\end{ytableau}^{\rm C}_{s} 
| 
\begin{array}{l}
j_1\in J_{{\rm C}}\cup\{\ovl{n+1}\},\ j_2,\cdots,j_k\in J_{{\rm C}},\\
{\rm if}\ j_1\neq \ovl{n+1}\ {\rm then}\ k\in[1,n-1]\ {\rm and}\ 1\leq j_1<j_2<\cdots<j_k\leq \overline{1},\ s\geq 1-P(k), \\
{\rm if}\ j_1=\ovl{n+1}\ {\rm then}\ k\in[1,n+1],\ \overline{n}\leq j_2<\cdots<j_k\leq\overline{1},\ s\geq 1-P(n).\\
\end{array}
\}
\]

\[
{\rm Tab}_{{\rm D},\iota}:=
\{ 
\begin{ytableau}
j_1 \\
j_2 \\
\vdots \\
j_k
\end{ytableau}^{\rm D}_{s} 
| 
\begin{array}{l}
j_1\in J_{{\rm D}}\cup\{\ovl{n+1}\},\ j_2,\cdots,j_k\in J_{{\rm D}}, \\
{\rm if}\ j_1\neq\ovl{n+1}\ {\rm then}\ k\in[1,n-2]\ {\rm and}\ j_1\ngeq j_2\ngeq\cdots\ngeq j_k,\ s\geq 1-P(k),\\
{\rm if}\ j_1=\ovl{n+1}\ {\rm and}\ k\ {\rm is\ even} \ {\rm then}\ k\in[1,n+1],\ \overline{n}\leq j_2<\cdots<j_k\leq\overline{1},
\ s\geq 1-P(n-1), \\
{\rm if}\ j_1=\ovl{n+1}\ {\rm and}\ k\ {\rm is\ odd} \ {\rm then}\ k\in[1,n+1],\ \overline{n}\leq j_2<\cdots<j_k\leq\overline{1},
\ s\geq 1-P(n).
\end{array}
\}
\]

\end{enumerate}
\end{defn}

\begin{rem}

Similar notations to Definition \ref{box-def} and \ref{tab-def} (i) can be found
in \cite{Nj, NN}.

\end{rem}

\begin{thm}\label{thm2}
For ${\rm X}={\rm A}$, ${\rm B}$, ${\rm C}$ or ${\rm D}$,
we suppose that $\iota$ is adapted to the Cartan matrix of type ${\rm X}$.
Then
\[
\Xi_{\iota}={\rm Tab}_{{\rm X},\iota}.
\]
\end{thm}

\begin{thm}\label{thm1}
In the setting of Theorem \ref{thm2}, $\iota$ satisfies the positivity condition.
\end{thm}

\begin{cor}\label{cor1}
In the setting of Theorem \ref{thm2},
we have
\[
{\rm Im}(\Psi_{\iota}) =\{\textbf{a}\in\mathbb{Z}^{\infty}_{\iota} | \varphi(\textbf{a})\geq0
\ {\rm for\ all}\ \varphi\in {\rm Tab}^{n}_{{\rm X},\iota},\ a_{m,i}=0\ {\rm for}\ m>n,\ i\in I \},
\]
where ${\rm Tab}^{n}_{{\rm X},\iota}:=
\{\begin{ytableau}
j_1 \\
\vdots \\
j_k
\end{ytableau}^{\rm X}_{s} \in{\rm Tab}_{{\rm X},\iota} | s\leq n \}$.
\end{cor}

\nd
The proofs of Theorem \ref{thm2}, \ref{thm1} and Corollary \ref{cor1} will be given in Sect. \ref{pr-sect}.

\begin{ex}

Let $\mathfrak{g}$ be the Lie algebra of type ${\rm A}_3$ and $\iota=(\cdots,3,1,2,3,1,2)$.
The sequence $\iota$ is adapted to the Cartan matrix of type ${\rm A}_3$. We get
$p_{2,1}=1$, $p_{3,2}=0$, $P(2)=P(3)=1$ and
\begin{eqnarray}
{\rm Tab}_{{\rm A},\iota}&=&\{\fbox{$j$}^{\rm A}_s| s\geq 1,\ j\in[1,4] \}
\cup\{
\begin{ytableau}
i \\
j 
\end{ytableau}^{\rm A}_{s} 
| 
\begin{array}{l}
 s\geq 0,\\
 1\leq i<j\leq 4. 
\end{array}
\}
\cup\{
\begin{ytableau}
i \\
j \\
k 
\end{ytableau}^{\rm A}_{s} 
| 
\begin{array}{l}
 s\geq 0,\\
 1\leq i<j<k\leq 4. 
\end{array}
\}\nonumber \\
&=& 
\{x_{s,1},\ x_{s+1,2}-x_{s+1,1},\ x_{s+1,3}-x_{s+2,2},\ -x_{s+2,3} | s\geq 1 \}\nonumber\\
&\cup&
\left\{x_{s+1,2},\ x_{s+1,3}-x_{s+2,2}+x_{s+1,1},\ x_{s+1,1}-x_{s+2,3},\ x_{s+1,3}-x_{s+2,1},\right.\nonumber\\
& &\left. x_{s+2,2}-x_{s+2,1}-x_{s+2,3},\ -x_{s+3,2} | s\geq 0 \right\}\label{ex-1}\\
&\cup&
\{x_{s+1,3},\ x_{s+2,2}-x_{s+2,3},\ x_{s+2,1}-x_{s+3,2},\ -x_{s+3,1} | s\geq 0 \}.\nonumber
\end{eqnarray}
By Theorem \ref{thm2}, we have $\Xi_{\iota}={\rm Tab}_{{\rm A},\iota}$.
The explicit form (\ref{ex-1}) means $\iota$ satisfies the positivity condition. Hence,
\[
{\rm Im}(\Psi_{\iota})=\Sigma_{\io}
=\{\textbf{x}\in\mathbb{Z}^{\infty}_{\iota} | \varphi(\textbf{x})\geq0, \ \ \forall\varphi\in\Xi_{\iota} \}.
\]
For $\textbf{x}=(\cdots,x_{2,3},x_{2,1},x_{2,2},x_{1,3},x_{1,1},x_{1,2})\in {\rm Im}(\Psi_{\iota})$,
combining inequalities $x_{s,1}\geq0$ $(s\geq1)$, $-x_{s+3,1}\geq0$ $(s\geq0)$,
we obtain
$x_{s+3,1}=0$ $(s\geq0)$.
Similarly, by $x_{s+1,2}\geq0$, $-x_{s+3,2}\geq0$ $(s\geq0)$,
we get $x_{s+3,2}=0$ $(s\geq0)$. We also get $x_{s+2,3}=0$ $(s\geq 1)$.
Hence, simplifying the inequalities, we obtain
\begin{eqnarray*}
{\rm Im}(\Psi_{\iota})&=&\left\{\textbf{x}\in \mathbb{Z}^{\infty}_{\iota} | x_{s,1}=x_{s,2}=x_{s,3}=0\ {\rm for}\ s\in\mathbb{Z}_{\geq3},\ 
\ x_{2,2}- x_{2,1}\geq x_{2,3}\geq0,\right.\\
& & \left. x_{1,3}+x_{1,1}\geq x_{2,2}\geq0,\ x_{1,1}\geq x_{2,3}\geq0,\ x_{1,3}\geq x_{2,1}\geq0,\ x_{1,2}\geq0 \right\}.
\end{eqnarray*}

\end{ex}

\begin{ex}

Let $\mathfrak{g}$ be the Lie algebra of type ${\rm C}_3$ and $\iota=(\cdots,3,1,2,3,1,2)$.
The sequence $\iota$ is adapted to the Cartan matrix of type ${\rm C}_3$. We get
$p_{2,1}=1$, $p_{3,2}=0$, $P(2)=P(3)=1$ and
\begin{eqnarray}
{\rm Tab}_{{\rm C},\iota}&=&\{\fbox{$j$}^{\rm C}_s| s\geq 1,\ 1\leq j\leq \ovl{1}\}
\cup\{
\begin{ytableau}
i \\
j 
\end{ytableau}^{\rm C}_{s} 
| 
\begin{array}{l}
 s\geq 0,\\
 1\leq i<j\leq \ovl{1}. 
\end{array}
\}
\cup\{
\begin{ytableau}
\ovl{4} \\
j_2 \\
\vdots\\
j_k 
\end{ytableau}^{\rm C}_{s} 
| 
\begin{array}{l}
 s\geq 0,\\
 k\in[1,4],\ 
 \ovl{3}\leq j_2<\cdots<j_k\leq \ovl{1}. 
\end{array}
\}\nonumber \\
&=& 
\{x_{s,1},\ x_{s+1,2}-x_{s+1,1},\ 2x_{s+1,3}-x_{s+2,2},\ x_{s+2,2}-2x_{s+2,3},\ x_{s+2,1}-x_{s+3,2},\ -x_{s+3,1}| s\geq 1 \}\nonumber\\
&\cup&
\left\{x_{s+1,2},\ 2x_{s+1,3}-x_{s+2,2}+x_{s+1,1},\ x_{s+1,1}+x_{s+2,2}-2x_{s+2,3},\ x_{s+1,1}+x_{s+2,1}-x_{s+3,2},
\right.\nonumber\\
& &\left. x_{s+1,1}-x_{s+3,1},\ 2x_{s+1,3}-x_{s+2,1},\ 2x_{s+2,2}-x_{s+2,1}-2x_{s+2,3},\ x_{s+2,2}-x_{s+3,2},\ 
x_{s+2,2}-x_{s+2,1} -x_{s+3,1}, \right. \nonumber \\
& & \left. 2x_{s+2,3}-2x_{s+3,2} + x_{s+2,1},\ 2x_{s+2,3}-x_{s+3,2} - x_{s+3,1},\ x_{s+2,1}-2x_{s+3,3},\right. \label{ex-2} \\
& & \left. x_{s+3,2}-2x_{s+3,3}-x_{s+3,1},\ -x_{s+4,2} | s\geq 0 \right\} \nonumber \\
&\cup&
\left\{x_{s+1,3},\ x_{s+2,2}-x_{s+2,3},\ x_{s+2,3}+x_{s+2,1}-x_{s+3,2},\ x_{s+2,3}-x_{s+3,1},\ x_{s+2,1}-x_{s+3,3}, \right. \nonumber\\
& & \left. x_{s+3,2}-x_{s+3,1}-x_{s+3,3},\ x_{s+3,3}-x_{s+4,2},\ -x_{s+4,3} | s\geq 0 \right\}.\nonumber
\end{eqnarray}
By Theorem \ref{thm2}, we have $\Xi_{\iota}={\rm Tab}_{{\rm C},\iota}$.
The explicit form (\ref{ex-2}) means $\iota$ satisfies the positivity condition. Hence,
\[
{\rm Im}(\Psi_{\iota})=\Sigma_{\io}
=\{\textbf{x}\in\mathbb{Z}^{\infty}_{\iota} | \varphi(\textbf{x})\geq0, \ \ \forall\varphi\in\Xi_{\iota} \}.
\]
Simplifying the inequalities, we obtain
\begin{eqnarray*}
& &{\rm Im}(\Psi_{\iota})\\
&=&
\left\{\textbf{x}\in \mathbb{Z}^{\infty}_{\iota} | x_{s,1}=x_{s,2}=x_{s,3}=0\ {\rm for}\ s\in\mathbb{Z}_{\geq4},\ 
\ x_{2,2}\geq x_{2,1}\geq0, \ 
2x_{2,3}\geq x_{3,2}\geq x_{3,1}\geq 0, \right.\\
& &\left. x_{3,2}-2x_{3,3}\geq0,\ 2x_{1,3}-x_{2,2}+x_{1,1}\geq0, \right. \\ 
& &\left. x_{1,1}+x_{2,2}-2x_{2,3}\geq0,\ 
 x_{1,1}+x_{2,1}-x_{3,2}\geq0,\ x_{1,1}\geq x_{3,1}\geq0,\ 
2x_{1,3}-x_{2,1}\geq0, \right.\\
& &\left. 2x_{2,2}-x_{2,1}-2x_{2,3}\geq0,\ x_{2,2}-x_{3,2}\geq0,\ x_{2,2}-x_{2,1}-x_{3,1}\geq0,\ 2x_{2,3}-2x_{3,2}+x_{2,1}\geq0, 
  \right.\\
& & \left. 2x_{2,3}- x_{3,2}-x_{3,1}\geq0,\ x_{2,1}- 2x_{3,3}\geq0,\ x_{3,2}- 2x_{3,3}-x_{3,1}\geq0, \right. \\
& & \left. x_{2,2}-x_{2,3}\geq0,\ x_{2,3}+x_{2,1}-x_{3,2}\geq0,\ x_{2,3}\geq x_{3,1}\geq0,\ x_{2,1}\geq x_{3,3}\geq0,\ 
x_{3,2}-x_{3,1}-x_{3,3}\geq0, \right. \\
& &\left. x_{1,2}\geq0,\ x_{1,3}\geq0  \right\}.
\end{eqnarray*}


\end{ex}

\section{Stability of ${\rm Tab}_{{\rm X},\iota}$}\label{Stabsect}

In what follows, we denote each tableau $\begin{ytableau}
j_1 \\
\vdots \\
j_k
\end{ytableau}^{\rm X}_{s}$ by $[j_1,\cdots,j_k]_s^{\rm X}$.
When we see the condition $j_l\neq t$ with $t\in\{1,2,\cdots,n,n+1,\ovl{n},\cdots,\ovl{1}\}$
for $[j_1,\cdots,j_k]_s^{\rm X}$ it means $j_l\neq t$ with $l\in[1,k]$ or $l>k$ or $l<1$.

\subsection{Stability of ${\rm Tab}_{{\rm A},\iota}$}

\begin{prop}\label{closednessA}

For each $T=[j_1,\cdots,j_k]_s^{\rm A}\in {\rm Tab}_{{\rm A},\iota}$, $m\in\mathbb{Z}_{\geq1}$ and $j\in I$,
\begin{eqnarray*}
\hspace{-15mm}& &S_{m,j}T=\\
\hspace{-15mm}& &\begin{cases}
[j_1,\cdots,j_{i-1},j+1,j_{i+1},\cdots,j_k]_s^{\rm A} & {\rm if}\ j_i=j,\ j_{i+1}\neq j+1,\ m=s+k-i+P(j)
\ {\rm for\ some}\ i\in[1,k],\\
[j_1,\cdots,j_{i-1},j,j_{i+1},\cdots,j_k]_s^{\rm A} & {\rm if}\ j_i=j+1,\ j_{i-1}\neq j,\ m=s+k-i+1+P(j)\ {\rm for\ some}\ i\in[1,k], \\
T & {\rm otherwise}.
\end{cases}
\end{eqnarray*}
In particular, the set ${\rm Tab}_{{\rm A},\iota}$ is closed under the action of $S_{m,j}$.
\end{prop}

\nd
{\it Proof.}

\nd
Let $T\in{\rm Tab}_{{\rm A},\iota}$ be the following:
\begin{equation}\label{A-pr01}
 T=
[j_1,\cdots,j_k]^{\rm A}_s 
=\sum^{k}_{i=1} \fbox{$j_i$}^{\rm A}_{s+k-i} 
=\sum^{k}_{i=1} (x_{s+k-i+P(j_i),j_i}-x_{s+k-i+1+P(j_i-1),j_{i}-1}) 
\end{equation}
with $s\geq 1-P(k)$, $k\in I$ and $1\leq j_1<\cdots<j_k\leq n+1$. Note that since
\begin{eqnarray*}
& &
\fbox{$j_i$}^{\rm A}_{s+k-i}
+
\fbox{$j_{i+1}$}^{\rm A}_{s+k-i-1}\\ 
&=&
x_{s+k-i+P(j_i),j_i}-x_{s+k-i+1+P(j_i-1),j_i-1}
+x_{s+k-i-1+P(j_{i+1}),j_{i+1}}-x_{s+k-i+P(j_{i+1}-1),j_{i+1}-1}, 
\end{eqnarray*}
if $j_{i+1}=j_{i}+1$ then we get
\begin{equation}\label{A-pr02}
\fbox{$j_i$}^{\rm A}_{s+k-i}
+
\fbox{$j_{i}+1$}^{\rm A}_{s+k-i-1}=
x_{s+k-i-1+P(j_{i}+1),j_{i}+1}-x_{s+k-i+1+P(j_i-1),j_i-1}.
\end{equation}
Hence if $j_{i+1}=j_{i}+1$ then the coefficients of $x_{\xi,j_i}$ $(\xi\in\mathbb{Z}_{\geq1})$ in $T$ are $0$.
It follows from (\ref{A-pr01}) and (\ref{A-pr02}) that
$x_{m,j}$ has non-zero coefficient in $T$ if and only if
the pair $(m,j)$ belongs to
\begin{multline*}
\{ (s+k-i+P(j_i),j_i) | i=1,2,\cdots,k,\ \ j_{i+1}> j_i+1\} \\
\cup \{(s+k-i+1+P(j_i-1),j_{i}-1) | i=1,2,\cdots,k,\ \ j_i-1> j_{i-1}\},
\end{multline*}
where we set $j_{k+1}=n+2$ and $j_0=0$.
Note that
by (\ref{pineq}), (\ref{pineq2}), $j_i\geq i$ and $s\geq 1-P(k)$,
we get
\begin{equation}\label{A-m1}
s+k-i+P(j_i)
\geq s+k-i+P(i)
\geq s+P(k) \geq 1,
\end{equation}
which yields
$s+k-i\geq 1-P(j_i)$ for $i\in[1,k]$.
If $(m,j)=(s+k-i+P(j_i),j_i)$ with $j_{i+1}> j_i+1$ then $x_{m,j}$ has coefficient $1$ in $T$
and by Lemma \ref{box-lem} (\ref{A-box}) and the definition of $S_{m,j}$,
\begin{eqnarray*}
S_{m,j} T
&=&T-\beta_{m,j}
= [j_1,\cdots,j_{i-1},j,j_{i+1},\cdots,j_k]^{\rm A}_s -\beta_{m,j}\\
&=&\fbox{$j_1$}^{\rm A}_{s+k-1}+\cdots+
\fbox{$j_{i-1}$}^{\rm A}_{s+k-i+1}+
\fbox{$j$}^{\rm A}_{s+k-i}+
\fbox{$j_{i+1}$}^{\rm A}_{s+k-i-1}+
\cdots +\fbox{$j_k$}^{\rm A}_{s} -\beta_{m,j}\\
&=&
\fbox{$j_1$}^{\rm A}_{s+k-1}+\cdots+
\fbox{$j_{i-1}$}^{\rm A}_{s+k-i+1}+
\fbox{$j+1$}^{\rm A}_{s+k-i}+
\fbox{$j_{i+1}$}^{\rm A}_{s+k-i-1}+
\cdots +\fbox{$j_k$}^{\rm A}_{s} \\
&=&
[j_1,\cdots,j_{i-1},j+1,j_{i+1},\cdots,j_k]^{\rm A}_s 
\in {\rm Tab}_{{\rm A},\iota}.
\end{eqnarray*}
If $(m,j)=(s+k-i+1+P(j_i-1),j_{i}-1)$ with $j_i-1> j_{i-1}$ then $x_{m,j}$ has coefficient $-1$ in $T$.
By (\ref{pineq}), (\ref{pineq2}), $s\geq 1-P(k)$ and 
$j_i-1>j_{i-1}\geq i-1$ so that $j_i-1\geq i$, we get
\begin{equation}\label{A-pr-0}
s+k-i+1+P(j_i-1)
\geq s+k-i+1+P(i)\geq s+P(k)+1 \geq 2. 
\end{equation}
Thus, owing to (\ref{Sk}) and Lemma \ref{box-lem} (\ref{A-box}),
\begin{eqnarray*}
S_{m,j} T&=&S_{m,j}
[j_1,\cdots,j_i,\cdots,j_k]^{\rm A}_s 
=
[j_1,\cdots,j_i,\cdots,j_k]^{\rm A}_s +\beta_{m-1,j}\\
&=&  
\fbox{$j_1$}^{\rm A}_{s+k-1}+\cdots+\fbox{$j_i$}^{\rm A}_{s+k-i}+\cdots +\fbox{$j_k$}^{\rm A}_{s}+\beta_{m-1,j}\\
&=&  
\fbox{$j_1$}^{\rm A}_{s+k-1}+\cdots+\fbox{$j+1$}^{\rm A}_{s+k-i}+\cdots +\fbox{$j_k$}^{\rm A}_{s}+\beta_{m-1,j}\\
&=&\fbox{$j_1$}^{\rm A}_{s+k-1}+\cdots+\fbox{$j$}^{\rm A}_{s+k-i}+\cdots +\fbox{$j_k$}^{\rm A}_{s}
=[j_1,\cdots,j,\cdots,j_k]^{\rm A}_s
\in {\rm Tab}_{{\rm A},\iota}.
\end{eqnarray*}
The definition of $S_{m,j}$ means that if $x_{m,j}$ is not a summand of $T$ then $S_{m,j}T=T$.
Thus, we get our claim. \qed

\subsection{Stabilities of ${\rm Tab}_{{\rm B},\iota}$ and ${\rm Tab}_{{\rm C},\iota}$}

\begin{prop}\label{closednessBC}
For each $T=[j_1,\cdots,j_k]_s^{\rm X}\in {\rm Tab}_{{\rm X},\iota}$ $({\rm X}={\rm B}$ or ${\rm C})$,
$j\in I$ and $m\in\mathbb{Z}_{\geq1}$, we consider the following conditions for the triple $(T,j,m)$:
\begin{enumerate}
\item[(1)] $j<n$ and there exists $i\in [1,k]$ such that $j_i=j$, $j_{i+1}\neq j+1$ and $m=s+k-i+P(j)$,
\item[(2)] $j<n$ and there exists $i'\in [1,k]$ such that $j_{i'}=\overline{j+1}$, $j_{i'+1}\neq \overline{j}$
and $m=s+k-i'+n-j+P(j)$,
\item[(3)] $j=n$ and there exists $i\in [1,k]$ such that $j_i=n$, $j_{i+1}\neq \overline{n}$ and $m=s+k-i+P(n)$,
\item[(4)] $j<n$ and there exists $i\in [1,k]$ such that $j_{i-1}\neq j$, $j_{i}= j+1$ and $m=s+k-i+1+P(j)$,
\item[(5)] $j<n$ and there exists $i'\in [1,k]$ such that $j_{i'-1}\neq\overline{j+1}$, $j_{i'}= \overline{j}$
and $m=s+k-i'+n-j+1+P(j)$,
\item[(6)] $j=n$ and there exists $i\in [1,k]$ such that $j_{i-1}\neq n$, $j_{i}= \overline{n}$ and $m=s+k-i+1+P(n)$,
\item[(7)] $j=n$, $j_1=\ovl{n+1}$, $j_2\neq\ovl{n}$ and $m=s+k-1+P(n)$.
\end{enumerate}
\begin{enumerate}
\item
If $k\in[1,n-1]$ and $j_1\neq \ovl{n+1}$ then we have
\begin{eqnarray*}
\hspace{-15mm}& &S_{m,j}T= \\
\hspace{-15mm}& &\begin{cases}
[j_1,\cdots,j_{i-1},j+1,j_{i+1},\cdots,j_k]_s^{\rm X} & {\rm if\ (1)\ holds\ and\ (2),\ (5)\ do\ not\ hold}, \\
[j_1,\cdots,j_{i'-1},\ovl{j},j_{i'+1},\cdots,j_k]_s^{\rm X} & {\rm if\ (2)\ holds\ and\ (1),\ (4)\ do\ not\ hold}, \\
[j_1,\cdots,j_{i-1},j+1,j_{i+1},\cdots,j_{i'-1},\ovl{j},j_{i'+1},\cdots,j_k]_s^{\rm X} & {\rm if\ (1)\ and\ (2)\ hold}, \\
[j_1,\cdots,j_{i-1},\ovl{n},j_{i+1}\cdots,,j_k]_s^{\rm X} & {\rm if\ (3)\ holds}, \\
[j_1,\cdots,j_{i-1},j,j_{i+1},\cdots,j_k]_s^{\rm X} & {\rm if\ (4)\ holds\ and\ (2),\ (5)\ do\ not\ hold}, \ \ \\
[j_1,\cdots,j_{i'-1},\ovl{j+1},j_{i'+1},\cdots,j_k]_s^{\rm X} & {\rm if\ (5)\ holds\ and\ (1),\ (4)\ do\ not\ hold}, \\
[j_1,\cdots,j_{i-1},j,j_{i+1},\cdots,j_{i'-1},\ovl{j+1},j_{i'+1},\cdots,j_k]_s^{\rm X} & {\rm if\ (4)\ and\ (5)\ hold}, \\
[j_1,\cdots,j_{i-1},n,j_{i+1}\cdots,,j_k]_s^{\rm X} & {\rm if\ (6)\ holds}, \\
T & {\rm otherwise}.
\end{cases}
\end{eqnarray*}
\item For each $T=[j_1,\cdots,j_n]_s^{\rm B}\in {\rm Tab}_{{\rm B},\iota}$ (so that $k=n$),
 $m\in\mathbb{Z}_{\geq1}$ and $j\in I$, we have
\[S_{m,j}T=
\begin{cases}
[j_1,\cdots,j_{i-1},j+1,j_{i+1},\cdots,j_{i'-1},\ovl{j},j_{i'+1},\cdots,j_n]_s^{\rm B} & {\rm if\ (1)\ and\ (2)\ hold}, \\
[j_1,\cdots,j_{i-1},\ovl{n},j_{i+1}\cdots,,j_n]_s^{\rm B} & {\rm if\ (3)\ holds}, \\
[j_1,\cdots,j_{i-1},j,j_{i+1},\cdots,j_{i'-1},\ovl{j+1},j_{i'+1},\cdots,j_n]_s^{\rm B} & {\rm if\ (4)\ and\ (5)\ hold}, \\
[j_1,\cdots,j_{i-1},n,j_{i+1}\cdots,,j_n]_s^{\rm B} & {\rm if\ (6)\ holds}, \\
T & {\rm otherwise}.
\end{cases}
\]
\item For each $T=[\ovl{n+1},j_2,j_3,\cdots,j_k]_s^{\rm C}\in {\rm Tab}_{{\rm C},\iota}$ with $k\in[1,n+1]$,
$m\in\mathbb{Z}_{\geq1}$ and $j\in I$, we have
\[S_{m,j}T=
\begin{cases}
[\ovl{n+1},j_2,\cdots,j_{i'-1},\ovl{j},j_{i'+1},\cdots,j_k]_s^{\rm C} & {\rm if\ (2)\ holds}, \\
[\ovl{n+1},j_2,\cdots,j_{i'-1},\ovl{j+1},j_{i'+1},\cdots,j_k]_s^{\rm C} & {\rm if\ (5)\ holds}, \\
[\ovl{n+1},j_3,\cdots,j_k]_s^{\rm C} & {\rm if\ (6)\ holds}, \\
[\ovl{n+1},\ovl{n},j_2,\cdots,j_k]_s^{\rm C} & {\rm if\ (7)\ holds}, \\
T & {\rm otherwise}.
\end{cases}
\]
\end{enumerate}
In particular, the sets ${\rm Tab}_{{\rm B},\iota}$
and ${\rm Tab}_{{\rm C},\iota}$
are closed under the action of $S_{m,j}$.
\end{prop}

\nd
{\it Proof.}

\nd
(i)
For $T=[j_1,\cdots,j_k]^{\rm X}_s\in {\rm Tab}_{{\rm X},\iota}$ with $k\in[1,n-1]$ and $j_1\neq \ovl{n+1}$,
let us consider the action of $S_{m,j}$ $(m\in\mathbb{Z}_{\geq1},\ j\in I)$.
Recall that 
$T=[j_1,\cdots,j_k]^{\rm X}_s=\sum_{i=1}^{k}\fbox{$j_i$}^{\rm X}_{s+k-i}$, and
by Definition \ref{box-def}, we obtain
\begin{equation}
\fbox{$j_i$}^{\rm X}_{s+k-i}=
\begin{cases}\label{BC-box}
c(j_i)x_{s+k-i+P(j_i),j_i}-x_{s+k-i+P(j_i-1)+1,j_i-1} & {\rm if}\ j_i\leq n,\\
x_{s+k-i+P(|j_i|-1)+n-|j_i|+1,|j_i|-1}-c(j_i)x_{s+k-i+P(|j_i|)+n-|j_i|+1,|j_i|} & {\rm if}\ j_i\geq \ovl{n},
\end{cases}
\end{equation}
where if $\mathfrak{g}$ is of type C and $j_i\in\{n,\ovl{n}\}$ then $c(j_i)=2$,
otherwise $c(j_i)=1$.
By $s\geq 1-P(k)$ and a similar calculation to (\ref{A-m1}), (\ref{A-pr-0}), for the 
left indices in (\ref{BC-box}), in the case $j_i\leq n$ it holds
\begin{equation}\label{BC-rev1}
s+k-i+P(j_i)\geq1,
\end{equation}
and if $j_i>j_{i-1}+1$ then
\begin{equation}\label{BC-rev2}
s+k-i+P(j_i-1)+1\geq2.
\end{equation}
Just as in (\ref{A-pr02}), if $j_i=j_{i-1}+1$ then $-x_{s+k-i+P(j_i-1)+1,j_i-1}$ is cancelled in $T$.
Using $i\leq k$, $s\geq 1-P(k)$, (\ref{pineq}) and (\ref{pineq2}), in the case $j_i\geq\ovl{n}$, we also see that
\begin{equation}\label{BCm1}
\begin{array}{l}
 s+k-i+P(|j_i|-1)+n-|j_i|+1\geq s+k-i+P(n)\geq 1+k-i+P(n)-P(k)\geq1,\\
 s+k-i+P(|j_i|)+n-|j_i|+1\geq s+ P(n) +1 \geq 2.
\end{array}
\end{equation}

Clearly, if $j$, $j+1$, $\ovl{j+1}$, $\ovl{j}\notin\{j_1,\cdots,j_k\}$
then the coefficient of $x_{m,j}$ in $T$ is $0$ so that $S_{m,j}T=T$. 
Hence, we consider the case at least one of $j$, $j+1$, $\ovl{j+1}$, $\ovl{j}$ belongs to
$\{j_1,\cdots,j_k\}$ and $x_{m,j}$ is a summand of $T$ with non-zero coefficient.
By (\ref{BC-box}) and a similar calculation to (\ref{A-pr02}),
the condition (1) in the claim is equivalent to that $j_i=j<n$ and $x_{m,j}$
is a summand of
$\fbox{$j_i$}^{\rm X}_{s+k-i}+
\fbox{$j_{i+1}$}^{\rm X}_{s+k-i-1}=
\fbox{$j$}^{\rm X}_{s+k-i}+
\fbox{$j_{i+1}$}^{\rm X}_{s+k-i-1}$
with a positive coefficient for some $i\in[1,k]$. The coefficient is $1$ except for the case
\begin{equation}\label{exbc1}
j_i=j=n-1,\ \ j_{i+1}=\ovl{n}.
\end{equation}
In the case (\ref{exbc1}), the coefficient of $x_{m,j}$ in $\fbox{$j_i$}^{\rm X}_{s+k-i}+
\fbox{$j_{i+1}$}^{\rm X}_{s+k-i-1}$ is $2$.

Note that if $\ovl{n}\leq j_{i'}$ then
\begin{eqnarray*}
& &
\fbox{$j_{i'}$}^{\rm X}_{s+k-i'}+
\fbox{$j_{i'+1}$}^{\rm X}_{s+k-i'-1}\\
&=&
x_{s+k-i'+P(|j_{i'}|-1)+n-|j_{i'}|+1,|j_{i'}|-1}-c(j_{i'})x_{s+k-i'+P(|j_{i'}|)+n-|j_{i'}|+1,|j_{i'}|}\\
& &+x_{s+k-i'+P(|j_{i'+1}|-1)+n-|j_{i'+1}|,|j_{i'+1}|-1}-x_{s+k-i'+P(|j_{i'+1}|)+n-|j_{i'+1}|,|j_{i'+1}|}.
\end{eqnarray*}
In particular, if $j_{i'}=\ovl{j+1}$ and $j_{i'+1}=\ovl{j}$ $(j\in[1,n-1])$ then
\begin{eqnarray}
\fbox{$\ovl{j+1}$}_{s+k-i'}^{\rm X}+\fbox{$\ovl{j}$}_{s+k-i'-1}^{\rm X}
&=&x_{s+k-i'+P(j)+n-j,j}-c(j+1)x_{s+k-i'+P(j+1)+n-j,j+1} \nonumber\\
& &+x_{s+k-i'+P(j-1)+n-j,j-1}-x_{s+k-i'+P(j)+n-j,j}\nonumber\\
&=&x_{s+k-i'+P(j-1)+n-j,j-1}-c(j+1)x_{s+k-i'+P(j+1)+n-j,j+1}.\label{BC-pr1}
\end{eqnarray}
Hence, the condition (2) in the claim is equivalent to that $\ovl{n}\leq j_{i'}=\ovl{j+1}$ and $x_{m,j}$
is a summand of
$\fbox{$j_{i'}$}^{\rm X}_{s+k-i'}+
\fbox{$j_{i'+1}$}^{\rm X}_{s+k-i'-1}=
\fbox{$\ovl{j+1}$}^{\rm X}_{s+k-i'}+
\fbox{$j_{i'+1}$}^{\rm X}_{s+k-i'-1}$
with the coefficient $1$ for some $i'\in[1,k]$.
Similarly, we can verify that the
condition (3) is equivalent to $j=j_i=n$ and $x_{m,n}$
is a summand of
$\fbox{$j_{i}$}^{\rm X}_{s+k-i}+
\fbox{$j_{i+1}$}^{\rm X}_{s+k-i-1}=
\fbox{$n$}^{\rm X}_{s+k-i}+
\fbox{$j_{i+1}$}^{\rm X}_{s+k-i-1}$
with the coefficient $c(n)$ for some $i\in I$.
We also see that the condition (4) (resp. (6)) holds if and only if
$j_{i}= j+1\leq n$ (resp. $j_{i}= \overline{n}=\ovl{j}$)
and $x_{m,j}$ is a summand of
$\fbox{$j_{i}$}^{\rm X}_{s+k-i}+
\fbox{$j_{i-1}$}^{\rm X}_{s+k-i+1}$
with the coefficient $-1$ (resp. $-c(n)$) for some $i\in I$.
The condition (5) is equivalent to $\ovl{n}<j_{i'}= \overline{j}$ and 
$x_{m,j}$ is a summand of
$\fbox{$j_{i'}$}^{\rm X}_{s+k-i'}+
\fbox{$j_{i'-1}$}^{\rm X}_{s+k-i'+1}$
with a negative coefficient for some $i'\in I$. The coefficient is $-1$ except for the case
\begin{equation}\label{exbc2}
j_{i'}=\ovl{n-1},\ \ j_{i'-1}=n.
\end{equation}
In the case $j_{i'}=\ovl{n-1}$, $j_{i'-1}=n$, the coefficient of $x_{m,j}$ in
$\fbox{$j_{i'}$}^{\rm X}_{s+k-i'}+
\fbox{$j_{i'-1}$}^{\rm X}_{s+k-i'+1}$ is $-2$.

Note that if (1) and (5) hold then
$m=s+k-i+P(j)=s+k-i'+n-j+1+P(j)$ so that $i'=i+n-j+1$.
In the case $j=n-1$, we have $i'=i+2$ and $j_{i+1}=n$ or $\ovl{n}$, which contradicts the assumption (1) and (5) hold.
In particular (\ref{exbc1}),  (\ref{exbc2}) do not hold.
Thus, $j<n-1$ so that $(i'-i)>2$.
By the above argument, the coefficient of $x_{m,j}$ in $T$ is $0$.
Similarly, if (2) and (4) hold then the coefficient of $x_{m,j}$ in $T$ is $0$.

Now, we suppose the coefficient $c_{m,j}$ of $x_{m,j}$ in $T$ is positive.
Under this assumption, at least one of (1), (2), (3) holds. 
Clearly, if (3) holds then neither (1) nor (2) hold. 
The conditions (1) and (4) (resp. (2) and (5)) do not hold simultaneously.
Thus, we need to consider the following four cases:

\vspace{2mm}

\nd
\underline{Case 1. (1) holds and (2), (5) do not hold}

\vspace{2mm}

\nd
In this case, since (2) does not hold,
if (\ref{exbc1}) holds then $j_{i+2}=\ovl{n-1}$, which yields
that
the coefficient of $x_{m,j}=x_{m,n-1}$ in $\fbox{$j_i$}^{\rm X}_{s+k-i}+
\fbox{$j_{i+1}$}^{\rm X}_{s+k-i-1}+
\fbox{$j_{i+2}$}^{\rm X}_{s+k-i-2}$ so that in $T$ is $1$ by (\ref{BC-pr1}).

Hence, in both the case (\ref{exbc1}) holds and the case (\ref{exbc1}) does not hold,
we obtain $c_{m,j}=1$ and $m=s+k-i+P(j)$.
By Lemma \ref{box-lem},
\begin{eqnarray*}
S_{m,j}T&=&S_{m,j}[j_1,\cdots,j_{i-1},j_i,j_{i+1},\cdots,j_k]^{\rm X}_s\\
&=& 
[j_1,\cdots,j_{i-1},j_i,j_{i+1},\cdots,j_k]^{\rm X}_s - \beta_{m,j}\\
&=&
\fbox{$j_1$}^{\rm X}_{s+k-1}+\cdots+\fbox{$j_i$}^{\rm X}_{s+k-i}+\cdots +\fbox{$j_k$}^{\rm X}_{s}-\beta_{m,j}\\
&=&
\fbox{$j_1$}^{\rm X}_{s+k-1}+\cdots+\fbox{$j$}^{\rm X}_{s+k-i}+\cdots +\fbox{$j_k$}^{\rm X}_{s}-\beta_{m,j}\\
&=&
\fbox{$j_1$}^{\rm X}_{s+k-1}+\cdots+\fbox{$j+1$}^{\rm X}_{s+k-i}+\cdots +\fbox{$j_k$}^{\rm X}_{s}=
[j_1,\cdots,j_{i-1},j+1,j_{i+1},\cdots,j_k]^{\rm X}_s \in {\rm Tab}_{{\rm X},\iota}.
\end{eqnarray*}

\vspace{2mm}

\nd
\underline{Case 2. (2) holds and (1), (4) do not hold}

\vspace{2mm}

In this case,
we obtain $c_{m,j}=1$ and $m=s+k-i'+n-j+P(j)$.
Due to Lemma \ref{box-lem},
\begin{eqnarray*}
S_{m,j}T&=&S_{m,j}[j_1,\cdots,j_{i'-1},j_{i'},j_{i'+1},\cdots,j_k]^{\rm X}_s\\
&=&
\fbox{$j_1$}^{\rm X}_{s+k-1}+\cdots+\fbox{$j_{i'}$}^{\rm X}_{s+k-i'}+\cdots +\fbox{$j_k$}^{\rm X}_{s}-\beta_{m,j}\\
&=&
\fbox{$j_1$}^{\rm X}_{s+k-1}+\cdots+\fbox{$\ovl{j+1}$}^{\rm X}_{s+k-i'}+\cdots +\fbox{$j_k$}^{\rm X}_{s}-\beta_{m,j}\\
&=&
\fbox{$j_1$}^{\rm X}_{s+k-1}+\cdots+\fbox{$\ovl{j}$}^{\rm X}_{s+k-i'}+\cdots +\fbox{$j_k$}^{\rm X}_{s}=
[j_1,\cdots,j_{i'-1},\ovl{j},j_{i'+1},\cdots,j_k]^{\rm X}_s \in {\rm Tab}_{{\rm X},\iota}.
\end{eqnarray*}

\vspace{2mm}

\nd
\underline{Case 3. Both (1) and (2) hold}

\vspace{2mm}

In this case, 
we have $m=s+k-i+P(j)$ and $m=s+k-i'+n-j+P(j)$ so that $i'=i+n-j$.
If $j=n-1$ then $j_i=n-1$, $j_{i+1}=\ovl{n}$ and $j_{i+2}\neq\ovl{n-1}$. Thus,
(\ref{exbc1}) holds, which yields
that
the coefficient of $x_{m,j}=x_{m,n-1}$ in $T$ is $2$.

If $j\leq n-2$ then the equation $i'=i+n-j$ means $(i'-i)\geq2$.
Hence, we obtain $c_{m,j}=2$ for $j\in[1,n-1]$.
By a similar argument to Case 1. and 2., it follows
\begin{eqnarray*}
S_{m,j}T&=&S_{m,j}[j_1,\cdots,j_{i-1},j,j_{i+1},\cdots, j_{i'-1},\overline{j+1},j_{i'+1},\cdots,j_k]^{\rm X}_s\\
&=&[j_1,\cdots,j_{i-1},j,j_{i+1},\cdots, j_{i'-1},\overline{j+1},j_{i'+1},\cdots,j_k]^{\rm X}_s - 2\beta_{m,j} \\
&=&[j_1,\cdots,j_{i-1},j+1,j_{i+1},\cdots, j_{i'-1},\overline{j},j_{i'+1},\cdots,j_k]^{\rm X}_s\in {\rm Tab}_{{\rm X},\iota}.
\end{eqnarray*}

\vspace{2mm}

\nd
\underline{Case 4. (3) holds}

\vspace{2mm} 

We obtain $c_{m,n}=1$ (if $\mathfrak{g}$ is of type B), $c_{m,n}=2$ (if $\mathfrak{g}$ is of type C)
and $m=s+k-i+P(n)$. Taking Lemma \ref{box-lem} into account, we have
\begin{eqnarray*}
S_{m,n}T&=&S_{m,n}[j_1,\cdots,j_{i-1},j_i,j_{i+1},\cdots,j_k]^{\rm X}_s\\
&=&
\fbox{$j_1$}^{\rm X}_{s+k-1}+\cdots+\fbox{$j_{i}$}^{\rm X}_{s+k-i}+\cdots +\fbox{$j_k$}^{\rm X}_{s}-c_{m,n}\beta_{m,n}\\
&=&
\fbox{$j_1$}^{\rm X}_{s+k-1}+\cdots+\fbox{$n$}^{\rm X}_{s+k-i}+\cdots +\fbox{$j_k$}^{\rm X}_{s}-c_{m,n}\beta_{m,n}\\
&=&
\fbox{$j_1$}^{\rm X}_{s+k-1}+\cdots+\fbox{$\ovl{n}$}^{\rm X}_{s+k-i}+\cdots +\fbox{$j_k$}^{\rm X}_{s}=
[j_1,\cdots,j_{i-1},\ovl{n},j_{i+1},\cdots,j_k]^{\rm X}_s \in {\rm Tab}_{{\rm X},\iota}.
\end{eqnarray*}

\nd
Next, we suppose the coefficient $c_{m,j}$ of $x_{m,j}$ in $T$ is negative.
Under this assumption, at least one of (4), (5), (6) holds. 
Because if (6) holds then neither (4) nor (5) hold, we need to consider the following four cases:

\vspace{2mm}

\nd
\underline{Case 5. (4) holds and (2), (5) do not hold}

\vspace{2mm}

In this case, by a similar argument to that appeared in Case 1., we obtain $c_{m,j}=-1$ and $m=s+k-i+1+P(j)$.
By (\ref{BC-rev2}), one can verify $m\geq2$.
Thus, by Lemma \ref{box-lem},
\begin{eqnarray*}
S_{m,j}T&=&S_{m,j}[j_1,\cdots,j_{i-1},j_i,j_{i+1},\cdots,j_k]^{\rm X}_s\\
&=& 
[j_1,\cdots,j_{i-1},j_i,j_{i+1},\cdots,j_k]^{\rm X}_s + \beta_{m-1,j}\\
&=&
\fbox{$j_1$}^{\rm X}_{s+k-1}+\cdots+\fbox{$j_i$}^{\rm X}_{s+k-i}+\cdots +\fbox{$j_k$}^{\rm X}_{s}+\beta_{m-1,j}\\
&=&
\fbox{$j_1$}^{\rm X}_{s+k-1}+\cdots+\fbox{$j+1$}^{\rm X}_{s+k-i}+\cdots +\fbox{$j_k$}^{\rm X}_{s}+\beta_{m-1,j}\\
&=&
\fbox{$j_1$}^{\rm X}_{s+k-1}+\cdots+\fbox{$j$}^{\rm X}_{s+k-i}+\cdots +\fbox{$j_k$}^{\rm X}_{s}=
[j_1,\cdots,j_{i-1},j,j_{i+1},\cdots,j_k]^{\rm X}_s \in {\rm Tab}_{{\rm X},\iota}.
\end{eqnarray*}

\vspace{2mm}

\nd
\underline{Case 6. (5) holds and (1), (4) do not hold}

\vspace{2mm}

In this case,
we obtain $c_{m,j}=-1$ and $m=s+k-i'+n-j+1+P(j)$.
By (\ref{BCm1}),
it is easy to see $m\geq2$.
Hence, by Lemma \ref{box-lem},
\begin{eqnarray*}
S_{m,j}T&=&S_{m,j}[j_1,\cdots,j_{i'-1},j_{i'},j_{i'+1},\cdots,j_k]^{\rm X}_s\\
&=& 
[j_1,\cdots,j_{i'-1},j_{i'},j_{i'+1},\cdots,j_k]^{\rm X}_s + \beta_{m-1,j}\\
&=&
\fbox{$j_1$}^{\rm X}_{s+k-1}+\cdots+\fbox{$j_{i'}$}^{\rm X}_{s+k-i'}+\cdots +\fbox{$j_k$}^{\rm X}_{s}+\beta_{m-1,j}\\
&=&
\fbox{$j_1$}^{\rm X}_{s+k-1}+\cdots+\fbox{$\ovl{j}$}^{\rm X}_{s+k-i'}+\cdots +\fbox{$j_k$}^{\rm X}_{s}+\beta_{m-1,j}\\
&=&
\fbox{$j_1$}^{\rm X}_{s+k-1}+\cdots+\fbox{$\ovl{j+1}$}^{\rm X}_{s+k-i'}+\cdots +\fbox{$j_k$}^{\rm X}_{s}=
[j_1,\cdots,j_{i'-1},\ovl{j+1},j_{i'+1},\cdots,j_k]^{\rm X}_s \in {\rm Tab}_{{\rm X},\iota}.
\end{eqnarray*}

\vspace{2mm}

\nd
\underline{Case 7. Both (4) and (5) hold}

\vspace{2mm}

In this case, we obtain $c_{m,j}=-2$. As in Case 5., 6., we have $m\geq2$ and
\begin{eqnarray*}
S_{m,j}T&=&S_{m,j}[j_1,\cdots,j_{i-1},j+1,j_{i+1},\cdots, j_{i'-1},\overline{j},j_{i'+1},\cdots,j_k]^{\rm X}_s\\
&=&[j_1,\cdots,j_{i-1},j+1,j_{i+1},\cdots, j_{i'-1},\overline{j},j_{i'+1},\cdots,j_k]^{\rm X}_s + 2\beta_{m-1,j} \\
&=&[j_1,\cdots,j_{i-1},j,j_{i+1},\cdots, j_{i'-1},\overline{j+1},j_{i'+1},\cdots,j_k]^{\rm X}_s\in {\rm Tab}_{{\rm X},\iota}.
\end{eqnarray*}

\vspace{2mm}

\nd
\underline{Case 8. (6) holds}

\vspace{2mm}

In this case, we obtain $c_{m,n}=-1$ (if $\mathfrak{g}$ is of type B), $c_{m,n}=-2$ (if $\mathfrak{g}$ is of type C)
and $m=s+k-i+1+P(n)$.
It follows by (\ref{BCm1}) that $m\geq2$.
Using Lemma \ref{box-lem}, we get
\begin{eqnarray*}
S_{m,n}T&=&S_{m,n}[j_1,\cdots,j_{i-1},j_i,j_{i+1},\cdots,j_k]^{\rm X}_s\\
&=& 
[j_1,\cdots,j_{i-1},j_i,j_{i+1},\cdots,j_k]^{\rm X}_s -c_{m,n}\beta_{m-1,n}\\
&=&
\fbox{$j_1$}^{\rm X}_{s+k-1}+\cdots+\fbox{$j_i$}^{\rm X}_{s+k-i}+\cdots +\fbox{$j_k$}^{\rm X}_{s}-c_{m,n}\beta_{m-1,n}\\
&=&
\fbox{$j_1$}^{\rm X}_{s+k-1}+\cdots+\fbox{$\ovl{n}$}^{\rm X}_{s+k-i}+\cdots +\fbox{$j_k$}^{\rm X}_{s}-c_{m,n}\beta_{m-1,n}\\
&=&
\fbox{$j_1$}^{\rm X}_{s+k-1}+\cdots+\fbox{$n$}^{\rm X}_{s+k-i}+\cdots +\fbox{$j_k$}^{\rm X}_{s}=
[j_1,\cdots,j_{i-1},n,j_{i+1},\cdots,j_k]^{\rm X}_s \in {\rm Tab}_{{\rm X},\iota}.
\end{eqnarray*}
The definition of $S_{m,j}$ means that if the coefficient $c_{m,j}$ of $x_{m,j}$ in $T$ is $0$ then
$S_{m,j}T=T\in {\rm Tab}_{{\rm X},\iota}$. Therefore, we obtain our claim.

\vspace{2mm}

\nd
(ii) 
We take an element $T=[j_1,\cdots,j_n]^{\rm B}_s\in {\rm Tab}_{{\rm B},\iota}$ and
let us consider the action of $S_{m,j}$ ($m\in\mathbb{Z}_{\geq1}$, $j\in I$).
If $j\in[1,n-1]$ then
by Definition \ref{tab-def} (ii) $(*)^{\rm B}_n$, one knows
$j_{\xi}=j$ and $j_{\xi+1}= j+1$ (resp. $j_{\xi}=\ovl{j+1}$ and $j_{\xi+1}= \ovl{j}$) for some $\xi\in[1,n-1]$ if and only if
$\ovl{j}$, $\ovl{j+1}\notin\{j_1,\cdots,j_n\}$ (resp. $j$, $j+1\notin\{j_1,\cdots,j_n\}$).
In these cases, 
by the same calculation as in (\ref{A-pr02}) and (\ref{BC-pr1}),
we see that the coefficient of $x_{r,j}$ $(r\in\mathbb{Z}_{\geq1})$ is $0$. In particular, $S_{m,j}T=T$ follows.

The condition (1) is equivalent to that
$j_i=j<n$, $j_{i+1}\neq j+1$ and $m=s+k-i+P(j)$.
In this case, we have $j_{i'}=\ovl{j+1}$ and $j_{i'+1}\neq\ovl{j}$ for some $i'$ ($i'\in[1,n-1]$) by 
Definition \ref{tab-def} (ii) $(*)^{\rm B}_n$.
We also get $\{|j_{i+1}|,|j_{i+2}|,\cdots,|j_{i'-1}|\}=\{j+2,j+3,\cdots,n\}$,
which yields $i'-1-i=n-j-1$ so that $m=s+k-i'+n-j+P(j)$.
Thus, we see that the condition (1) is equivalent to (2).
Similarly, we also see that
the condition (4) is equivalent to (5).

Accordingly, if the coefficient $c_{m,j}$ of $x_{m,j}$ in $T$ $(j\in I)$ is positive then
either Case 3. or 4. in the proof of (i) happens. 
If $c_{m,j}$ is negative then either Case 7. or 8. in the proof of (i) happens. 
By the same argument as in (i), we obtain our claim (ii).

\vspace{2mm}

\nd
(iii) 
We take an element $T=[\ovl{n+1},j_2,\cdots,j_k]^{\rm C}_s\in {\rm Tab}_{{\rm C},\iota}$ and
let us consider the action of $S_{m,j}$ ($m\in\mathbb{Z}_{\geq1}$, $j\in I$).
It is easy to see $\ovl{n}\leq j_2<\cdots<j_k\leq\ovl{1}$ by the definition of ${\rm Tab}_{{\rm C},\iota}$ (Definition \ref{tab-def}).
Clearly, if $j<n$ and $\ovl{j+1}$, $\ovl{j}\notin\{j_2,\cdots,j_k\}$ then $S_{m,j}T=T$. 
Thus, we can explicitly write $T$ as
\begin{eqnarray}
T&=&\fbox{$\ovl{n+1}$}^{\rm C}_{s+k-1} + \sum^{k}_{i=2} \fbox{$j_i$}^{\rm C}_{s+k-i}\nonumber\\
&=& x_{s+k-1+P(n),n} + \sum^{k}_{i=2} (x_{s+k-i+P(|j_i|-1)+n-|j_i|+1,|j_i|-1}-c(|j_i|)x_{s+k-i+P(|j_i|)+n-|j_i|+1,|j_i|}),\qquad \ \ \label{T-ex}
\end{eqnarray}
where $c(n)=2$ and $c(t)=1$ for $t\in[1,n-1]$. Using $s\geq 1-P(n)$ and
a similar calculation to (\ref{BCm1}), we see that
\begin{equation}\label{BC-rev3}
s+k-1+P(n)\geq1,\ \ s+k-i+P(|j_i|-1)+n-|j_i|+1\geq1,\ \ s+k-i+P(|j_i|)+n-|j_i|+1\geq2.
\end{equation}
Note that it follows from (\ref{BC-pr1}) that for $j\in[1,n-1]$,
if $j_{\xi}=\ovl{j+1}$, $j_{\xi+1}=\ovl{j}$ with some $\xi\in[2,k]$
then the coefficient of $x_{l,j}$ in $T$ is $0$ for all $l\in\mathbb{Z}_{\geq1}$.

As seen in the proof of (i), the condition (2) (resp. (5)) is equivalent to that
$\ovl{n}\leq j_{i'}=\ovl{j+1}$ (resp. $\ovl{n}<j_{i'}=\ovl{j}$) and $x_{m,j}$
is a summand of
$\fbox{$j_{i'}$}^{\rm C}_{s+k-i'}+
\fbox{$j_{i'+1}$}^{\rm C}_{s+k-i'-1}$ (resp. $\fbox{$j_{i'}$}^{\rm C}_{s+k-i'}+
\fbox{$j_{i'-1}$}^{\rm C}_{s+k-i'+1}$)
with the coefficient $1$ (resp. $-1$) for some $i'\in[2,k]$.

Note that $\fbox{$\ovl{n+1}$}_{s+k-1}^{\rm C}=x_{s+k-1+P(n),n}$ and
\begin{eqnarray}
\fbox{$\ovl{n+1}$}_{s+k-1}^{\rm C}+\fbox{$\ovl{n}$}_{s+k-2}^{\rm C}
&=&x_{s+k-1+P(n),n}
+x_{s+k-1+P(n-1),n-1}-2x_{s+k-1+P(n),n}\nonumber\\
&=&x_{s+k-1+P(n-1),n-1}-x_{s+k-1+P(n),n}, \label{BC-pr2}
\end{eqnarray}
Thus for $T=[\ovl{n+1},j_2,\cdots,j_k]^{\rm C}_s$,
the triple $(T,j,m)$ satisfies the condition (6) (resp. (7))
if and only if $j=n$, $j_2=\ovl{n}$ (resp. $j_2\neq\ovl{n}$) and
$x_{m,n}$ is a summand in
$\fbox{$\ovl{n+1}$}_{s+k-1}^{\rm C}+\fbox{$j_2$}_{s+k-2}^{\rm C}$ so that in $T$
with the coefficient $-1$ (resp. $1$).

Hence, if the coefficient $c_{m,j}$ of $x_{m,j}$ in $T$ is non-zero then
one of (2), (5), (6), (7) holds. In the case (2) or (5), by the same way as in Case 2., 6. of the proof of (i),
we can prove 
\[
S_{m,j}T=
\begin{cases} 
[\ovl{n+1},j_2,\cdots,j_{i'-1},\ovl{j},j_{i'+1},\cdots,j_k]_s^{\rm C} & {\rm if\ (2)\ holds}, \\
[\ovl{n+1},j_2,\cdots,j_{i'-1},\ovl{j+1},j_{i'+1},\cdots,j_k]_s^{\rm C} & {\rm if\ (5)\ holds}.
\end{cases}
\]
In the case (6), it holds $i=2\leq k$, $m=s+k-1+P(n)$ and the coefficient of $x_{m,n}$ in $T$ is $-1$.
By $s\geq 1-P(n)$, it is easily shown that $m\geq2$.
Using Lemma \ref{box-lem}, we get
\begin{eqnarray*}
S_{m,n}T&=&S_{m,n}[\ovl{n+1},\ovl{n},j_3,\cdots,j_k]^{\rm C}_s\\
&=& 
[\ovl{n+1},\ovl{n},j_3,\cdots,j_k]^{\rm C}_s +\beta_{m-1,n}\\
&=&
\fbox{$\ovl{n+1}$}^{\rm C}_{s+k-1}+\fbox{$\ovl{n}$}^{\rm C}_{s+k-2}+\fbox{$j_3$}^{\rm C}_{s+k-3}+\cdots +\fbox{$j_k$}^{\rm C}_{s}+\beta_{m-1,n}\\
&=&
\fbox{$\ovl{n+1}$}^{\rm C}_{s+k-2}+\fbox{$j_3$}^{\rm C}_{s+k-3}+\cdots +\fbox{$j_k$}^{\rm C}_{s}=
[\ovl{n+1},j_3,\cdots,j_k]^{\rm C}_s \in {\rm Tab}_{{\rm C},\iota}.
\end{eqnarray*}
In the case (7), using Lemma \ref{box-lem}, we can prove $S_{m,n}T=[\ovl{n+1},\ovl{n},j_2,\cdots,j_k]^{\rm C}_s\in {\rm Tab}_{{\rm C},\iota}$.

Finally, by the definition of $S_{m,j}$, if the coefficient $c_{m,j}$ of $x_{m,j}$ in $T$ is $0$ then
$S_{m,j}T=T\in {\rm Tab}_{{\rm C},\iota}$. Therefore, we obtain our claim. \qed

\subsection{Stability of ${\rm Tab}_{{\rm D},\iota}$}

\begin{prop}\label{closednessD}
\begin{enumerate}
\item 
For each $T=[j_1,\cdots,j_k]_s^{\rm D}\in {\rm Tab}_{{\rm D},\iota}$, $j\in I$ and
$m\in\mathbb{Z}_{\geq1}$, we consider the following conditions for the triple $(T,j,m)$:
\begin{enumerate}
\item[(1)] $j<n$ and there exists $i\in [1,k]$ such that $j_i=j$, $j_{i+1}\neq j+1$ and $m=s+k-i+P(j)$,
\item[(2)] $j<n$ and there exists $i'\in [1,k]$ such that $j_{i'}=\overline{j+1}$, $j_{i'+1}\neq \overline{j}$, $n$
and $m=s+k-i'+n-j-1+P(j)$,
\item[(3)] $j<n$ and there exists $i\in [1,k]$ such that $j_{i}= j+1$, $j_{i-1}\neq j$, $\ovl{n}$ and $m=s+k-i+1+P(j)$,
\item[(4)] $j<n$ and there exists $i'\in [1,k]$ such that $j_{i'}= \overline{j}$, $j_{i'-1}\neq\ovl{j+1}$
and $m=s+k-i'+n-j+P(j)$.
\end{enumerate}
If $k\in [1,n-2]$, $j<n$ and $j_1\neq \ovl{n+1}$ then
\begin{eqnarray*}
\hspace{-15mm}& &S_{m,j}T= \\
\hspace{-15mm}& &\begin{cases}
[j_1,\cdots,j_{i-1},j+1,j_{i+1},\cdots,j_k]_s^{\rm D} & {\rm if\ (1)\ holds\ and\ (2),\ (4)\ do\ not\ hold}, \\
[j_1,\cdots,j_{i'-1},\ovl{j},j_{i'+1},\cdots,j_k]_s^{\rm D} & {\rm if\ (2)\ holds\ and\ (1),\ (3)\ do\ not\ hold}, \\
[j_1,\cdots,j_{i-1},j+1,j_{i+1},\cdots,j_{i'-1},\ovl{j},j_{i'+1},\cdots,j_k]_s^{\rm D} & {\rm if\ (1)\ and\ (2)\ hold}, \\
[j_1,\cdots,j_{i-1},j,j_{i+1},\cdots,j_k]_s^{\rm D} & {\rm if\ (3)\ holds\ and\ (2),\ (4)\ do\ not\ hold}, \ \ \\
[j_1,\cdots,j_{i'-1},\ovl{j+1},j_{i'+1},\cdots,j_k]_s^{\rm D} & {\rm if\ (4)\ holds\ and\ (1),\ (3)\ do\ not\ hold}, \\
[j_1,\cdots,j_{i-1},j,j_{i+1},\cdots,j_{i'-1},\ovl{j+1},j_{i'+1},\cdots,j_k]_s^{\rm D} & {\rm if\ (3)\ and\ (4)\ hold}, \\
T & {\rm otherwise}.
\end{cases}
\end{eqnarray*}
\item We consider the following conditions for the pair $(T,m)$:
\begin{enumerate}
\item[(5)] there exists $i\in [1,k]$ such that $j_i=n-1$, $j_{i+1}\neq \overline{n}$, $\ovl{n-1}$ and $m=s+k-i+P(n)$,
\item[(6)] there exists $i\in [1,k]$ such that $j_i=n$, $j_{i+1}\neq \ovl{n},\ \overline{n-1}$ and $m=s+k-i+P(n)$,
\item[(7)] there exists $i\in [1,k]$ such that $j_{i}= \overline{n}$, $j_{i-1}\neq n-1$, $n$ and $m=s+k-i+1+P(n)$,
\item[(8)] there exists $i\in [1,k]$ such that $j_{i}= \overline{n-1}$, $j_{i-1}\neq n-1$, $n$ and $m=s+k-i+1+P(n)$.
\end{enumerate}
If $k\in [1,n-2]$ and $j_1\neq \ovl{n+1}$ then
\begin{eqnarray*}
S_{m,n}T&=&\begin{cases}
[j_1,\cdots,j_{i-1},\ovl{n},j_{i+1}\cdots,,j_k]_s^{\rm D} & {\rm if\ (5)\ holds}, \\
[j_1,\cdots,j_{i-1},\ovl{n-1},j_{i+1}\cdots,,j_k]_s^{\rm D} & {\rm if\ (6)\ holds}, \\
[j_1,\cdots,j_{i-1},n-1,j_{i+1}\cdots,,j_k]_s^{\rm D} & {\rm if\ (7)\ holds}, \\
[j_1,\cdots,j_{i-1},n,j_{i+1}\cdots,,j_k]_s^{\rm D} & {\rm if\ (8)\ holds}, \\
T & {\rm otherwise}.
\end{cases}
\end{eqnarray*}
\item We consider the conditions (2), (4) in (i) and the following conditions for the triple $(T,j,m)$:
\begin{enumerate}
\item[(9)] $j=n$, $j_1=\ovl{n+1}$, $j_2\neq\ovl{n}$, $\ovl{n-1}$ and $m=s+k-1+P(n)$,
\item[(10)] $j=n$, $j_1=\ovl{n+1}$, $j_2=\ovl{n}$, $j_3=\ovl{n-1}$ and $m=s+k-2+P(n)$.
\end{enumerate}
For each $T=[\ovl{n+1},j_2,j_3,\cdots,j_k]_s^{\rm D}\in {\rm Tab}_{{\rm D},\iota}$, $j\in I$ and
$m\in\mathbb{Z}_{\geq1}$, we have
\[S_{m,j}T=
\begin{cases}
[\ovl{n+1},j_2,\cdots,j_{i'-1},\ovl{j},j_{i'+1},\cdots,j_k]_s^{\rm D} & {\rm if\ (2)\ holds}, \\
[\ovl{n+1},j_2,\cdots,j_{i'-1},\ovl{j+1},j_{i'+1},\cdots,j_k]_s^{\rm D} & {\rm if\ (4)\ holds}, \\
[\ovl{n+1},\ovl{n},\ovl{n-1},j_2,\cdots,j_k]_s^{\rm D} & {\rm if\ (9)\ holds}, \\
[\ovl{n+1},j_4,\cdots,j_k]_s^{\rm D} & {\rm if\ (10)\ holds}, \\
T & {\rm otherwise}.
\end{cases}
\]
\end{enumerate}

In particular, the set ${\rm Tab}_{{\rm D},\iota}$ is closed under the action of $S_{m,j}$.
\end{prop}

\nd
{\it Proof.}

\nd
(i)
For $T=[j_1,\cdots,j_k]^{\rm D}_s\in {\rm Tab}_{{\rm D},\iota}$ with $k\in[1,n-2]$ and $j_1\neq \ovl{n+1}$,
let us consider the action of $S_{m,j}$ $(m\in\mathbb{Z}_{\geq1},\ j\in I)$.
Recall that 
$T=[j_1,\cdots,j_k]^{\rm D}_s=\sum_{i=1}^{k}\fbox{$j_i$}^{\rm D}_{s+k-i}$, and
by Definition \ref{box-def}, we obtain
\begin{equation}\label{D-box}
\fbox{$j_i$}^{\rm D}_{s+k-i}=
\begin{cases}
x_{s+k-i+P(j_i),j_i}-x_{s+k-i+P(j_i-1)+1,j_i-1} & {\rm if}\ j_i\in[1,n-2]\cup\{n\},\\
x_{s+k-i+P(n-1),n-1}+x_{s+k-i+P(n),n}-x_{s+k-i+P(n-2)+1,n-2} & {\rm if}\ j_i=n-1,\\
x_{s+k-i+P(n-1),n-1}-x_{s+k-i+P(n)+1,n} & {\rm if}\ j_i= \ovl{n}, \\ 
x_{s+k-i+P(n-2)+1,n-2}-x_{s+k-i+P(n-1)+1,n-1}
-x_{s+k-i+P(n)+1,n} & {\rm if}\ j_i= \ovl{n-1},\\
x_{s+k-i+P(|j_i|-1)+n-|j_i|,|j_i|-1}-x_{s+k-i+P(|j_i|)+n-|j_i|,|j_i|} & {\rm if}\ j_i\geq \ovl{n-2}.
\end{cases}
\end{equation}
Since we know $1-P(k)\leq s$ and $P(k)\leq P(n-1)$, $P(n)$, one can check as in (\ref{A-m1}), (\ref{BCm1})
that all the left indices appearing in (\ref{D-box}) are positive.
In the case $j_i\in[1,n]$, just as in (\ref{BC-rev2}), we see that if $j_i>j_{i-1}+1$ then
\begin{equation}\label{D-rev1}
s+k-i+P(j_i-1)+1\geq2.
\end{equation}
We can also check
\begin{equation}\label{D-rev2}
s+k-i+P(n)+1\geq2,\ \ s+k-i+P(n-1)+1\geq2,\ \ s+k-i+P(|j_i|)+n-|j_i|\geq2.
\end{equation}
By Definition \ref{box-def}, it is easy to see that for $l\in\mathbb{Z}_{\geq 1-P(n-1)}$,
\begin{eqnarray}
\fbox{$n-1$}^{\rm D}_{l+1}+\fbox{$\ovl{n}$}^{\rm D}_{l}
&=&
x_{l+1+P(n-1),n-1}+x_{l+1+P(n),n}-x_{l+P(n-2)+2,n-2}
+x_{l+P(n-1),n-1}-x_{l+P(n)+1,n}\nonumber\\
&=&x_{l+1+P(n-1),n-1}-x_{l+P(n-2)+2,n-2}
+x_{l+P(n-1),n-1},\label{D-pr-a}
\end{eqnarray}
\begin{eqnarray}
\fbox{$n$}^{\rm D}_{l+1}+\fbox{$\ovl{n-1}$}^{\rm D}_{l}
&=&
x_{l+1+P(n),n}-x_{l+P(n-1)+2,n-1}
+x_{l+P(n-2)+1,n-2}-x_{l+P(n-1)+1,n-1}-x_{l+P(n)+1,n} \nonumber\\
&=&
x_{l+P(n-2)+1,n-2}-x_{l+P(n-1)+1,n-1}
-x_{l+P(n-1)+2,n-1},\label{D-pr-c}
\end{eqnarray}
\begin{eqnarray}
\fbox{$n$}^{\rm D}_{l+1}+\fbox{$\ovl{n}$}^{\rm D}_{l}
&=&
x_{l+1+P(n),n}-x_{l+P(n-1)+2,n-1}
+x_{l+P(n-1),n-1}-x_{l+P(n)+1,n}\nonumber\\
&=&-x_{l+P(n-1)+2,n-1}
+x_{l+P(n-1),n-1}.\label{D-pr-d}
\end{eqnarray}
For $l\in\mathbb{Z}_{\geq 1-P(n-2)}$,
\begin{eqnarray}
\fbox{$n-1$}^{\rm D}_{l+1}+\fbox{$\ovl{n-1}$}^{\rm D}_{l}
&=&
x_{l+1+P(n-1),n-1}+x_{l+1+P(n),n}-x_{l+P(n-2)+2,n-2}\nonumber\\
& &+x_{l+1+P(n-2),n-2}-x_{l+1+P(n-1),n-1}
-x_{l+1+P(n),n}\nonumber\\
&=&-x_{l+P(n-2)+2,n-2}
+x_{l+1+P(n-2),n-2}.\label{D-pr-b}
\end{eqnarray}
For $l\in\mathbb{Z}_{\geq 1-P(n)}$,
\begin{eqnarray}
\fbox{$\ovl{n}$}^{\rm D}_{l+1}+\fbox{$n$}^{\rm D}_{l}
&=&
x_{l+1+P(n-1),n-1}-x_{l+2+P(n),n}
+x_{l+P(n),n}-x_{l+1+P(n-1),n-1}
\nonumber\\
&=&-x_{l+2+P(n),n}
+x_{l+P(n),n}.\label{D-pr-e}
\end{eqnarray}
By a similar argument to the proof of Proposition \ref{closednessBC} (i), (\ref{D-pr-a}) and (\ref{D-pr-b}),
we see that the condition (1) in the claim is equivalent to
that $j_i=j<n$ and $x_{m,j}$
is a summand of
$\fbox{$j_i$}^{\rm D}_{s+k-i}+
\fbox{$j_{i+1}$}^{\rm D}_{s+k-i-1}$
for some $i\in [1,k]$
with the coefficient $1$ except for the cases
\begin{equation}\label{exd1}
j_i=j=n-2,\ \ j_{i+1}=\ovl{n-1},
\end{equation}
\begin{equation}\label{exd2}
j_i=j=n-1,\ \ j_{i+1}=\ovl{n-1}.
\end{equation}
In the case (\ref{exd1}) (resp. (\ref{exd2})), the coefficient of $x_{m,j}$ in $\fbox{$j_i$}^{\rm D}_{s+k-i}+
\fbox{$j_{i+1}$}^{\rm D}_{s+k-i-1}$ is $2$ (resp. $0$).

By a similar argument to the proof of Proposition \ref{closednessBC} (i) and (\ref{D-pr-a})-(\ref{D-pr-e}),
we also see that the condition (2) (resp. (3)) in the claim is equivalent to
that $j_{i'}=\ovl{j+1}\leq \ovl{n}$ (resp. $j_i=j+1\leq n$) and $x_{m,j}$
is a summand of
$\fbox{$j_{i'}$}^{\rm D}_{s+k-i'}+
\fbox{$j_{i'+1}$}^{\rm D}_{s+k-i'-1}$
(resp. $\fbox{$j_i$}^{\rm D}_{s+k-i}+
\fbox{$j_{i-1}$}^{\rm D}_{s+k-i+1}$)
with the coefficient $1$ (resp. $-1$) for some $i$, $i'\in[1,k]$.

We also see that the condition (4)
in the claim is equivalent to
that $j_{i'}=\ovl{j}> \ovl{n}$ and $x_{m,j}$
is a summand of $\fbox{$j_{i'}$}^{\rm D}_{s+k-i'}+
\fbox{$j_{i'-1}$}^{\rm D}_{s+k-i'+1}$ for some $i'\in[1,k]$
with coefficient $-1$ except for the cases
\begin{equation}\label{exd3}
j=n-2,\ j_{i'}=\ovl{n-2},\ \ j_{i'-1}=n-1,
\end{equation}
\begin{equation}\label{exd4}
j=n-1,\ j_{i'}=\ovl{n-1},\ \ j_{i'-1}=n-1.
\end{equation}
In the case (\ref{exd3}) (resp. (\ref{exd4})), the coefficient of $x_{m,j}$ in $\fbox{$j_{i'}$}^{\rm D}_{s+k-i'}+
\fbox{$j_{i'-1}$}^{\rm D}_{s+k-i'+1}$ is $-2$ (resp. $0$).

If (1) and (4) hold then $m=s+k-i+P(j)=s+k-i'+n-j+P(j)$
so that $i'=i+n-j$. In the case $j=n-1$, we have $i'=i+1$
and the coefficient of $x_{m,j}=x_{m,n-1}$ in $T$
is $0$ by (\ref{D-pr-b}). In the case $j=n-2$,
we get $i'=i+2$, $j_i=n-2$, $j_{i+1}\in\{n,\ovl{n}\}$ and $j_{i+2}=\ovl{n-2}$.
By a direct calculation, we see that the coefficient of $x_{m,j}=x_{m,n-2}$
is $0$ in $\fbox{$j_i$}^{\rm D}_{s+k-i}+
\fbox{$j_{i+1}$}^{\rm D}_{s+k-i-1}+
\fbox{$j_{i+2}$}^{\rm D}_{s+k-i-2}$. It follows from $j_1\ngeq j_2\ngeq\cdots\ngeq j_k$
that the coefficient of $x_{m,j}=x_{m,n-2}$ in $T$ is $0$.
In the case $j<n-2$, we get $i'-i>2$ and 
the coefficient of $x_{m,j}$ in $T$ is $0$ by the above argument.
Similarly, if (2) and (3) hold then the coefficient of $x_{m,j}$ in $T$ is $0$.

Now, we suppose the coefficient $c_{m,j}$ of $x_{m,j}$ in $T$ is positive.
At least one of (1), (2) holds.

\vspace{2mm}

\nd
\underline{Case 1. (1) holds and (2), (4) do not hold}

\vspace{2mm}

In this case, since (4) does not hold, the case (\ref{exd2}) does not hold.
Since (2) does not hold if (\ref{exd1}) holds then we obtain $j_{i+2}=\ovl{n-2}$.
Thus, the coefficient of $x_{m,j}$ in $T$ is $1$ by a similar calculation to (\ref{BC-pr1}).
By the same argument as in Case 1. of the proof of Proposition \ref{closednessBC}, in both the case
(\ref{exd1}) holds and the case (\ref{exd1}) does not hold,
we obtain $c_{m,j}=1$, $m=s+k-i+P(j)$ and
$S_{m,j}T=
[j_1,\cdots,j_{i-1},j+1,j_{i+1},\cdots,j_k]^{\rm D}_s \in {\rm Tab}_{{\rm D},\iota}$.

\vspace{2mm}

\nd
\underline{Case 2. (2) holds and (1), (3) do not hold}

\vspace{2mm}

In this case,
we obtain $c_{m,j}=1$, $m=s+k-i'+n-j-1+P(j)$ and
\[
S_{m,j}T=
[j_1,\cdots,j_{i'-1},\ovl{j},j_{i'+1},\cdots,j_k]^{\rm D}_s \in {\rm Tab}_{{\rm D},\iota}.
\]

\vspace{2mm}

\nd
\underline{Case 3. Both (1) and (2) hold}

\vspace{2mm}

In this case, the case (\ref{exd2}) does not hold. In both the case
(\ref{exd1}) holds and the case (\ref{exd1}) does not hold, 
we obtain $c_{m,j}=2$, $m=s+k-i+P(j)$ and $m=s+k-i'+n-j-1+P(j)$.
By a similar argument to Case 1. and 2., it follows
\[
S_{m,j}T=[j_1,\cdots,j_{i-1},j+1,j_{i+1},\cdots, j_{i'-1},\overline{j},j_{i'+1},\cdots,j_k]^{\rm D}_s\in {\rm Tab}_{{\rm D},\iota}.
\]

\nd
Next, we suppose the coefficient $c_{m,j}$ of $x_{m,j}$ in $T$ is negative.
Then at least one of (3), (4) holds.

\vspace{2mm}

\nd
\underline{Case 4. (3) holds and (2), (4) do not hold}

\vspace{2mm}

In this case, 
by a similar argument to that appeared in Case 5. of the proof of Proposition \ref{closednessBC} and Case 1.,
we obtain $c_{m,j}=-1$, $m=s+k-i+1+P(j)\geq2$ and
$S_{m,j}T=
[j_1,\cdots,j_{i-1},j,j_{i+1},\cdots,j_k]^{\rm D}_s \in {\rm Tab}_{{\rm D},\iota}$.

\vspace{2mm}

\nd
\underline{Case 5. (4) holds and (1), (3) do not hold}

\vspace{2mm}

In this case,
we obtain $c_{m,j}=-1$, $m=s+k-i'+n-j+P(j)= s+k-i'+1 + (n-1)-j+P(j)\geq s+k-i'+1+P(n-1)\geq2$ and
$S_{m,j}T=
[j_1,\cdots,j_{i'-1},\ovl{j+1},j_{i'+1},\cdots,j_k]^{\rm D}_s \in {\rm Tab}_{{\rm D},\iota}$.

\vspace{2mm}

\nd
\underline{Case 6. (3) and (4) hold}

\vspace{2mm}

In this case,
we obtain $c_{m,j}=-2$, $m\geq2$ and
\[
S_{m,j}T=
[j_1,\cdots,j_{i-1},j,j_{i+1},\cdots,j_{i'-1},\ovl{j+1},j_{i'+1},\cdots,j_k]^{\rm D}_s \in {\rm Tab}_{{\rm D},\iota}.
\]

\nd
(ii) The condition (5) (resp. (6), (7), (8)) in the claim is equivalent to
that $j_i=n-1$ (resp. $j_i=n$, $j_i=\ovl{n}$, $j_{i}=\ovl{n-1}$) and $x_{m,n}$
is a summand of
$\fbox{$j_i$}^{\rm D}_{s+k-i}+
\fbox{$j_{i+1}$}^{\rm D}_{s+k-i-1}$
(resp. $\fbox{$j_i$}^{\rm D}_{s+k-i}+
\fbox{$j_{i+1}$}^{\rm D}_{s+k-i-1}$, 
$\fbox{$j_i$}^{\rm D}_{s+k-i}+
\fbox{$j_{i-1}$}^{\rm D}_{s+k-i+1}$, $\fbox{$j_i$}^{\rm D}_{s+k-i}+
\fbox{$j_{i-1}$}^{\rm D}_{s+k-i+1}$)
with the coefficient $1$ (resp. $1$, $-1$, $-1$) for some $i\in[1,k]$.
Note that since $m$, $s$, $k$ and $n$ are fixed, two of (5)-(8) do not hold at the same time.

\vspace{3mm}

\nd
We suppose the coefficient $c_{m,n}$ of $x_{m,n}$ in $T$ is positive.

\vspace{3mm}

\nd
\underline{Case 1. (5) holds}

\vspace{2mm} 

In this case, we obtain $c_{m,n}=1$
and $m=s+k-i+P(n)$.
Taking Lemma \ref{box-lem} into account, we have
\begin{eqnarray*}
S_{m,n}T&=&S_{m,n}[j_1,\cdots,j_{i-1},j_i,j_{i+1},\cdots,j_k]^{\rm D}_s\\
&=&
\fbox{$j_1$}^{\rm D}_{s+k-1}+\cdots+\fbox{$j_{i}$}^{\rm D}_{s+k-i}+\cdots +\fbox{$j_k$}^{\rm D}_{s}-\beta_{m,n}\\
&=&
\fbox{$j_1$}^{\rm D}_{s+k-1}+\cdots+\fbox{$n-1$}^{\rm D}_{s+k-i}+\cdots +\fbox{$j_k$}^{\rm D}_{s}-\beta_{m,n}\\
&=&
\fbox{$j_1$}^{\rm D}_{s+k-1}+\cdots+\fbox{$\ovl{n}$}^{\rm D}_{s+k-i}+\cdots +\fbox{$j_k$}^{\rm D}_{s}=
[j_1,\cdots,j_{i-1},\ovl{n},j_{i+1},\cdots,j_k]^{\rm D}_s \in {\rm Tab}_{{\rm D},\iota}.
\end{eqnarray*}

\vspace{2mm}

\nd
\underline{Case 2. (6) holds}

\vspace{2mm} 

We obtain $c_{m,n}=1$
and $m=s+k-i+P(n)$. 
Taking Lemma \ref{box-lem} into account, we have
\begin{eqnarray*}
S_{m,n}T&=&S_{m,n}[j_1,\cdots,j_{i-1},j_i,j_{i+1},\cdots,j_k]^{\rm D}_s\\
&=&
\fbox{$j_1$}^{\rm D}_{s+k-1}+\cdots+\fbox{$j_{i}$}^{\rm D}_{s+k-i}+\cdots +\fbox{$j_k$}^{\rm D}_{s}-\beta_{m,n}\\
&=&
\fbox{$j_1$}^{\rm D}_{s+k-1}+\cdots+\fbox{$n$}^{\rm D}_{s+k-i}+\cdots +\fbox{$j_k$}^{\rm D}_{s}-\beta_{m,n}\\
&=&
\fbox{$j_1$}^{\rm D}_{s+k-1}+\cdots+\fbox{$\ovl{n-1}$}^{\rm D}_{s+k-i}+\cdots +\fbox{$j_k$}^{\rm D}_{s}=
[j_1,\cdots,j_{i-1},\ovl{n-1},j_{i+1},\cdots,j_k]^{\rm D}_s \in {\rm Tab}_{{\rm D},\iota}.
\end{eqnarray*}

\vspace{3mm}

\nd
Next, we suppose the coefficient $c_{m,n}$ of $x_{m,n}$ in $T$ is negative.

\nd
\underline{Case 3.  (7) holds}

\vspace{2mm}

We obtain $c_{m,n}=-1$, $m\geq2$ and
$S_{m,n}T=
[j_1,\cdots,j_{i-1},n-1,j_{i+1},\cdots,j_k]^{\rm D}_s \in {\rm Tab}_{{\rm D},\iota}$.

\vspace{2mm}

\nd
\underline{Case 4.  (8) holds}

\vspace{2mm}

In this case,
we obtain $c_{m,n}=-1$, $m\geq2$ and
$S_{m,n}T=
[j_1,\cdots,j_{i-1},n,j_{i+1},\cdots,j_k]^{\rm D}_s \in {\rm Tab}_{{\rm D},\iota}$.

\vspace{4mm}

\nd
(iii) We take an element $T=[\ovl{n+1},j_2,\cdots,j_k]^{\rm D}_s\in {\rm Tab}_{{\rm D},\iota}$ and
let us consider the action of $S_{m,j}$ ($m\in\mathbb{Z}_{\geq1}$, $j\in I$).
It is easy to see $\ovl{n}\leq j_2<\cdots<j_k\leq\ovl{1}$ by the definition of ${\rm Tab}_{{\rm D},\iota}$ (Definition \ref{tab-def} (ii)).
Thus, we can explicitly write $T$ as
\[
T=\fbox{$\ovl{n+1}$}^{\rm D}_{s+k-1} + \sum^{k}_{i=2} \fbox{$j_i$}^{\rm D}_{s+k-i}.
\]
Recall that
\begin{equation}\label{D-box-ag}
\fbox{$j_i$}^{\rm D}_{s+k-i}=
\begin{cases}
x_{s+k-i+P(n),n} & {\rm if}\ j_i=\ovl{n+1}\\
x_{s+k-i+P(n-1),n-1}-x_{s+k-i+P(n)+1,n} & {\rm if}\ j_i= \ovl{n}, \\ 
x_{s+k-i+P(n-2)+1,n-2}-x_{s+k-i+P(n-1)+1,n-1}
-x_{s+k-i+P(n)+1,n} & {\rm if}\ j_i= \ovl{n-1},\\
x_{s+k-i+P(|j_i|-1)+n-|j_i|,|j_i|-1}-x_{s+k-i+P(|j_i|)+n-|j_i|,|j_i|} & {\rm if}\ j_i\geq \ovl{n-2}.
\end{cases}
\end{equation}
In the case $j_i\geq \ovl{n-2}$, it is easy to see
\begin{eqnarray*}
& & s+k-i+P(|j_i|-1)+n-|j_i|\\
& &=s+k-i+P(|j_i|-1)+(n-1)-|j_i|+1\geq s+k-i+P(n-1),\ s+k-i+P(n).
\end{eqnarray*}
Since $s\geq 1-P(n-1)$ or $s\geq 1-P(n)$ hold by the definition of ${\rm Tab}_{{\rm D},\iota}$, the left indices
$s+k-i+P(|j_i|-1)+n-|j_i|$ appearing in (\ref{D-box-ag}) are positive. Similarly, one can check
\begin{equation}\label{D-rev3}
s+k-i+P(|j_i|)+n-|j_i|\geq2.
\end{equation}

In the case $j_i= \ovl{n-1}$, one can check the left index $s+k-i+P(n-2)+1$ in (\ref{D-box-ag}) is positive by $s\geq 1-P(n-1)$ or $s\geq 1-P(n)$ and
the values of $p_{n-1,n-2}$, $p_{n,n-2}$ are $1$ or $0$. 
In the case $k$ is even, it holds $s\geq 1-P(n-1)$, which implies $s+k-i+P(n-1)+1\geq2$.
If $i=k$ then $i=k=2$ and $-x_{s+k-i+P(n)+1,n}$  
is cancelled in $T=\fbox{$j_1$}^{\rm D}_{s+k-1}+\fbox{$j_2$}^{\rm D}_{s+k-2}=\fbox{$\ovl{n+1}$}^{\rm D}_{s+k-1}
+\fbox{$\ovl{n-1}$}^{\rm D}_{s+k-2}$.
If $k>i$ then $s+k-i+P(n)+1\geq 2$.

In the case $k$ is odd, it holds $s\geq 1-P(n)$, which implies $s+k-i+P(n)+1\geq2$.
If $i=k$ then $i=k=3$, $j_2=\ovl{n}$ and $-x_{s+k-i+P(n-1)+1,n-1}$ is cancelled 
in $T=\fbox{$j_1$}^{\rm D}_{s+k-1}+\fbox{$j_2$}^{\rm D}_{s+k-2}+\fbox{$j_3$}^{\rm D}_{s+k-3}
=\fbox{$\ovl{n+1}$}^{\rm D}_{s+k-1}+\fbox{$\ovl{n}$}^{\rm D}_{s+k-2}+\fbox{$\ovl{n-1}$}^{\rm D}_{s+k-3}$
 (see (\ref{D-pr-f})). If $i<k$ then $s+k-i+P(n-1)+1\geq2$.

In the case $j_i= \ovl{n}$ we have $i=2$ and  $\fbox{$j_2$}^{\rm D}_{s+k-2}=\fbox{$\ovl{n}$}^{\rm D}_{s+k-2}=x_{s+k-2+P(n-1),n-1}-x_{s+k-1+P(n),n}$.
We see that $-x_{s+k-1+P(n),n}$ is cancelled in $\fbox{$j_1$}^{\rm D}_{s+k-1}+\fbox{$j_2$}^{\rm D}_{s+k-2}=\fbox{$\ovl{n+1}$}^{\rm D}_{s+k-1}+\fbox{$\ovl{n}$}^{\rm D}_{s+k-2}$.
If $k$ is even then $s\geq 1-P(n-1)$ by the definition of ${\rm Tab}_{{\rm D},\iota}$.
Thus we obtain left index $s+k-i+P(n-1)$ is positive.
If $k$ is odd then we obtain $s\geq 1-P(n)$ and $2=i<k$.
Hence we get $s+k-i+P(n-1)$ is positive.
Similarly, in the case $j_i=\ovl{n+1}$, we see that the left indices $s+k-i+P(n)$ is positive.

Note that it follows from a similar argument to (\ref{BC-pr1}) that for $j\in[1,n-1]$,
if $j_{\xi}=\ovl{j+1}$, $j_{\xi+1}=\ovl{j}$ with some $\xi\in[2,k-1]$
then the coefficient of $x_{l,j}$ in $T$ is $0$ for all $l\in\mathbb{Z}_{\geq1}$.
As seen in the proof of (i), the condition (2) (resp. (4)) is equivalent to that
$\ovl{n}\leq j_{i'}=\ovl{j+1}$ (resp. $\ovl{n}<j_{i'}=\ovl{j}$) and $x_{m,j}$
is a summand of
$\fbox{$j_{i'}$}^{\rm D}_{s+k-i'}+
\fbox{$j_{i'+1}$}^{\rm D}_{s+k-i'-1}$ (resp. $\fbox{$j_{i'}$}^{\rm D}_{s+k-i'}+
\fbox{$j_{i'-1}$}^{\rm D}_{s+k-i'+1}$)
with the coefficient $1$ (resp. $-1$) for some $i'\in[2,k]$.

By a direct calculation, we get $\fbox{$\ovl{n+1}$}_{s+k-1}^{\rm D}=x_{s+k-1+P(n),n}$ and
\begin{eqnarray*}
\fbox{$\ovl{n+1}$}_{s+k-1}^{\rm D}+\fbox{$\ovl{n}$}_{s+k-2}^{\rm D}
&=&x_{s+k-1+P(n),n}
+x_{s+k-2+P(n-1),n-1}-x_{s+k-1+P(n),n}\\
&=&x_{s+k-2+P(n-1),n-1}, 
\end{eqnarray*}
\begin{eqnarray*}
\fbox{$\ovl{n+1}$}_{s+k-1}^{\rm D}+\fbox{$\ovl{n-1}$}_{s+k-2}^{\rm D}
&=&x_{s+k-1+P(n),n}
+x_{s+k-1+P(n-2),n-2}-x_{s+k-1+P(n-1),n-1}-x_{s+k-1+P(n),n}\\
&=&x_{s+k-1+P(n-2),n-2}-x_{s+k-1+P(n-1),n-1}, 
\end{eqnarray*}
\begin{eqnarray}
\fbox{$\ovl{n+1}$}_{s+k-1}^{\rm D}+\fbox{$\ovl{n}$}_{s+k-2}^{\rm D}+\fbox{$\ovl{n-1}$}_{s+k-3}^{\rm D}
&=&x_{s+k-1+P(n),n}
+x_{s+k-2+P(n-1),n-1}-x_{s+k-1+P(n),n}\nonumber \\
& &+x_{s+k-2+P(n-2),n-2}-x_{s+k-2+P(n-1),n-1}-x_{s+k-2+P(n),n}\nonumber\\
&=&x_{s+k-2+P(n-2),n-2}-x_{s+k-2+P(n),n}\label{D-pr-f}.
\end{eqnarray}
Thus for $T=[\ovl{n+1},j_2,\cdots,j_k]$,
the triple $(T,j,m)$ satisfies the condition (9) (resp. (10))
if and only if $j=n$ and
$x_{m,n}$ is a summand of
$\fbox{$\ovl{n+1}$}_{s+k-1}^{\rm D}+\fbox{$j_2$}_{s+k-2}^{\rm D}+\fbox{$j_3$}_{s+k-3}^{\rm D}$
with the coefficient $1$ (resp. $-1$).

Hence, if the coefficient $c_{m,j}$ of $x_{m,j}$ in $T$ is non-zero then
one of (2), (4), (9), (10) holds. In the case (2) or (4), by the same way as in Case 2., 5. of the proof of (i),
we can prove 
\[
S_{m,j}T=
\begin{cases} 
[\ovl{n+1},j_2,\cdots,j_{i'-1},\ovl{j},j_{i'+1},\cdots,j_k]_s^{\rm D} & {\rm if\ (2)\ holds}, \\
[\ovl{n+1},j_2,\cdots,j_{i'-1},\ovl{j+1},j_{i'+1},\cdots,j_k]_s^{\rm D} & {\rm if\ (4)\ holds}.
\end{cases}
\]
In the case (9), it holds $m=s+k-1+P(n)\geq1$.
The coefficient of $x_{m,n}$ in $T$ is $1$.
Using Lemma \ref{box-lem}, we get
\begin{eqnarray*}
S_{m,n}T&=&S_{m,n}[\ovl{n+1},j_2,j_3,\cdots,j_k]^{\rm D}_s\\ 
&=& 
[\ovl{n+1},j_2,j_3,\cdots,j_k]^{\rm D}_s -\beta_{m,n}\\
&=&
\fbox{$\ovl{n+1}$}^{\rm D}_{s+k-1}+\fbox{$j_2$}^{\rm D}_{s+k-2}+\fbox{$j_3$}^{\rm D}_{s+k-3}
+\cdots +\fbox{$j_k$}^{\rm D}_{s}-\beta_{m,n}\\
&=&
\fbox{$\ovl{n+1}$}^{\rm D}_{s+k+1}+
\fbox{$\ovl{n}$}^{\rm D}_{s+k}+
\fbox{$\ovl{n-1}$}^{\rm D}_{s+k-1}+
\fbox{$j_2$}^{\rm D}_{s+k-2}+\cdots +\fbox{$j_k$}^{\rm D}_{s}=
[\ovl{n+1},\ovl{n},\ovl{n-1},j_2,\cdots,j_k]^{\rm D}_s \in {\rm Tab}_{{\rm D},\iota}.
\end{eqnarray*}
In the case (10), it follows $k\geq3$ and $m=s+k-2+P(n)$.
If $k$ is odd then we obtain $s\geq 1-P(n)$ and $m\geq2$.
If $k$ is even then we get $s\geq 1-P(n-1)$, $k\geq4$, which implies $m\geq2$.
Using Lemma \ref{box-lem}, we can prove
$S_{m,n}T=[\ovl{n+1},j_4,\cdots,j_k]^{\rm D}_s\in {\rm Tab}_{{\rm D},\iota}$.

Finally, by the definition of $S_{m,j}$, if the coefficient $c_{m,j}$ of $x_{m,j}$ in $T$ is $0$ then
$S_{m,j}T=T\in {\rm Tab}_{{\rm D},\iota}$. Therefore, we obtain our claim. \qed

\section{Proof of Theorem \ref{thm2} and \ref{thm1}}\label{pr-sect}

Let $P_{\rm X}=\bigoplus_{j\in I}\mathbb{Z}\Lambda_j$ be the weight lattice of the Lie algebra $\mathfrak{g}$ of type X,
where $\Lambda_j$ is the $j$-th fundamental weight.
We define a partial order in $P_{\rm X}$ as follows: For $\lambda$, $\mu\in P_{\rm X}$,
$\lambda \geq \mu $ if and only if $\lambda - \mu \in \bigoplus_{j\in I}\mathbb{Z}_{\geq0}\alpha_j$.
Let ${\rm X}^L$ be the type of the Langlands dual Lie algebra of $\mathfrak{g}$,
that is, ${\rm A}^L={\rm A}$, ${\rm B}^L={\rm C}$, ${\rm C}^L={\rm B}$ and ${\rm D}^L={\rm D}$.  
We define $\mathbb{Z}$-submodule $(\ZZ^{\ify})^*$ of $(\QQ^{\ify})^*$ generated by
$x_{s,j}$ ($s\in\mathbb{Z}_{\geq1}$, $j\in I$),
where $x_{s,j}$ is the notation in (\ref{rewrite}).
When we treat type X case one will use the $\mathbb{Z}$-linear map
${\rm wt} : (\ZZ^{\ify})^*\rightarrow P_{{\rm X}^L}$
defined as ${\rm wt}(x_{s,j}):=\Lm_j$ for any $s\in\mathbb{Z}_{\geq1}$ and $j\in I$.
The definitions (\ref{betak}), (\ref{rewrite}) of $\beta_{s,j}$ and Definition \ref{adapt} imply that 
\begin{equation}\label{be-al}
{\rm wt}(\beta_{s,j})=\al_{j}\in P_{{\rm X}^L}.
\end{equation}

\subsection{Type A case}\label{Asubsect}

In this subsection, we write $\fbox{$j$}^{\rm A}_{s}$ as $\fbox{$j$}_{s}$
 and $[j_1,\cdots,j_k]_s^{\rm A}$ as $[j_1,\cdots,j_k]_s$.

\vspace{3mm}

\nd
{\it Proof of Theorem \ref{thm2} for type A case}

\nd
First, let us prove the inclusion $\Xi_{\iota}\subset{\rm Tab}_{{\rm A},\iota}$.
It is easy to see $x_{s,i}=\fbox{$i$}_{s-P(i)} +\fbox{$i-1$}_{s-P(i)+1}+\cdots+\fbox{$1$}_{s-P(i)+i-1}=[1,\cdots,i-1,i]_{s-P(i)}\in
 {\rm Tab}_{{\rm A},\iota}$ $(i\in I, s\in\mathbb{Z}_{\geq1})$.
Since we know ${\rm Tab}_{{\rm A},\iota}$ is closed under the operators $S_{m,j}$ $(j\in I, m\in\mathbb{Z}_{\geq1})$
(Proposition \ref{closednessA}),
it follows $\Xi_{\iota}\subset{\rm Tab}_{{\rm A},\iota}$.

Next, let us show the inclusion ${\rm Tab}_{{\rm A},\iota}\subset \Xi_{\iota}$.
For a fixed $k\in[1,n]$ and each $[j_1,\cdots,j_k]_s\in {\rm Tab}_{{\rm A},\iota}$, we show
$[j_1,\cdots,j_k]_s\in \Xi_{\iota}$
by induction on the value $j_1+\cdots+j_k$.
By $1\leq j_1<\cdots<j_k\leq n+1$, the minimal value of $j_1+\cdots+j_k$
is $1+2+\cdots+k$. In this case, we can easily check that
$j_1=1$, $j_2=2$, $\cdots$, $j_k=k$
and
$[1,2,\cdots,k]_s=x_{s+P(k),k}\in \Xi_{\iota}$ for $s\geq 1-P(k)$.

Next, we assume $j_1+\cdots+j_k>1+2+\cdots+k$, 
so that $j_i>j_{i-1}+1$ for some $i\in[1,k]$, where we set $j_0=0$.
Taking the explicit form $[j_1,\cdots,j_k]_s 
=\sum^{k}_{i=1} (x_{s+k-i+P(j_i),j_i}-x_{s+k-i+1+P(j_i-1),j_{i}-1})$ into account,
the coefficient of $x_{s+k-i+P(j_i-1),j_{i}-1}$ in $[j_1,\cdots,j_i-1,\cdots,j_k]_s$ is $1$.
Therefore, putting $m:=s+k-i+P(j_i-1)$, it follows by Lemma \ref{box-lem} (i) that
\begin{eqnarray}
& &S_{m,j_{i}-1}[j_1,\cdots,j_i-1,\cdots,j_k]_s=
[j_1,\cdots,j_i-1,\cdots,j_k]_s-\beta_{m,j_i-1} \nonumber \\
&=& \fbox{$j_1$}_{s+k-1}+\cdots+\fbox{$j_i-1$}_{s+k-i}+\cdots +\fbox{$j_k$}_{s}-\beta_{m,j_i-1}
=[j_1,\cdots,j_i,\cdots,j_k]_s.\label{A-pr1}
\end{eqnarray}
Since $j_i>j_{i-1}+1$, we obtain $[j_1,\cdots,j_i-1,\cdots,j_k]_s\in{\rm Tab}_{{\rm A},\iota}$.
By the induction assumption,
we get $[j_1,\cdots,j_i-1,\cdots,j_k]_s\in \Xi_{\iota}$, which implies
$[j_1,\cdots,j_i-1,\cdots,j_k]_s=S_{l_p}\cdots S_{l_2}S_{l_1}x_{l_{0}}$ with some $l_0$, $l_1,\ \cdots,\ l_p\in\mathbb{Z}_{\geq1}$.
In conjunction with (\ref{A-pr1}), the conclusion $[j_1,\cdots,j_i,\cdots,j_k]_s\in \Xi_{\iota}$ follows. \qed

\vspace{3mm}

\nd
{\it Proof of Theorem \ref{thm1} for type A case}

\vspace{2mm}

Each element in ${\rm Tab}_{{\rm A},\iota}$ is in the following form:
\begin{equation}\label{A-pr-03}
T=\sum^{k}_{i=1} (x_{s+k-i+P(j_i),j_i}-x_{s+k-i+1+P(j_i-1),j_{i}-1})
\end{equation}
with $1\leq j_1<\cdots<j_k\leq n+1$, $k\in I$ and $s\geq 1-P(k)$.
All we need to prove is each summand $x_{s+k-i+1+P(j_i-1),j_{i}-1}$ in (\ref{A-pr-03})
is cancelled or $s+k-i+1+P(j_i-1)\geq2$ (Remark \ref{pos-rem}). Note that $j_i\geq i$ $(1\leq i\leq k)$. 
If $j_i>i$ so that $j_i-1\geq i$ then by the same calculation as in (\ref{A-pr-0}), we have $s+k-i+1+P(j_i-1)\geq2$.

If $j_i=i$ then we have $j_1=1$, $j_2=2,\cdots,j_{i-1}=i-1$.
By (\ref{A-pr02}), if $j_i=j_{i-1}+1$ then 
$x_{s+k-i+1+P(j_i-1),j_{i}-1}$ is cancelled.\qed

\subsection{Type B and C cases}

In this subsection, we supposed $X={\rm B}$ or C.
Note that Lemma \ref{box-lem} and (\ref{be-al}) mean for $s\in\mathbb{Z}$ in Lemma \ref{box-lem},
\begin{equation}\label{wt-rel}
\begin{array}{l}
{\rm wt}(\fbox{$j$}^{\rm X}_s) > {\rm wt}(\fbox{$j+1$}^{\rm X}_s)\ (1\leq j\leq n-1),\ \ \ 
\ {\rm wt}(\fbox{$\ovl{j}$}^{\rm X}_s) > {\rm wt}(\fbox{$\ovl{j-1}$}^{\rm X}_s)\quad (2\leq j\leq n), \\
 {\rm wt}(\fbox{$n$}^{\rm X}_s) > {\rm wt}(\fbox{$\ovl{n}$}^{\rm X}_s),\ \ 
{\rm wt}(\fbox{$\ovl{n+1}$}^{\rm C}_s) > {\rm wt}(\fbox{$\ovl{n+1}$}^{\rm C}_{s+1}) + {\rm wt}(\fbox{$\ovl{n}$}^{\rm C}_{s}).
\end{array}
\end{equation}

\begin{lem}\label{BClem-1}
Let $T=[j_1,\cdots,j_k]^{\rm X}_s$ be an element of ${\rm Tab}_{{\rm X},\iota}$ and $t\in[1,n-1]$.
We suppose that 
there exists $i'\in[1,k]$ such that
\begin{equation}\label{BClem-1-as}
j_{i'}=\ovl{t+1},\qquad j_{i'+1}\neq \ovl{t}.
\end{equation}
Putting $M:=s+k-i'+n-t+P(t)$,
we have $M\in\mathbb{Z}_{\geq1}$ and the following. 
\begin{enumerate}
\item[$(1)$] 
The coefficient of $x_{M,t}$ in $T$ is $0$ or $1$ or $2$. 
\item[$(2)$] The coefficient of $x_{M,t}$ in $T$ is $0$
if and only if $j_{i'-n+t+1}=t+1$ and $j_{i'-n+t}\neq t$.
\item[$(3)$] The coefficient of $x_{M,t}$ in $T$ is $2$
if and only if $j_{i'-n+t+1}\neq t+1$ and $j_{i'-n+t}=t$. 
\end{enumerate}
Furthermore, under the assumption (\ref{BClem-1-as}), 
if the coefficient of $x_{M,t}$ in $[j_1,\cdots,j_k]^{\rm X}_s$ is $1$ then
$S_{M,t}[j_1,\cdots,j_{i'},\cdots,j_k]^{\rm X}_s=S_{M,t}[j_1,\cdots,\ovl{t+1},\cdots,j_k]^{\rm X}_s
=[j_1,\cdots,\ovl{t},\cdots,j_k]^{\rm X}_s$.
\end{lem}

\nd
{\it Proof.}

The triple $(T,t,M)$ satisfies the condition (2) of Proposition \ref{closednessBC}.
Thus the condition (5) of Proposition \ref{closednessBC} does not hold.
Note that the conditions (1) and (4) of Proposition \ref{closednessBC} do not hold  simultaneously.
Hence, one of the following three cases happens:

\vspace{2mm}

\nd
\underline{(I) $(T,t,M)$ does not satisfy (1), (4)}

\vspace{2mm}

\nd
\underline{(II) $(T,t,M)$ satisfies (4) and does not satisfy (1)}

\vspace{2mm}

\nd
\underline{(III) $(T,t,M)$ satisfies (1) and does not satisfy (4)}

\vspace{3mm}

By the argument in the proof of
Proposition \ref{closednessBC} (i), we see that the coefficient of $x_{M,t}$
is $1$ (resp. $0$, $2$) in the case (I) (resp. (II), (III)).
In the case (I),
the coefficient of $x_{M,t}$ in $T$ is $1$ and $S_{M,t}T=S_{M,t}[j_1,\cdots,j_{i'},\cdots,j_k]^{\rm X}_s
=S_{M,t}[j_1,\cdots,\ovl{t+1},\cdots,j_k]^{\rm X}_s
=[j_1,\cdots,\ovl{t},\cdots,j_k]^{\rm X}_s$ by the argument in Case 2. of the proof of
Proposition \ref{closednessBC} (i).

The coefficient of $x_{M,t}$ in $T$ is $0$ if and only if the above (II) holds,
that is, $j_i=t+1$, $j_{i-1}\neq t$ and $M=s+k-i+1+P(t)$ with some $i\in[1,k]$.
If (II) holds then combining with $M:=s+k-i'+n-t+P(t)$, we get $i=i'-n+t+1$.
Conversely, if $j_{i'-n+t+1}=t+1$ and $j_{i'-n+t}\neq t$ then
the triple $(T,i'-n+t+1,M)$ satisfies the condition (4) of Proposition \ref{closednessBC},
that is, (II) holds.

The coefficient of $x_{M,t}$ in $T$ is $2$ if and only if the above (III) holds,
that is, $j_i=t$, $j_{i+1}\neq t+1$ and $M=s+k-i+P(t)$. If (III) holds then combining with
$M=s+k-i'+n-t+P(t)$, we have $i=i'-n+t$. 
Conversely, if $j_{i'-n+t+1}\neq t+1$ and $j_{i'-n+t}=t$ then
the triple $(T,i'-n+t,M)$ satisfies the condition (1) of Proposition \ref{closednessBC}, that is,
(III) holds. \qed

\vspace{3mm}

By a similar way to the proof of Lemma \ref{BClem-1}, we can verify the following lemma.

\begin{lem}\label{BClem-11}
Let $[j_1,\cdots,j_k]^{\rm X}_s$ be an element of ${\rm Tab}_{{\rm X},\iota}$ and $t\in[1,n-1]$.
We suppose that there exists $i\in[1,k]$ such that
\begin{equation}\label{BClem-2-as}
j_i=t,\qquad j_{i+1}\neq t+1.
\end{equation}
Putting $M'=s+k-i+P(t)$, 
we have $M'\in\mathbb{Z}_{\geq1}$ and the following.
\begin{enumerate}
\item[$(1)$] 
The coefficient of $x_{M',t}$ in $[j_1,\cdots,j_k]^{\rm X}_s$ is $0$ or $1$ or $2$. 
\item[$(2)$] The coefficient of $x_{M',t}$ in $[j_1,\cdots,j_k]^{\rm X}_s$ is $0$
if and only if $j_{i+n-t+1}=\ovl{t}$ with $j_{i+n-t}\neq \ovl{t+1}$.
\item[$(3)$] The coefficient of $x_{M',t}$ in $[j_1,\cdots,j_k]^{\rm X}_s$ is $2$
if and only if $j_{i+n-t+1}\neq \ovl{t}$ with $j_{i+n-t}= \ovl{t+1}$. 
\end{enumerate}
Furthermore, under the assumption (\ref{BClem-2-as}), 
if the coefficient of $x_{M',t}$ in $[j_1,\cdots,j_k]^{\rm X}_s$ is $1$ then
$S_{M',t}[j_1,\cdots,j_{i},\cdots,j_k]^{\rm X}_s=S_{M',l}[j_1,\cdots,t,\cdots,j_k]^{\rm X}_s
=[j_1,\cdots,t+1,\cdots,j_k]^{\rm X}_s$.
\end{lem}

\begin{defn}
In the case Lemma \ref{BClem-1} (2) or \ref{BClem-11} (2), we say
$[j_1,\cdots,j_k]^{\rm X}_s$ has a $t$-cancelling pair.
In the case Lemma \ref{BClem-1} (3) or \ref{BClem-11} (3),
we say $[j_1,\cdots,j_k]^{\rm X}_s$ has a $t$-double pair.
\end{defn}

\begin{rem}\label{BCrem}
Let $[j_1,\cdots,j_k]^{\rm X}_s$ be an element of ${\rm Tab}_{{\rm X},\iota}$ and $t\in[1,n-1]$.
Lemma \ref{BClem-1} and \ref{BClem-11} imply the following: 
If $j_{m'}=t$, $j_{m}=\ovl{t}$ and $m-m'\neq n-t+1$
then $T$ does not have $t$-cancelling pair.
If $j_{m'}=t$, $j_{m}=\ovl{t+1}$ and $m-m'\neq n-t$
then $[j_1,\cdots,j_k]^{\rm X}_s$ does not have $t$-double pair.
\end{rem}

\begin{lem}\label{BClem-111}
Let $T=[j_1,\cdots,j_k]^{\rm X}_s$ be an element of ${\rm Tab}_{{\rm X},\iota}$.
We suppose that there exists $i\in[1,k]$ such that
\[
j_i=n,\qquad j_{i+1}\neq \ovl{n}.
\]
Putting $M:=s+k-i+P(n)$,
we have $M\in\mathbb{Z}_{\geq1}$
and
$S_{M,n}[j_1,\cdots,j_{i},\cdots,j_k]^{\rm X}_s=S_{M,n}[j_1,\cdots,n,\cdots,j_k]^{\rm X}_s
=[j_1,\cdots,\ovl{n},\cdots,j_k]^{\rm X}_s$.
\end{lem}

\nd
{\it Proof.}

Since the triple $(T,n,M)$ satisfies the condition (3) of Proposition \ref{closednessBC},
our claim is an easy consequence of the same proposition. \qed

\begin{lem}\label{BClem-2}
Let $\mathfrak{g}$ be of type X (X$=$B or C).
For $k\in[1,n-1]$, we have
\begin{equation}\label{bck}
\{S_{l_p}\cdots S_{l_1}x_{r,k} | p\in\mathbb{Z}_{\geq0},\ r,\ l_1,\cdots,l_p\in \mathbb{Z}_{\geq1} \}
=\{ [j_1,\cdots,j_k]^{\rm X}_s\in{\rm Tab}_{{\rm X},\iota} | j_1,\cdots,j_k\in J_{\rm X} \}.
\end{equation}
\end{lem}

\nd
{\it Proof.}

\nd
Let $\Xi_{\iota,k}$
be the set of the left-hand side of (\ref{bck})
and ${\rm Tab}_{{\rm X},\iota,k}$
be the set of the right-hand side of (\ref{bck}).
First, we prove $\Xi_{\iota,k}\subset {\rm Tab}_{{\rm X},\iota,k}$.
Combining $x_{r,k}=\fbox{$k$}^{\rm X}_{r-P(k)} +\fbox{$k-1$}^{\rm X}_{r-P(k)+1}+\cdots+\fbox{$1$}^{\rm X}_{r-P(k)+k-1}=[1,\cdots,k-1,k]^{\rm X}_{r-P(k)}
\in {\rm Tab}_{{\rm X},\iota,k}$ $(k\in I, r\in\mathbb{Z}_{\geq1})$ with Proposition \ref{closednessBC} (i),
we get $\Xi_{\iota,k}\subset {\rm Tab}_{{\rm X},\iota,k}$.

Next, let us prove ${\rm Tab}_{{\rm X},\iota,k}\subset \Xi_{\iota,k}$.
For $T\in {\rm Tab}_{{\rm X},\iota,k}$, we will show $T\in \Xi_{\iota,k}$
by induction on the weight of $T$. 
The elements which have the highest weight (that is, the weight which is maximum in ${\rm wt}({\rm Tab}_{{\rm X},\iota,k})$) are
 $T=[1,2,\cdots,k]^{\rm X}_s$ $(s\geq 1-P(k))$ by (\ref{wt-rel}).
In this case, we have $T=x_{s+P(k),k}\in \Xi_{\iota,k}$.

Now we assume $T=[j_1,j_2,\cdots,j_k]^{\rm X}_s\neq [1,2,\cdots,k]^{\rm X}_s$ $(s\geq 1-P(k))$.
If $j_k\leq n$ then we can verify $T\in \Xi_{\iota,k}$
by the same way as in the proof of type A case in \ref{Asubsect}.
Thus we suppose that $1\leq j_1<\cdots<j_{l-1}\leq n$ and
$\overline{n}\leq j_{l}<\cdots<j_k\leq \overline{1}$
with some integer $l$ $(1\leq l\leq k)$. We can write
$j_l=\overline{t}$ with some $t\in I$.

\vspace{2mm}

\nd
\underline{Case 1. $j_l=\ovl{t}>\ovl{n-1}$}

\vspace{2mm}

First, let us consider the case $\ovl{n-1}<\ovl{t}$.
Let $T_0\in {\rm Tab}_{{\rm X},\iota,k}$ be the element obtained from $T$
by replacing $j_l=\overline{t}$ with $\ovl{t+1}$. 
We can verify ${\rm wt}(T_0)>{\rm wt}(T)$ by (\ref{wt-rel}).

\vspace{2mm}

\nd
\underline{Case 1-1. $T_0$ has neither $t$-cancelling nor $t$-double pair}

\vspace{2mm}

In this case, we have $S_{M,t}T_0=T$, where $M$ is the integer in Lemma \ref{BClem-1}.
Since we know ${\rm wt}(T_0)>{\rm wt}(T)$, using the induction assumption, we obtain
$T_0=S_{l_p}\cdots S_{l_1}x_{\xi,k}$ with some $p\in\mathbb{Z}_{\geq0},\ \xi,\ l_1,\cdots,l_p\in \mathbb{Z}_{\geq1}$,
which yields $T=S_{M,t}S_{l_p}\cdots S_{l_1}x_{\xi,k}\in \Xi_{\iota,k}$.

\vspace{2mm}

\nd
\underline{Case 1-2. $T_0$ has a $t$-cancelling pair}

\vspace{2mm}

In this case, $j_{l-n+t+1}=t+1$ and $j_{l-n+t}\neq t$ by Lemma \ref{BClem-1}.
Instead of $T_0$, we consider the element $T'$ obtained from $T$ by
replacing  $j_{l-n+t+1}=t+1$ with $t$.
Using Lemma \ref{BClem-11},
we see that the coefficient of $x_{M',t}$ ($M'$ is the integer in Lemma \ref{BClem-11}) in $T'$ is equal to $1$
and $S_{M',t}T'=T$. By (\ref{wt-rel}), we get
${\rm wt}(T')>{\rm wt}(T)$.
It follows from the induction assumption, we obtain
$T'=S_{l_p}\cdots S_{l_1}x_{\xi,k}$ and $T=S_{M',t}S_{l_p}\cdots S_{l_1}x_{\xi,k}\in \Xi_{\iota,k}$
with some $p\in\mathbb{Z}_{\geq0},\ \xi,\ l_1,\cdots,l_p\in \mathbb{Z}_{\geq1}$.

\vspace{2mm}

\nd
\underline{Case 1-3. $T_0$ has a $t$-double pair}

\vspace{2mm}

In this case, Lemma \ref{BClem-1} claims $j_{l-n+t}=t$
and $j_{l-n+t+1}\neq t+1$. Let us put $J:=j_{l-n+t+1}$.
Since $t<n-1$ so that $l-n+t+1<l$, we obtain $t+2\leq J\leq n$.
Let $T''$ be the element obtained from $T$ by replacing
$j_{l-n+t+1}=J$ with $J-1$. Lemma \ref{BClem-11}
shows that the coefficient of $x_{M'',J-1}$ ($M''$ is an integer) 
in $T''$ is $1$ and $T=S_{M'',J-1}T''$.
Hence, $T\in \Xi_{\iota,k}$ follows by ${\rm wt}(T'')>{\rm wt}(T)$ and the induction assumption.

\vspace{2mm}

\nd
\underline{Case 2. $j_l\in\{\ovl{n},\ovl{n-1}\}$ and $(j_{l-1},j_l)\neq(n,\ovl{n})$}

\vspace{2mm}

By the above argument, we may assume $j_l=\ovl{t}=\ovl{n-1}$ or $j_l=\ovl{t}=\ovl{n}$.
If $j_l=\ovl{n-1}$ and $j_{l-1}\neq n-1$
(resp. $j_l=\ovl{n}$ and $j_{l-1}\neq n$) then
considering the element $T_1$ obtained from $T$ by replacing
$j_l=\ovl{n-1}$ with $\ovl{n}$ (resp. $j_l=\ovl{n}$ with $n$),
we obtain $T=S_{M_1,n-1}T_1$ (resp. $T=S_{M_1,n}T_1$) with some integer $M_1$
by Lemma \ref{BClem-1} (resp. Lemma \ref{BClem-111}),
which yields $T\in \Xi_{\iota,k}$ by ${\rm wt}(T_1)>{\rm wt}(T)$ and the induction assumption.

Hence, we assume $(j_{l-1},j_l)=(n-1,\ovl{n-1})$ or $(j_{l-1},j_l)=(n,\ovl{n})$.
As a consequence of a direct calculation using Definition \ref{box-def} (ii), (iii), one get
\begin{equation}\label{BC-pr1-1}
\fbox{$n-1$}^{\rm X}_{\xi+1}+\fbox{$\ovl{n-1}$}^{\rm X}_{\xi}=\fbox{$n$}^{\rm X}_{\xi+1}+\fbox{$\ovl{n}$}^{\rm X}_{\xi}
\end{equation}
for any $\xi\in \mathbb{Z}_{\geq1-P(k)}$.
Thus, we may also assume $(j_{l-1},j_l)=(n,\ovl{n})$.

\vspace{2mm}

\nd
\underline{Case 3. $(j_{l-1},j_l)=(n,\ovl{n})$ and $l=k$}

\vspace{2mm}

Now we consider the element $T=[j_1,\cdots,j_k]^{\rm X}_s$. 
Thus the condition $l=k$ implies that in the tableau description
\[
T=
\begin{ytableau}
j_1 \\
j_2 \\
\vdots \\
j_k
\end{ytableau}^{\rm X}_{s},
\]
there are no boxes under $\fbox{$j_k$}_s=\fbox{$j_l$}_s$.

\vspace{2mm}

\nd
\underline{Case 3-1. $j_{l-2}\leq n-3$ or $l=2$}

\vspace{2mm}

Let $T_2$ be the element in ${\rm Tab}_{{\rm X},\iota,k}$ obtained from $T$ by replacing 
$(j_{l-1},j_l)=(n,\ovl{n})$ with $(n-2,\ovl{n-1})$.
Lemma \ref{BClem-11} says that
there exists an integer $M_2$ such that 
$S_{M_2,n-2}T_2=S_{M_2,n-2}[j_1,\cdots,j_{l-2},n-2,\ovl{n-1}]^{\rm X}_s
=[j_1,\cdots,j_{l-2},n-1,\ovl{n-1}]^{\rm X}_s=[j_1,\cdots,j_{l-2},n,\ovl{n}]^{\rm X}_s=T$,
where we use (\ref{BC-pr1-1}) in the third equality.
By ${\rm wt}(T_2)>{\rm wt}([j_1,\cdots,j_{l-2},n-1,\ovl{n-1}]^{\rm X}_s)=
{\rm wt}([j_1,\cdots,j_{l-2},n,\ovl{n}]^{\rm X}_s)={\rm wt}(T)$,
using the induction assumption, we get $T_{2}\in\Xi_{\iota,k}$, which means
$T\in\Xi_{\iota,k}$.

\vspace{2mm}

\nd
\underline{Case 3-2. $n-2\leq j_{l-2}\leq n-1$ and $2<l$}

\vspace{2mm}

Since we supposed $k<n$ and $k=l$,
there exists $l'\in[1,l-2]$ such that
$j_{l'-1}+1<j_{l'}$, where we set $j_0=0$.
Let $T_2'$ be the element in ${\rm Tab}_{{\rm X},\iota,k}$ obtained from $T$ by replacing 
$j_{l'}$ with $j_{l'}-1$. Because we supposed $j_l=j_k=\ovl{n}$ is the unique barred index
in $\{j_1,\cdots,j_k\}$, the element $T_2'$ has neither 
$(j_{l'}-1)$-cancelling nor $(j_{l'}-1)$-double pair.
Keeping Lemma \ref{BClem-11} in mind,
we see that $S_{M_2',j_{l'}-1}T_2'=T$ with some integer $M_2'$.
Using the induction assumption, we obtain $T_2'\in\Xi_{\iota,k}$ and
$T\in\Xi_{\iota,k}$.

\vspace{2mm}

\nd
\underline{Case 4. $(j_{l-1},j_l)=(n,\ovl{n})$, $l<k$ and $\ovl{n-1}<j_{l+1}$}

\vspace{2mm}

We put $p:=|j_{l+1}|$.

\vspace{2mm}

\nd
\underline{Case 4-1. $p=n-2$}

\vspace{2mm}

Let $T^{\dagger}$ be the element in ${\rm Tab}_{{\rm X},\iota,k}$ obtained from $T$
by replacing 
$j_{l+1}=\ovl{n-2}$ with $\ovl{n-1}$.
By Lemma \ref{BClem-1},
$T^{\dagger}$ has neither $(n-2)$-cancelling nor $(n-2)$-double pair, and
we see that
$T=S_{M^{\dagger},n-2} T^{\dagger}$ with some integer $M^{\dagger}$.
Taking ${\rm wt}(T^{\dagger})>{\rm wt}(T)$ into account
that the induction assumption means $T^{\dagger}\in\Xi_{\iota,k}$.
Therefore, it can be shown that $T\in\Xi_{\iota,k}$.

\vspace{2mm}

\nd
\underline{Case 4-2. $p=n-3$}

\vspace{2mm}

If $j_{l-2}\leq n-3$ or $l=2$ then
we can prove $T\in\Xi_{\iota,k}$
as in the Case 3-1.
If $j_{l-2}=n-2$ then the element $T_3$ obtained from
$T$ by replacing $j_{l+1}=\ovl{n-3}$ with $\ovl{n-2}$
has neither $(n-3)$-cancelling nor $(n-3)$-double pair.
Hence, using the induction assumption and Lemma \ref{BClem-1}, we obtain $T_3\in\Xi_{\iota,k}$ and
then $T\in\Xi_{\iota,k}$.

If $j_{l-2}=n-1$ then the element $T_3'$ obtained from $T$
by replacing $j_{l+1}=\ovl{n-3}$ with $\ovl{n-2}$
has neither $(n-3)$-cancelling nor $(n-3)$-double pair by Lemma \ref{BClem-1}.
It follows from the induction assumption that 
$T_3'\in\Xi_{\iota,k}$ and then
$T\in\Xi_{\iota,k}$.

\vspace{2mm}

\nd
\underline{Case 4-3. $p<n-3$}

\vspace{2mm}

Let $T_3''$ be the element 
obtained from
$T$ by replacing $j_{l+1}=\ovl{p}$ with $\ovl{p+1}$. Note that ${\rm wt}(T_3'')>{\rm wt}(T)$.
If $T_3''$ has neither $p$-cancelling nor $p$-double pair then
Lemma \ref{BClem-1} means that $T=S_{M_3'',p}T_3''$ with some integer $M_3''$.
In conjunction with the induction assumption, we get $T\in \Xi_{\iota,k}$.

If $T_3''$ has a $p$-cancelling pair then 
$j_{l-n+p+2}=p+1$ and $j_{l-n+p+1}\neq p$
(Lemma \ref{BClem-1} (2)). In this case,
we take the element $T^c_3$ obtained from $T$
by replacing $j_{l-n+p+2}=p+1$ with $p$,
which has neither $p$-cancelling nor $p$-double pair (Lemma \ref{BClem-11}).
Note that ${\rm wt}(T^c_3)>{\rm wt}(T)$.
Combining the induction assumption and Lemma \ref{BClem-11},
we get $T\in \Xi_{\iota,k}$.

If $T_3''$ has a $p$-double pair then
$j_{l-n+p+2}\neq p+1$ and $j_{l-n+p+1}= p$
(Lemma \ref{BClem-1} (3)). 
By $p<n-3$, we have $l-n+p+2\leq l-2$.
We put $J_3:=j_{l-n+p+2}$.
Taking the element $T^d_3$ obtained from $T$
by replacing $j_{l-n+p+2}=J_3$ with $J_3-1$,
we obtain ${\rm wt}(T^d_3)>{\rm wt}(T)$ and $T=S_{M^d_3,J_3-1}T^d_3$ with some integer $M^d_3$
by the same way as in the Case 1-3.
The induction assumption implies
$T\in \Xi_{\iota,k}$.

\vspace{2mm}

\nd
\underline{Case 5. $(j_{l-1},j_l)=(n,\ovl{n})$, $l<k$ and $\ovl{n-1}=j_{l+1}$} 

\vspace{2mm}

If $j_{l-2}\leq n-2$ or $l=2$ then the element $T_4$
obtained from $T$ by replacing $j_{l-1}=n$ with $n-1$
has neither 
$(n-1)$-cancelling nor $(n-1)$-double pair (Lemma \ref{BClem-11}).
It follows $T=S_{M_4,n-1}T_4$. Using the induction assumption,
we get $T_4\in \Xi_{\iota,k}$ and $T\in \Xi_{\iota,k}$.
Hence, we may assume $j_{l-2}=n-1$.

In this case, there exists an integer $p$ $(1\leq p\leq n-4)$
such that $j_{l-3}=n-2$, $j_{l-4}=n-3$, $\cdots$
, $j_{l-p-1}=n-p$ and $j_{l-p-2}\neq n-p-1$.
Let $T_4'$ be the element obtained from $T$
by replacing $j_{l-p-1}=n-p$ with $j'_{l-p-1}:=n-p-1$.

Thus, all we need to show is that
$T_4'$ has neither $(n-p-1)$-cancelling nor 
$(n-p-1)$-double pair. If we prove it
then the conclusion $T\in \Xi_{\iota,k}$
follows from the induction assumption and
Lemma \ref{BClem-11}.
We take the integer $p'$ $(1\leq p'\leq n-4)$ that $j_{l+2}=\ovl{n-2}$,
$j_{l+3}=\ovl{n-3},\cdots,j_{l+p'}=\ovl{n-p'}$
and $j_{l+p'+1}\neq \ovl{n-p'-1}$.
We must check that
\begin{enumerate}
\item[(1)] in the case $p'=p$, the pair $j'_{l-p-1}=n-p-1$ and $j_{l+p}=\ovl{n-p}$
is not a $(n-p-1)$-double pair,
\item[(2)] in the case $j_t=\ovl{n-p-1}$ or $\ovl{n-p}$ ($t$ is an integer satisfying $l+p'<t$),
the pair $j'_{l-p-1}=n-p-1$ and $j_{t}$ is neither $(n-p-1)$-cancelling nor 
$(n-p-1)$-double pair.
\end{enumerate}
Let us check (1). Since $p\geq 1$, we see
$T_4'$ does not have $(n-p-1)$-double pair (Remark \ref{BCrem}).

Finally, let us check (2).
We supposed $j_t=\ovl{n-p-1}$ or $\ovl{n-p}$
($l+p'<t$). By $t-(l-p-1)>p'+p+1\geq p+2$ and Remark \ref{BCrem},
we see that $T_4'$ has neither $(n-p-1)$-cancelling nor 
$(n-p-1)$-double pair.

The above argument says $T\in\Xi_{\iota,k}$. \qed

\begin{lem}\label{BClem-3}
\begin{enumerate}
\item In the case $\mathfrak{g}$ is of type ${\rm B}$,
\begin{equation}\label{bn}
\{S_{l_p}\cdots S_{l_1}x_{r,n} |p\in\mathbb{Z}_{\geq0},\ r,\ l_1,\cdots,l_p\in \mathbb{Z}_{\geq1} \}
=\{ [j_1,\cdots,j_n]^{\rm B}_s\in {\rm Tab}_{{\rm B},\iota} | j_1,\cdots,j_n\in J_{\rm B} \}.
\end{equation}
\item In the case $\mathfrak{g}$ is of type ${\rm C}$,
\begin{equation}\label{cn}
\{S_{l_p}\cdots S_{l_1}x_{r,n} |p\in\mathbb{Z}_{\geq0},\ r,\ l_1,\cdots,l_p\in \mathbb{Z}_{\geq1} \}
=\{ [j_1,\cdots,j_t]^{\rm C}_s\in {\rm Tab}_{{\rm C},\iota} | j_1=\overline{n+1} \}.
\end{equation}
\end{enumerate}
\end{lem}

\nd
{\it Proof.}

\nd
(i) Let $\Xi_{\iota,n}$
be the set of the left-hand side of (\ref{bn})
and ${\rm Tab}_{{\rm B},\iota,n}$
be the set of the right-hand side of (\ref{bn}).
First, let us prove $\Xi_{\iota,n}\subset {\rm Tab}_{{\rm B},\iota,n}$.
By $x_{r,n}=[1,2,\cdots,n]^{\rm B}_{r-P(n)}\in {\rm Tab}_{{\rm B},\iota,n}$ and Proposition \ref{closednessBC} (ii),
we get $\Xi_{\iota,n}\subset {\rm Tab}_{{\rm B},\iota,n}$.

Next, let us prove ${\rm Tab}_{{\rm B},\iota,n}\subset\Xi_{\iota,n}$. For each $T\in{\rm Tab}_{{\rm B},\iota,n}$,
we show $T=[j_1,\cdots,j_n]^{\rm B}_s\in \Xi_{\iota,n}$ using the induction on ${\rm wt}(T)$.
It is clear that $[1,2,\cdots,n]^{\rm B}_s=x_{s+P(n),n}\in \Xi_{\iota,n}$ for $s\geq 1-P(n)$,
so we may assume that $T\neq [1,2,\cdots,n]^{\rm B}_s$.
Then we can take $i\in[0,n-1]$ such that 
$j_1=1$, $j_2=2$, $\cdots$, $j_i=i$ and
$j_{i+1}\neq j_i+1$ (we set $j_0=0$).
By Definition \ref{tab-def} (ii) $(*)^{\rm B}_n$,
it follows $j_{i+1}\in[2,n]\cup \{\ovl{n}\}$.
We put $p:=j_{i+1}$. 
If $p\in[2,n]$ then
the definition of ${\rm Tab}_{{\rm B},\iota,n}$ implies that
$\ovl{p-1}\in\{j_1,\cdots,j_n\}$ and $\ovl{p}\notin\{j_1,\cdots,j_n\}$.
Let $T'\in {\rm Tab}_{{\rm B},\iota,n}$ be the element obtained from $T$ by
replacing $j_{i+1}=p$, $\ovl{p-1}$ with $p-1$, $\ovl{p}$, respectively.
Note that ${\rm wt}(T')>{\rm wt}(T)$.
Using Proposition \ref{closednessBC} (ii),
we see that there exists an integer $M$ such that 
$S_{M,p-1}T'=T'-2\beta_{M,p-1}=T$. In conjunction with
the induction assumption, we have $T\in{\rm Tab}_{{\rm B},\iota,n}$.
If $p=\ovl{n}$ then it follows
from the definition of ${\rm Tab}_{{\rm B},\iota,n}$ that $n\notin\{j_1,\cdots,j_n\}$.
Let $T''\in {\rm Tab}_{{\rm B},\iota,n}$ be the element obtained from $T$ by
replacing $j_{i+1}=\ovl{n}$ with $n$.
Taking an integer $M''$ properly, we get $S_{M'',n}T''=T$
by Lemma \ref{BClem-111}.
In conjunction with
the induction assumption, we have $T\in{\rm Tab}_{{\rm B},\iota,n}$.

\vspace{2mm}

\nd
(ii) Let $\Xi_{\iota,n}$
be the set of the left-hand side of (\ref{cn})
and
${\rm Tab}_{{\rm C},\iota,n}$
be the set of the right-hand side of (\ref{cn}).
First, let us prove $\Xi_{\iota,n}\subset{\rm Tab}_{{\rm C},\iota,n}$.
Since it is clear that $x_{r,n}=[\ovl{n+1}]^{\rm C}_{r-P(n)}\in{\rm Tab}_{{\rm C},\iota,n}$ for any $r\in\mathbb{Z}_{\geq1}$,
the claim $\Xi_{\iota,n}\subset{\rm Tab}_{{\rm C},\iota,n}$ follows from Proposition \ref{closednessBC}.

Next, let us show ${\rm Tab}_{{\rm C},\iota,n}\subset \Xi_{\iota,n}$.
For each $T=[j_1,\cdots,j_t]^{\rm C}_s\in{\rm Tab}_{{\rm C},\iota,n}$ with $j_1=\ovl{n+1}$, we prove $T\in \Xi_{\iota,n}$
by using the induction on ${\rm wt}(T)$.
It is clear $[\ovl{n+1}]^{\rm C}_s=x_{s+P(n),n}\in\Xi_{\iota,n}$ for $s\geq 1-P(n)$,
so we may assume ${\rm wt}(T)<{\rm wt}([\ovl{n+1}]^{\rm C}_s)$ so that $2\leq t$.
If $j_2=\ovl{n}$ then taking the element $T'=[\ovl{n+1},j_3,\cdots,j_t]^{\rm C}_s$,
we get $S_{M',n}T'=T$ with some integer $M'$ by 
Proposition \ref{closednessBC} (iii).
Because we can verify ${\rm wt}(T')>{\rm wt}(T)$,
the induction assumption means $T'\in \Xi_{\iota,n}$.
Since we know $S_{M',n}T'=T$, it follows $T'\in \Xi_{\iota,n}$.
If $j_2>\ovl{n}$ then taking the element $T''=[\ovl{n+1},\ovl{|j_2|+1},j_3,\cdots,j_t]^{\rm C}_s$,
we get $S_{M'',|j_2|}T''=S_{M'',|j_2|}[\ovl{n+1},\ovl{|j_2|+1},j_3,\cdots,j_t]^{\rm C}_s
=[\ovl{n+1},\ovl{|j_2|},j_3,\cdots,j_t]^{\rm C}_s=[\ovl{n+1},j_2,j_3,\cdots,j_t]^{\rm C}_s=T$ with some integer $M''$
by Lemma \ref{BClem-1}. By ${\rm wt}(T'')>{\rm wt}(T)$
and the induction assumption, we also get $T\in \Xi_{\iota,n}$.
Consequently, we completed the proof of ${\rm Tab}_{{\rm C},\iota,n}\subset \Xi_{\iota,n}$. \qed

\vspace{2mm}

\nd
{\it Proof of Theorem \ref{thm2} for type B,C cases}

\vspace{2mm}

We use the notation in the proofs of Lemma \ref{BClem-2}, \ref{BClem-3}.
Since
$\Xi_{\iota}=\coprod_{i\in I} \Xi_{\iota,i}$ and ${\rm Tab}_{{\rm X},\iota}=\coprod_{i\in I}{\rm Tab}_{{\rm X},\iota,i}$,
the claim in Theorem \ref{thm2}
is an easy consequence of Lemma \ref{BClem-2} and \ref{BClem-3}. \qed

\vspace{2mm}

\nd
{\it Proof of Theorem \ref{thm1} for type B,C cases}

\vspace{2mm}

\nd
\underline{In the case of type B}

\vspace{2mm}

We take any $T=[j_1,\cdots,j_k]^{\rm B}_s\in \Xi_{\iota}={\rm Tab}_{{\rm B},\iota}$ $(s\geq 1-P(k))$.
Recall that $T=[j_1,\cdots,j_k]^{\rm B}_s=\sum_{i\in [1,k]} \fbox{$j_i$}^{\rm B}_{s+k-i}$.
Taking Remark \ref{pos-rem}, (\ref{BC-box}), (\ref{BC-rev2}) and (\ref{BCm1}) into account,
we get our claim for type B case.

\vspace{2mm}

\nd
\underline{In the case of type C} 

\vspace{2mm}

Just as in the proof of type B,
we can verify
each $T=\sum_{l\in\mathbb{Z}_{\geq1}}c_lx_l\in\Xi_{\iota,i}={\rm Tab}_{{\rm C},\iota,i}$ ($i\in[1,n-1]$)
satisfies that $c_l\geq0$ if $l^{(-)}=0$.
We take an arbitrary $T=[\ovl{n+1},j_2,\cdots,j_k]^{\rm C}_s\in\Xi_{\iota,n}={\rm Tab}_{{\rm C},\iota,n}$
($k\in[1,n+1]$, $s\geq 1-P(n)$).
Using (\ref{T-ex}) and (\ref{BC-rev3}), we get our desired result for $i=n$.\qed

\subsection{Type D case}

In this subsection, we suppose $\mathfrak{g}$ is of type ${\rm D}_n$, and
the notation $\fbox{$j$}_{s}$ means $\fbox{$j$}^{\rm D}_{s}$
 and $[j_1,\cdots,j_k]_s$ means $[j_1,\cdots,j_k]_s^{\rm D}$.
Note that by Lemma \ref{box-lem} and (\ref{be-al}), for $s\in\mathbb{Z}$ in Lemma \ref{box-lem},
\begin{equation}\label{wt-rel2}
\begin{array}{l}
{\rm wt}(\fbox{$j$}_s) > {\rm wt}(\fbox{$j+1$}_s)\ (1\leq j\leq n-1),\ \ \ {\rm wt}(\fbox{$\ovl{j}$}_s) > {\rm wt}(\fbox{$\ovl{j-1}$}_s)\quad (2\leq j\leq n), \\
 {\rm wt}(\fbox{$n$}_s) > {\rm wt}(\fbox{$\ovl{n-1}$}_s),\ \  {\rm wt}(\fbox{$n-1$}_s) > {\rm wt}(\fbox{$\ovl{n}$}_s), \\
{\rm wt}(\fbox{$\ovl{n+1}$}_s) > {\rm wt}(\fbox{$\ovl{n+1}$}_{s+2})+
{\rm wt}(\fbox{$\ovl{n}$}_{s+1})+{\rm wt}(\fbox{$\ovl{n-1}$}_s).
\end{array}
\end{equation}

\begin{lem}\label{Dlem-1}
Let $T=[j_1,\cdots,j_k]_s$ be an element of ${\rm Tab}_{{\rm D},\iota}$ and $t\in[1,n-1]$.
We suppose that there exists $i\in[1,k]$ such that
\begin{equation}\label{Dlem-1-as}
j_i=\ovl{t+1},\qquad j_{i+1}\neq \ovl{t},\ n.
\end{equation}
Putting $M:=s+k-i+P(t)+n-t-1$,
we have $M\in\mathbb{Z}_{\geq1}$ and
\begin{enumerate}
\item[$(1)$] 
The coefficient of $x_{M,t}$ in $[j_1,\cdots,j_k]_s$ is $0$ or $1$ or $2$. 
\item[$(2)$] The coefficient of $x_{M,t}$ in $[j_1,\cdots,j_k]_s$ is $0$
if and only if $j_{i-n+t+2}=t+1$ with $j_{i-n+t+1}\neq t$, $\ovl{n}$.
\item[$(3)$] The coefficient of $x_{M,t}$ in $[j_1,\cdots,j_k]_s$ is $2$
if and only if $j_{i-n+t+2}\neq t+1$ with $j_{i-n+t+1}= t$. 
\end{enumerate}
Furthermore, under the assumption (\ref{Dlem-1-as}), 
if the coefficient of $x_{M,t}$ in $[j_1,\cdots,j_k]_s$ is $1$ then
$S_{M,t}[j_1,\cdots,j_{i},\cdots,j_k]_s=S_{M,t}[j_1,\cdots,\ovl{t+1},\cdots,j_k]_s=[j_1,\cdots,\ovl{t},\cdots,j_k]_s$.
\end{lem}

\nd
{\it Proof.}

The triple $(T,t,M)$ satisfies the condition (2) of Proposition \ref{closednessD}.
Thus the condition (4) of Proposition \ref{closednessD} does not hold.
Note that the conditions (1) and (3) of Proposition \ref{closednessD} do not hold  simultaneously.
Hence, one of the following three cases happens:

\vspace{2mm}

\nd
\underline{(I) $(T,t,M)$ does not satisfy (1), (3)}

\vspace{2mm}

\nd
\underline{(II) $(T,t,M)$ satisfies (3) and does not satisfy (1)}

\vspace{2mm}

\nd
\underline{(III) $(T,t,M)$ satisfies (1) and does not satisfy (3)}

\vspace{2mm}

\nd
By a similar argument to the proof of Lemma \ref{BClem-1}, we get our claim. \qed

\vspace{3mm}

\nd
We can also verify the following lemma as in Lemma \ref{BClem-1}, \ref{Dlem-1}.

\begin{lem}\label{Dlem-11}
Let $[j_1,\cdots,j_k]_s$ be an element of ${\rm Tab}_{{\rm D},\iota}$ and $t\in[1,n-1]$.
We suppose that there exists $i\in[1,k]$ such that
\begin{equation}\label{Dlem-2-as}
j_i=t,\qquad j_{i+1}\neq t+1.
\end{equation}
Putting $M'=s+k-i+P(t)$, we have $M'\in\mathbb{Z}_{\geq1}$ and
\begin{enumerate}
\item[$(1)$] 
The coefficient of $x_{M',t}$ in $[j_1,\cdots,j_k]_s$ is $0$ or $1$ or $2$. 
\item[$(2)$] The coefficient of $x_{M',t}$ in $[j_1,\cdots,j_k]_s$ is $0$
if and only if $j_{i+n-t}=\ovl{t}$ with $j_{i+n-t-1}\neq \ovl{t+1}$.
\item[$(3)$] The coefficient of $x_{M',t}$ in $[j_1,\cdots,j_k]_s$ is $2$
if and only if $j_{i+n-t}\neq \ovl{t}$, $n$ with $j_{i+n-t-1}= \ovl{t+1}$. 
\end{enumerate}
Furthermore, under the assumption (\ref{Dlem-2-as}), 
if the coefficient of $x_{M',t}$ in $[j_1,\cdots,j_k]_s$ is $1$ then
$S_{M',t}[j_1,\cdots,j_{i},\cdots,j_k]_s=S_{M',t}[j_1,\cdots,t,\cdots,j_k]_s=[j_1,\cdots,t+1,\cdots,j_k]_s$.
\end{lem}

\begin{defn}
In the case Lemma \ref{Dlem-1} (2) or \ref{Dlem-11} (2), we say
$[j_1,\cdots,j_k]_s$ has an $t$-cancelling pair.
In the case Lemma \ref{Dlem-1} (3) or \ref{Dlem-11} (3),
we say $[j_1,\cdots,j_k]_s$ has an $t$-double pair.
\end{defn}

\begin{lem}\label{Dlem-111}
Let $T=[j_1,\cdots,j_k]_s$ be an element of ${\rm Tab}_{{\rm D},\iota}$.
\begin{enumerate}
\item We suppose that there exists $i\in[1,k]$ such that
\begin{equation}\label{111a}
j_i=n-1,\qquad j_{i+1}\neq \ovl{n},\ \ovl{n-1}.
\end{equation}
Putting $M:=s+k-i+P(n)$, we have $M\in\mathbb{Z}_{\geq1}$.
The coefficient of $x_{M,n}$ in $[j_1,\cdots,j_k]_s$ is $1$
and
$S_{M,n}[j_1,\cdots,j_{i},\cdots,j_k]_s=S_{M,n}[j_1,\cdots,n-1,\cdots,j_k]_s=[j_1,\cdots,\ovl{n},\cdots,j_k]_s$.
\item We suppose that there exists $i\in[1,k]$ such that
\begin{equation}\label{111b}
j_i=n,\qquad j_{i+1}\neq \ovl{n},\ \ovl{n-1}.
\end{equation}
Putting $M:=s+k-i+P(n)$, we have $M\in\mathbb{Z}_{\geq1}$.
The coefficient of $x_{M,n}$ in $[j_1,\cdots,j_k]_s$ is $1$ and
$S_{M,n}[j_1,\cdots,j_{i},\cdots,j_k]_s=S_{M,n}[j_1,\cdots,n,\cdots,j_k]_s=[j_1,\cdots,\ovl{n-1},\cdots,j_k]_s$.
\end{enumerate}
\end{lem}

\nd
{\it Proof.}

The assumption (\ref{111a}) (resp. (\ref{111b}))
and $M=s+k-i+P(n)$ is equivalent to the condition (5) (resp. (6)) of Proposition \ref{closednessD}
for $(T,M)$. Thus, our claim is an easy consequence of Proposition \ref{closednessD}. \qed

\begin{lem}\label{Dlem-2}
For $k\in[1,n-2]$, we have
\begin{equation}\label{dk}
\{S_{l_p}\cdots S_{l_1}x_{r,k} | p\in\mathbb{Z}_{\geq0},\ r,\ l_1,\cdots,l_p\in \mathbb{Z}_{\geq1} \}
=\{ [j_1,\cdots,j_k]_s\in {\rm Tab}_{{\rm D},\iota} | j_1\neq \ovl{n+1} \}.
\end{equation}
\end{lem}

\nd
{\it Proof.}

\nd
Let $\Xi_{\iota,k}$ be the set of the left-hand side of (\ref{dk})
and ${\rm Tab}_{{\rm D},\iota,k}$ be the set of the right-hand side of (\ref{dk}).
First, we prove $\Xi_{\iota,k}\subset {\rm Tab}_{{\rm D},\iota,k}$.
Since $x_{r,k}=\fbox{$k$}_{r-P(k)} +\fbox{$k-1$}_{r-P(k)+1}+\cdots+\fbox{$1$}_{r-P(k)+k-1}=[1,\cdots,k-1,k]_{r-P(k)}
\in {\rm Tab}_{{\rm D},\iota,k}$ $(r\in\mathbb{Z}_{\geq1})$, the inclusion $\Xi_{\iota,k}\subset {\rm Tab}_{{\rm D},\iota,k}$
follows by Proposition \ref{closednessD}.

Next, we prove ${\rm Tab}_{{\rm D},\iota,k}\subset \Xi_{\iota,k}$.
For $T\in {\rm Tab}_{{\rm D},\iota,k}$, we will show $T\in \Xi_{\iota,k}$
by induction on the weight of $T$. 
The elements which have the highest weight are
 $T=[1,2,\cdots,k]_s$ $(s\geq 1-P(k))$ by (\ref{wt-rel2}).
In this case, we have $T=x_{s+P(k),k}\in \Xi_{\iota,k}$.

Now we assume $T=[j_1,j_2,\cdots,j_k]_s\neq [1,2,\cdots,k]_s$ $(s\geq 1-P(k))$.
If $j_k\leq n-1$ then we can verify $T\in \Xi_{\iota,k}$
by the same way as in the proof of type A case. 
Hence, we suppose $j_k\in \{n,\ovl{n},\ovl{n-1},\cdots,\ovl{1}\}$.

\vspace{2mm}

\nd
\underline{Case 1. $j_k\in\{n,\ovl{n}\}$}

\vspace{2mm}

We consider the case $j_k\in\{n,\ovl{n}\}$.
One can take some integer $M\in[1,k]$
such that $1\leq j_1<\cdots<j_{k-M}\leq n-1$ and $j_{k-M+1}$, $j_{k-M+2}$, $\cdots$, $j_k\in\{n,\ovl{n}\}$.
By $j_{k-M+1}\ngeq j_{k-M+2}\ngeq\cdots\ngeq j_k$, the sequence $j_{k-M+1}$, $j_{k-M+2}$, $\cdots$, $j_k$
is an alternative sequence of $n$ and $\ovl{n}$. We consider the following two cases:

\vspace{2mm}

\nd
\underline{Case 1-1. $j_{k-M}=k-M$ or $k=M$}

\vspace{2mm}

In the case $k>M$, we have $j_1=1$, $j_2=2$, $\cdots$, $j_{k-M}=k-M$ by  $1\leq j_1<\cdots<j_{k-M}=k-M$,
in the case $k=M$, we see that $j_{1}$, $j_{2}$, $\cdots$, $j_k$
is an alternative sequence of $n$ and $\ovl{n}$.
By $k\leq n-2$, we can check $k-M\leq n-3$.
If $j_{k-M+1}=n$ (resp. $j_{k-M+1}=\ovl{n}$) then 
taking the element $T_n\in {\rm Tab}_{{\rm D},\iota,k}$ obtained from $T$
by replacing $j_{k-M+1}=n$ (resp. $j_{k-M+1}=\ovl{n}$) with $n-1$,
we get $S_{M_n,n-1}T_n=T$ (resp. $S_{M_n,n}T_n=T$) with some integer $M_n$ by Lemma \ref{Dlem-11}, \ref{Dlem-111}.
Using ${\rm wt}(T_n)>{\rm wt}(T)$ and the induction assumption, we get $T_n\in \Xi_{\iota,k}$.
Combining this result with 
$S_{M_n,n-1}T_n=T$ (resp. $S_{M_n,n}T_n=T$),
we obtain $T\in \Xi_{\iota,k}$.

\vspace{2mm}

\nd
\underline{Case 1-2. $j_{k-M}>k-M>0$}

\vspace{2mm}

There is an integer $r\in[1,k-M]$ such that
$j_{r-1}+1<j_r$, where we set $j_0=0$.
Let $T_r\in {\rm Tab}_{{\rm D},\iota,k}$ be the element 
obtained from $T$ by replacing $j_r$ with $j_r-1$.
It follows from Lemma \ref{Dlem-11}
that $S_{M_r,j_r-1}T_r=T$ with some integer $M_r$.
By the induction assumption, we get 
$T_r\in \Xi_{\iota,k}$. Thus, we also obtain
$T\in \Xi_{\iota,k}$.

Therefore, we may suppose that $j_{l-1}< \ovl{n-1}\leq j_{l}$
with some integer $l$ $(1\leq l\leq k)$. We can write
$j_l=\overline{i}$ with some $i\in I$.

\vspace{2mm}

\nd
\underline{Case 2. $j_l=\ovl{i}>\ovl{n-2}$}

\vspace{2mm}

Let $T_0\in {\rm Tab}_{{\rm D},\iota,k}$ be the element obtained from $T$
by replacing $j_l=\overline{i}$ with $\ovl{i+1}$. 
We can verify ${\rm wt}(T_0)>{\rm wt}(T)$ by (\ref{wt-rel2}).

\vspace{2mm}

\nd
\underline{Case 2-1. $T_0$ has neither $i$-cancelling nor $i$-double pair}

\vspace{2mm}

In this case, we have $S_{M,i}T_0=T$, where $M$ is an integer in Lemma \ref{Dlem-1}.
Since we know ${\rm wt}(T_0)>{\rm wt}(T)$, using the induction assumption, we obtain
$T_0=S_{l_p}\cdots S_{l_1}x_{s,k}$ with some $p\in\mathbb{Z}_{\geq0},\ s,\ l_1,\cdots,l_p\in \mathbb{Z}_{\geq1}$,
which yields $T=S_{M,i}S_{l_p}\cdots S_{l_1}x_{s,k}\in \Xi_{\iota,k}$.

\vspace{2mm}

\nd
\underline{Case 2-2. $T_0$ has an $i$-cancelling pair}

\vspace{2mm}

In this case, Lemma \ref{Dlem-1} means $j_{l-n+i+2}=i+1$ and $j_{l-n+i+1}\neq i$.
Instead of $T_0$, we consider the element $T'$ obtained from $T$ by
replacing  $j_{l-n+i+2}=i+1$ with $i$.
Using Lemma \ref{Dlem-11},
we see that the coefficient of $x_{M',i}$ ($M'$ is an integer of Lemma \ref{Dlem-11} (2)) in $T'$ is equal to $1$
and $S_{M',i}T'=T$. By (\ref{wt-rel2}), we get
${\rm wt}(T')>{\rm wt}(T)$.
It follows from the induction assumption, we obtain
$T'=S_{l_p}\cdots S_{l_1}x_{s,k}$ and $T=S_{M,i}S_{l_p}\cdots S_{l_1}x_{s,k}\in \Xi_{\iota,k}$
with some $p\in\mathbb{Z}_{\geq0},\ s,\ l_1,\cdots,l_p\in \mathbb{Z}_{\geq1}$.

\vspace{2mm}

\nd
\underline{Case 2-3. $T_0$ has an $i$-double pair}

\vspace{2mm}

In this case, Lemma \ref{Dlem-1} means $j_{l-n+i+1}=i$
and $j_{l-n+i+2}\neq i+1$. 
Since we supposed $i<n-2$ so that $l-n+i+2<l$, 
putting $J:=j_{l-n+i+2}$,
we have $i+2\leq J< \ovl{n-1}$.
If $J\in[1,n]$ (resp. $J=\ovl{n}$)
we can take $T''$ by replacing
$j_{l-n+i+2}=J$ with $J-1$ (resp. $n-1$) from $T$.
Lemma \ref{Dlem-11}, \ref{Dlem-111}
show that the coefficient of $x_{M'',J-1}$ (resp. $x_{M'',n}$)  ($M''$ is an integer) 
in $T''$ is $1$ and $T=S_{M'',J-1}T''$ (resp. $T=S_{M'',n}T''$).
Hence, $T\in \Xi_{\iota,k}$ follows by ${\rm wt}(T'')>{\rm wt}(T)$ and the induction assumption.

\vspace{2mm}

\nd
\underline{Case 3. $j_l=\ovl{i}=\ovl{n-2}$}

\vspace{2mm}

Next, let us consider the case $\ovl{i}=\ovl{n-2}$.
We can take some integer $M\in[0,l-1]$
such that $1\leq j_1<\cdots<j_{l-M-1}\leq n-1$ and $j_{l-M}$, $j_{l-M+1}$, $\cdots$, $j_{l-1}\in\{n,\ovl{n}\}$.

\vspace{2mm}

\nd
\underline{Case 3-1. $j_{l-M-1}\neq n-2$}

\vspace{2mm}

Let $T_3\in {\rm Tab}_{{\rm D},\iota,k}$ be the element obtained from $T$
by replacing $j_l=\overline{n-2}$ with $\ovl{n-1}$. 
We can verify ${\rm wt}(T_3)>{\rm wt}(T)$ by (\ref{wt-rel2}).
Lemma \ref{Dlem-1} says that there exists some integer $M_3\in\mathbb{Z}_{\geq1}$
such that $S_{M_3,n-2}T_3=T$. Using the induction assumption,
it follows $T_3\in \Xi_{\iota,k}$ and $T\in \Xi_{\iota,k}$.

\vspace{2mm}

\nd
\underline{Case 3-2. $j_{l-M-1}=n-2$}

\vspace{2mm}

If $M\geq1$ then we can prove $T\in \Xi_{\iota,k}$
just as in Case 3-1. 
If $M=0$ then $(j_{l-1},j_l)=(n-2,\ovl{n-2})$.
Using Definition \ref{box-def} (iv), one get
\begin{equation}\label{D-pr1-1}
\fbox{$n-2$}_{t+1}+\fbox{$\ovl{n-2}$}_{t}=\fbox{$n-1$}_{t+1}+\fbox{$\ovl{n-1}$}_{t}
\end{equation}
for any $t\in \mathbb{Z}_{\geq1-P(k)}$.
Thus, our proof is reduced to the following Case 4.

\vspace{2mm}

\nd
\underline{Case 4. $j_l=\ovl{i}=\ovl{n-1}$}

\vspace{2mm}

Just as in Case 3.,
we can take an integer $M\in[0,l-1]$
such that $1\leq j_1<\cdots<j_{l-M-1}\leq n-1$ and $j_{l-M}$, $j_{l-M+1}$, $\cdots$, $j_{l-1}\in\{n,\ovl{n}\}$.
If $j_{l-1}=n$ (resp. $j_{l-1}\neq n$) then we consider the element
$T_4\in {\rm Tab}_{{\rm D},\iota,k}$ obtained from $T$ replacing $j_l=\ovl{n-1}$
with $\ovl{n}$ (resp. $n$). By ${\rm wt}(T_4)>{\rm wt}(T)$,
the induction assumption means $T_4\in \Xi_{\iota,k}$.
By Lemma \ref{Dlem-1} and \ref{Dlem-111},
we get $S_{M_4,n-1}T_4=T$ (resp. $S_{M_4,n}T_4=T$) with some integer $M_4\in\mathbb{Z}_{\geq1}$.
Therefore, it follows $T\in \Xi_{\iota,k}$. \qed

\begin{lem}\label{Dlem-3}
In the case $\mathfrak{g}$ is of type ${\rm D}$, we have the following:
\begin{enumerate}
\item
\begin{equation}\label{dnn}
\{S_{l_p}\cdots S_{l_1}x_{r,n-1} |p\in\mathbb{Z}_{\geq0},\ r,\ l_1,\cdots,l_p\in \mathbb{Z}_{\geq1} \}
=\{ [j_1,j_2,\cdots,j_t]_s\in {\rm Tab}_{{\rm D},\iota} | j_1=\overline{n+1},\ t\ {\rm is\ even.} \}.
\end{equation}
\item
\begin{equation}\label{dn}
\{S_{l_p}\cdots S_{l_1}x_{r,n} |p\in\mathbb{Z}_{\geq0},\ r,\ l_1,\cdots,l_p\in \mathbb{Z}_{\geq1} \}
=\{ [j_1,j_2,\cdots,j_t]_s\in {\rm Tab}_{{\rm D},\iota} | j_1=\overline{n+1},\ t\ {\rm is\ odd.} \}.
\end{equation}
\end{enumerate}
\end{lem}

\nd
{\it Proof.}

\nd
(i) Let $\Xi_{\iota,n-1}$ (resp. $\Xi_{\iota,n}$) be the set of the left-hand side of (\ref{dnn}) (resp. (\ref{dn}))
and
${\rm Tab}_{{\rm D},\iota,n-1}$ (resp. ${\rm Tab}_{{\rm D},\iota,n}$) be the set of the right-hand side of (\ref{dnn}) (resp. (\ref{dn})).
First, let us prove $\Xi_{\iota,n-1}\subset{\rm Tab}_{{\rm D},\iota,n-1}$.
Since it is clear that $x_{r,n-1}=[\ovl{n+1}]_{r-P(n-1)+1}+[\ovl{n}]_{r-P(n-1)}\in{\rm Tab}_{{\rm D},\iota,n-1}$
for any $r\in\mathbb{Z}_{\geq1}$,
the inclusion $\Xi_{\iota,n-1}\subset{\rm Tab}_{{\rm D},\iota,n-1}$ follows by Proposition \ref{closednessD}.

Next, we turn to the proof of ${\rm Tab}_{{\rm D},\iota,n-1}\subset\Xi_{\iota,n-1}$.
By using Definition \ref{box-def}, we get
\begin{eqnarray}
\fbox{$\ovl{n+1}$}_{l+1}+\fbox{$\ovl{n}$}_{l}
&=&x_{l+P(n)+1,n}
+x_{l+P(n-1),n-1}-x_{l+P(n)+1,n}\nonumber\\
&=&x_{l+P(n-1),n-1},\label{D-pr1}
\end{eqnarray}
for $l\in\mathbb{Z}_{\geq 1-P(n-1)}$.
For each $T=[j_1,j_2,\cdots,j_t]_s\in{\rm Tab}_{{\rm D},\iota,n-1}$ with $j_1=\ovl{n+1}$, we prove $T\in \Xi_{\iota,n-1}$
by using the induction on ${\rm wt}(T)$.
The calculation (\ref{D-pr1}) means $[\ovl{n+1},\ovl{n}]_s=x_{s+P(n-1),n-1}\in\Xi_{\iota,n-1}$ for $s\geq 1-P(n-1)$,
so we may assume ${\rm wt}(T)<{\rm wt}([\ovl{n+1},\ovl{n}]_s)$.
If $j_2=\ovl{n}$ and $j_3> \ovl{n-1}$ then taking the element $T'=[\ovl{n+1},\ovl{n},\ovl{|j_3|+1},j_4,\cdots,j_t]_s$,
we get $S_{M',|j_3|}T'=T$ with some integer $M'$ by Lemma \ref{Dlem-1}.
Because of ${\rm wt}(T')>{\rm wt}(T)$,
and the induction assumption, we see that $T'\in \Xi_{\iota,n-1}$.
Since we know $S_{M',|j_3|}T'=T$, it follows $T\in \Xi_{\iota,n-1}$.
If $j_2=\ovl{n}$ and $j_3= \ovl{n-1}$ then taking the element $T''=[\ovl{n+1},j_4,\cdots,j_t]_s$,
we get $S_{M'',n}T''=S_{M'',n}[\ovl{n+1},j_4,\cdots,j_t]_s
=[\ovl{n+1},\ovl{n},\ovl{n-1},j_4,\cdots,j_t]_s=T$ with some integer $M''$
by $\fbox{$\ovl{n+1}$}_{s+t-1}=x_{s+t-1+P(n),n}$ and Lemma \ref{box-lem} (\ref{D-box5}).
By ${\rm wt}(T'')>{\rm wt}(T)$
and the induction assumption, we also get $T\in \Xi_{\iota,n-1}$.
If $j_2> \ovl{n}$ then taking the element $T^{\dagger}$
obtained from $T$ by replacing $j_2$ with $\ovl{|j_2|+1}$,
we have 
$S_{M^{\dagger},|j_2|}T^{\dagger}=T$ with some integer $M^{\dagger}$ by Lemma \ref{Dlem-1}.
It follows from the induction assumption that
$T^{\dagger}\in \Xi_{\iota,n-1}$ and $T\in \Xi_{\iota,n-1}$.
Therefore, we completed the proof of ${\rm Tab}_{{\rm D},\iota,n-1}\subset \Xi_{\iota,n-1}$,
which yields ${\rm Tab}_{{\rm D},\iota,n-1}= \Xi_{\iota,n-1}$.

We can prove (ii) by a similar way to (i). \qed

\vspace{2mm}

\nd
{\it Proof of Theorem \ref{thm2} for type D case}

\vspace{2mm}

Using the notation in the proofs of Lemma \ref{Dlem-2}, \ref{Dlem-3},
we obtain $\Xi_{\iota}=\coprod_{i\in I} \Xi_{\iota,i}$ and ${\rm Tab}_{{\rm D},\iota}=\coprod_{i\in I}{\rm Tab}_{{\rm D},\iota,i}$.
Hence,
the claim in Theorem \ref{thm2} is an easy consequence of the same lemmas. \qed

\vspace{2mm}

\nd
{\it Proof of Theorem \ref{thm1} for type D case}

\vspace{2mm}

\nd
We take an element $T=[j_1,\cdots,j_k]^{\rm D}_s\in {\rm Tab}_{{\rm D},\iota}$.

\nd
\underline{(I) The case $k\in[1,n-2]$ and $j_1\neq \ovl{n+1}$}

\vspace{2mm}

We suppose $k\in[1,n-2]$ and $j_1\neq \ovl{n+1}$.
Recall that 
$T=[j_1,\cdots,j_k]^{\rm D}_s=\sum_{i=1}^{k}\fbox{$j_i$}^{\rm D}_{s+k-i}$. 
Considering (\ref{D-box}), (\ref{D-rev1}) and (\ref{D-rev2}),
we see that if $x_{l,j}$ is a summand of $T$ with a negative coefficient then $2\leq l$.

\vspace{2mm}

\nd
\underline{(II) The case $k\in[1,n+1]$ and $j_1= \ovl{n+1}$}

\vspace{2mm}

Recall that the following (\ref{D-box-ag}):
\[
\fbox{$j_i$}^{\rm D}_{s+k-i}=
\begin{cases}
x_{s+k-i+P(n),n} & {\rm if}\ j_i=\ovl{n+1}\\
x_{s+k-i+P(n-1),n-1}-x_{s+k-i+P(n)+1,n} & {\rm if}\ j_i= \ovl{n}, \\ 
x_{s+k-i+P(n-2)+1,n-2}-x_{s+k-i+P(n-1)+1,n-1}
-x_{s+k-i+P(n)+1,n} & {\rm if}\ j_i= \ovl{n-1},\\
x_{s+k-i+P(|j_i|-1)+n-|j_i|,|j_i|-1}-x_{s+k-i+P(|j_i|)+n-|j_i|,|j_i|} & {\rm if}\ j_i\geq \ovl{n-2}.
\end{cases}
\]
By the argument in the beginning of proof of Proposition \ref{closednessD}(iii),
if $x_{l,j}$ is a summand of $T$ with a negative coefficient then $2\leq l$. \qed

\subsection{Proof of Corollary \ref{cor1}}

As in Sect. \ref{Stabsect}, we denote each tableau $\begin{ytableau}
j_1 \\
\vdots \\
j_k
\end{ytableau}^{\rm X}_{s}$ by $[j_1,\cdots,j_k]_s^{\rm X}$.

\vspace{2mm}

\nd
{\it Proof of Corollary \ref{cor1} for type A case}

\vspace{2mm}

\nd
For $k\in[1,n]$, $s\in\mathbb{Z}_{\geq1-P(k)}$ and $t\in\mathbb{Z}_{\geq1-P(n-k+1)}$,
there exist elements $[1,2,\cdots,k]^{\rm A}_s$, $[k+1,\cdots,n,n+1]^{\rm A}_t\in {\rm Tab}_{{\rm A},\iota}$ and
\begin{equation}\label{cor-pr1}
[1,2,\cdots,k]^{\rm A}_s=x_{s+P(k),k},
\end{equation}
\begin{equation}\label{cor-pr1a}
 [k+1,\cdots,n,n+1]^{\rm A}_t=-x_{t+n-k+P(k)+1,k}.
\end{equation}
By
Theorem \ref{polyhthm}, \ref{thm2} and \ref{thm1},
we have
\[
{\rm Im}(\Psi_{\iota}) =\Sigma_{\io}=\{\textbf{x}\in\mathbb{Z}^{\infty}_{\iota} | \varphi(\textbf{x})\geq0
\ {\rm for\ all}\ \varphi\in {\rm Tab}_{{\rm A},\iota}\}.
\]
For $\textbf{x}=(x_{l,j})_{l\in\mathbb{Z}_{\geq1}, j\in I}\in {\rm Im}(\Psi_{\iota})$,
using (\ref{cor-pr1}), we get $x_{l,k}\geq0$ for any $l\in\mathbb{Z}_{\geq1}$.
By $P(k)<k$, it holds $t+n-k+P(k)+1\leq t+n$.
Since $t\geq1-P(n-k+1)$ so that $t\geq1$,
taking (\ref{cor-pr1a}) into account, we see $x_{l',k}\leq0$ for $l'\in\mathbb{Z}_{\geq n+1}$.
Thus, $x_{l,k}=0$ for $l\in\mathbb{Z}$ with $l>n$ and $k\in I$.
It follows from Definition \ref{box-def} (i), \ref{tab-def} (i) that
$[j_1,\cdots,j_k]^{\rm A}_s \in {\rm Tab}_{{\rm A},\iota}$ $(s>n)$ is a $\mathbb{Z}$-linear combination of
$x_{l,j}$ $(l>n,\ j\in I)$.

\vspace{2mm}

\nd
{\it Proof of Corollary \ref{cor1} for type B case}

\vspace{2mm}

\nd
For $k\in[1,n]$ and $s\in\mathbb{Z}$ with $s\geq 1-P(k)$,
there exist elements $[1,2,\cdots,k]^{\rm B}_s$, $[\ovl{k},\cdots,\ovl{2},\ovl{1}]^{\rm B}_s\in {\rm Tab}_{{\rm B},\iota}$ and
\begin{equation}\label{cor-pr2}
[1,2,\cdots,k]^{\rm B}_s=x_{s+P(k),k},\ \ [\ovl{k},\cdots,\ovl{2},\ovl{1}]^{\rm B}_s=-x_{s+P(k)+n,k}.
\end{equation}
By
Theorem \ref{polyhthm}, \ref{thm2}, \ref{thm1}, (\ref{cor-pr2}) and $s\geq 1-P(k)$,
if $\textbf{x}=(x_{l,j})_{l\in\mathbb{Z}_{\geq1}, j\in I}\in {\rm Im}(\Psi_{\iota})$ then
 we get
$x_{l,k}=0$ for $l\in\mathbb{Z}$ with $l>n$ and $k\in I$.
It follows from Definition \ref{box-def} (ii), \ref{tab-def} (i) that
$[j_1,\cdots,j_k]^{\rm B}_s \in {\rm Tab}_{{\rm B},\iota}$ $(s>n)$ is a $\mathbb{Z}$-linear combination of
$x_{l,j}$ $(l>n,\ j\in I)$.

\vspace{2mm}

\nd
{\it Proof of Corollary \ref{cor1} for type C case}

\vspace{2mm}

\nd
For $k\in[1,n-1]$ and $s\in\mathbb{Z}$ with $s\geq 1-P(k)$,
there exist elements $[1,2,\cdots,k]^{\rm C}_s$, $[\ovl{k},\cdots,\ovl{2},\ovl{1}]^{\rm C}_s\in {\rm Tab}_{{\rm C},\iota}$ and
\[
[1,2,\cdots,k]^{\rm C}_s=x_{s+P(k),k},\ \ [\ovl{k},\cdots,\ovl{2},\ovl{1}]^{\rm C}_s=-x_{s+P(k)+n,k}.
\]
As in the type B case, we see that for $\textbf{x}=(x_{l,j})_{l\in\mathbb{Z}_{\geq1}, j\in I}\in {\rm Im}(\Psi_{\iota})$,
 we get
$x_{l,k}=0$ for $l\in\mathbb{Z}$ with $l>n$ and $k\in [1,n-1]$. We also see that for $s\in\mathbb{Z}$ with $s\geq 1-P(n)$,
\begin{equation}\label{cor-pr3}
\fbox{$\overline{n+1}$}^{\rm C}_{s}
:=x_{s+P(n),n}\in {\rm Tab}_{{\rm C},\iota},\qquad [\ovl{n+1},\ovl{n},\cdots,\ovl{1}]_s=-x_{s+n+P(n),n}\in {\rm Tab}_{{\rm C},\iota}.
\end{equation}
Thus, 
if $\textbf{x}=(x_{l,j})_{l\in\mathbb{Z}, j\in I}\in {\rm Im}(\Psi_{\iota})$ then
$x_{l,n}=0$ for $l\in\mathbb{Z}$ with $l>n$.
Taking \ref{box-def} (iii), \ref{tab-def} (i) into account that
$[j_1,\cdots,j_k]^{\rm C}_s \in {\rm Tab}_{{\rm C},\iota}$ $(s>n)$ is a $\mathbb{Z}$-linear combination of
$x_{l,j}$ $(l>n,\ j\in I)$.

\vspace{2mm}

\nd
{\it Proof of Corollary \ref{cor1} for type D case}

\vspace{2mm}

\nd
For $k\in[1,n-2]$ and $s\in\mathbb{Z}$ with $s\geq 1-P(k)$, we see that
\[
[1,2,\cdots,k]^{\rm D}_s=x_{s+P(k),k}\in {\rm Tab}_{{\rm D},\iota},
\ \ [\ovl{k},\cdots,\ovl{2},\ovl{1}]^{\rm D}_s=-x_{s+P(k)+n-1,k}\in {\rm Tab}_{{\rm D},\iota},
\]
and for $s\in\mathbb{Z}$ with $s\geq 1-P(n-2)$ so that $s\geq 1-P(n-1)$ and $s\geq 1-P(n)$, it follows
\[
[\ovl{n+1}]^{\rm D}_s=x_{s+P(n),n}\in {\rm Tab}_{{\rm D},\iota},\quad 
[\ovl{n+1},\ovl{n}]^{\rm D}_s=x_{s+P(n-1),n-1}\in {\rm Tab}_{{\rm D},\iota},\quad 
\]
\[
[\ovl{n+1},\ovl{n},\cdots,\ovl{1}]^{\rm D}_s=-x_{s+n+P(n)-1,n}\in {\rm Tab}_{{\rm D},\iota},\quad 
[\ovl{n+1},\ovl{n-1},\ovl{n-2},\cdots,\ovl{1}]^{\rm D}_s=-x_{s+n+P(n-1)-1,n-1}\in {\rm Tab}_{{\rm D},\iota}. 
\]
Thus, 
if $\textbf{x}=(x_{l,j})_{l\in\mathbb{Z}, j\in I}\in {\rm Im}(\Psi_{\iota})$ then
$x_{l,n}=0$ for $l\in\mathbb{Z}$ with $l>n$.
Taking \ref{box-def} (iv), \ref{tab-def} (i) into account that
$[j_1,\cdots,j_k]^{\rm D}_s \in {\rm Tab}_{{\rm D},\iota}$ $(s>n)$ is a $\mathbb{Z}$-linear combination of
$x_{l,j}$ $(l>n,\ j\in I)$. \qed


\begin{thebibliography}{9}







\bibitem{HL}D.Hernandez, B.Leclerc,
Quantum Grothendieck rings and derived Hall algebras,
J. Reine Angew. Math. 701, 77–-126 (2015). 

\bibitem{H1}
A.Hoshino,
Polyhedral realizations of crystal bases for quantum algebras of finite types,
J. Math. Phys. 46, no. 11, 113514,  31 pp, (2005). 


\bibitem{H2}
A.Hoshino, 
Polyhedral realizations of crystal bases for quantum algebras of classical affine types,
J. Math. Phys. 54, no. 5, 053511, 28 pp  (2013). 


\bibitem{Kac} V. G. Kac,
{\it Infinite-dimensional Lie algebras}, 
third edition. Cambridge University Press, Cambridge, xxii+400 pp, (1990).  

\bibitem{K0} M.Kashiwara, Crystalling the $q$-analogue of universal 
              enveloping algebras, Comm. Math. Phys.,
            {\it 133}, 249--260 (1990).

\bibitem{K1} M.Kashiwara,
 On crystal bases of the $q$-analogue of universal enveloping algebras,
	Duke Math. J., {\it 63} (2), 465--516 (1991).

\bibitem{K3}M.Kashiwara, 
The crystal base and Littelmann's refined Demazure character formula,
Duke Math. J., 71, no 3,  839--858 (1993).

\bibitem{L}
G.Lusztig, Canonical bases arising from quantized enveloping algebras, 
J. Amer. Math. Soc. 3, no. 2, 447--498 (1990). 



\bibitem{Nj}H. Nakajima, $t$-analogs of $q$-characters of quantum affine algebras of 
type $A_n$, $D_n$, Combinatorial and geometric representation theory (Seoul, 2001),
Contemp. Math, 325, AMS, Providence, RI, 141--160 (2003).

\bibitem{NN}
W.Nakai, T.Nakanishi,
Paths, tableaux and $q$-characters of quantum affine algebras: The $C_n$ case, 
J. Phys. A: Math. Gen. 39, no. 9, 2083–-2115, (2006). 


\bibitem{N99}T.Nakashima, Polyhedral realizations of crystal bases for integrable highest weight modules,
J. Algebra, vol.219, no. 2, 571–-597, (1999). 

\bibitem{NZ}
T.Nakashima,  A.Zelevinsky, Polyhedral realizations of 
crystal bases for quantized Kac-Moody algebras,
Adv. Math. {\bf 131}, no. 1, 253--278, (1997). 

\bibitem{R}
R.B\'{e}dard, On commutation classes of reduced words in Weyl groups,
European J. Combin. 20, no. 6, 483–-505, (1999).


\end{thebibliography}
\end{document}